\documentclass[PhD]{dukethesis2}
\usepackage{amsmath,amsfonts,amsthm,graphicx,maple2e}
\usepackage{amssymb}
\usepackage{amsmath,amsfonts}

\def\flecha{\longrightarrow}
\def\implies{\Longrightarrow}
\def\maparrow#1{\stackrel{#1}{\longrightarrow}}

\def\mtx#1#2#3#4{\left(\begin{array}{cc} #1 & #2 \\ #3 & #4 \end{array} \right)}
\def\row#1#2{\left(\begin{array}{cc} #1 & #2 \end{array} \right)}

\def\col#1#2{\left(\begin{array}{c} #1 \\ #2 \end{array} \right)}

\def\bmtxtth#1#2#3#4#5#6{\left[\begin{array}{ccc} #1 & #2 & #3
\\ #4 & #5 & #6 \end{array} \right]} 
\def\bMtx#1#2#3#4#5#6#7#8#9{\left[\begin{array}{ccc} #1 & #2 & #3 \\
#4 & #5 & #6 \\ #7 & #8 & #9 \end{array} \right]} 

\def\bRow#1#2#3{\left[\begin{array}{ccc} #1 & #2 & #3 \end{array}
\right]}

\def\Bcol#1#2#3{\left[\begin{array}{c} #1 \\ #2 \\ #3 \end{array} \right]}

\def\Tr{\operatorname{Tr}}
\def\deg{\operatorname{deg}}
\def\Spec{\operatorname{Spec}}

\def\OO{\mathcal{O}}
\def\MM{\mathcal{M}}
\def\FF{\mathcal{F}}

\def\WW{\mathcal{W}}
\def\N{\mathcal{N}}
\def\M{\mathcal{M}}

\def\CC{\mathbb{C}}

\def\PP{\mathbb{P}}

\def\ZZ{\mathbb{Z}}
\def\BB{\mathbb{B}}

\def\for{\mathrm{for}}
\def\and{\mathrm{and}}
\def\any{\mathrm{any}}
\def\where{\mathrm{where}}

\def\st{\mathrm{s.t.}}
\def\even{\mathrm{even}}
\def\odd{\mathrm{odd}}

\def\rhotil{\tilde{\rho}}

\def\skp{\\[.2in]}
\def\skop{\;\;\;\;\;\;\;\;}

\def\til{\widetilde}
\def\Hat{\widehat}

\def\part#1{\dfrac{\partial}{\partial #1}}
\def\d{\mathrm{d}}

\def\12{\dfrac{1}{2}}

\def\C{\mathcal{C}}
\def\S{\mathcal{S}}

\makeatletter

\def\diagram{\leftwidth=\z@ \rightwidth=\z@ \topheight=\z@
\botheight=\z@ \setbox\@picbox\hbox\bgroup}

\def\enddiagram{\egroup\wd\@picbox\rightwidth\unitlength
\ht\@picbox\topheight\unitlength \dp\@picbox\botheight\unitlength
\hskip\leftwidth\unitlength\box\@picbox}

\def\bfig{\begin{diagram}}
\def\efig{\end{diagram}}
\newcount\wideness \newcount\leftwidth \newcount\rightwidth
\newcount\highness \newcount\topheight \newcount\botheight

\def\ratchet#1#2{\ifnum#1<#2 \global #1=#2 \fi}

\def\putbox(#1,#2)#3{%
\horsize{\wideness}{#3} \divide\wideness by 2
{\advance\wideness by #1 \ratchet{\rightwidth}{\wideness}}
{\advance\wideness by -#1 \ratchet{\leftwidth}{\wideness}}
\vertsize{\highness}{#3} \divide\highness by 2
{\advance\highness by #2 \ratchet{\topheight}{\highness}}
{\advance\highness by -#2 \ratchet{\botheight}{\highness}}
\put(#1,#2){\makebox(0,0){$#3$}}}

\def\putlbox(#1,#2)#3{%
\horsize{\wideness}{#3}
{\advance\wideness by #1 \ratchet{\rightwidth}{\wideness}}
{\ratchet{\leftwidth}{-#1}}
\vertsize{\highness}{#3} \divide\highness by 2
{\advance\highness by #2 \ratchet{\topheight}{\highness}}
{\advance\highness by -#2 \ratchet{\botheight}{\highness}}
\put(#1,#2){\makebox(0,0)[l]{$#3$}}}

\def\putrbox(#1,#2)#3{%
\horsize{\wideness}{#3}
{\ratchet{\rightwidth}{#1}}
{\advance\wideness by -#1 \ratchet{\leftwidth}{\wideness}}
\vertsize{\highness}{#3} \divide\highness by 2
{\advance\highness by #2 \ratchet{\topheight}{\highness}}
{\advance\highness by -#2 \ratchet{\botheight}{\highness}}
\put(#1,#2){\makebox(0,0)[r]{$#3$}}}

\def\adjust[#1]{} 

\newcount \coefa
\newcount \coefb
\newcount \coefc
\newcount\tempcounta
\newcount\tempcountb
\newcount\tempcountc
\newcount\tempcountd
\newcount\xext
\newcount\yext
\newcount\xoff
\newcount\yoff
\newcount\gap%
\newcount\arrowtypea
\newcount\arrowtypeb
\newcount\arrowtypec
\newcount\arrowtyped
\newcount\arrowtypee
\newcount\height
\newcount\width
\newcount\xpos
\newcount\ypos
\newcount\run
\newcount\rise
\newcount\arrowlength
\newcount\halflength
\newcount\arrowtype
\newdimen\tempdimen
\newdimen\xlen
\newdimen\ylen
\newsavebox{\tempboxa}%
\newsavebox{\tempboxb}%
\newsavebox{\tempboxc}%

\newdimen\w@dth

\def\setw@dth#1#2{\setbox\z@\hbox{$#1$}\w@dth=\wd\z@
\setbox\@ne\hbox{$#2$}\ifnum\w@dth<\wd\@ne \w@dth=\wd\@ne \fi
\advance\w@dth by 1.2em}

\def\t@^#1_#2{\def\n@one{#1}\def\n@two{#2}\mathrel{\setw@dth{#1}{#2}
\mathop{\hbox to \w@dth{\rightarrowfill}}\limits
\ifx\n@one\empty\else ^{\box\z@}\fi
\ifx\n@two\empty\else _{\box\@ne}\fi}}
\def\t@@^#1{\@ifnextchar_ {\t@^{#1}}{\t@^{#1}_{}}}
\def\to{\@ifnextchar^ {\t@@}{\t@@^{}}}

\def\t@left^#1_#2{\def\n@one{#1}\def\n@two{#2}\mathrel{\setw@dth{#1}{#2}
\mathop{\hbox to \w@dth{\leftarrowfill}}\limits
\ifx\n@one\empty\else ^{\box\z@}\fi
\ifx\n@two\empty\else _{\box\@ne}\fi}}
\def\t@@left^#1{\@ifnextchar_ {\t@left^{#1}}{\t@left^{#1}_{}}}
\def\toleft{\@ifnextchar^ {\t@@left}{\t@@left^{}}}

\def\two@^#1_#2{\def\n@one{#1}\def\n@two{#2}\mathrel{\setw@dth{#1}{#2}
\mathop{\vcenter{\hbox to \w@dth{\rightarrowfill}\kern-1.7ex
         \hbox to \w@dth{\rightarrowfill}}%
       }\limits
\ifx\n@one\empty\else ^{\box\z@}\fi
\ifx\n@two\empty\else _{\box\@ne}\fi}}
\def\tw@@^#1{\@ifnextchar_ {\two@^{#1}}{\two@^{#1}_{}}}
\def\two{\@ifnextchar^ {\tw@@}{\tw@@^{}}}

\def\tofr@^#1_#2{\def\n@one{#1}\def\n@two{#2}\mathrel{\setw@dth{#1}{#2}
\mathop{\vcenter{\hbox to \w@dth{\rightarrowfill}\kern-1.7ex
         \hbox to \w@dth{\leftarrowfill}}%
       }\limits
\ifx\n@one\empty\else ^{\box\z@}\fi
\ifx\n@two\empty\else _{\box\@ne}\fi}}
\def\t@fr@^#1{\@ifnextchar_ {\tofr@^{#1}}{\tofr@^{#1}_{}}}
\def\tofro{\@ifnextchar^ {\t@fr@}{\t@fr@^{}}}

\def\mon{\mathop{\m@th\hbox to
      14.6\P@{\lasyb\char'51\hskip-2.1\P@$\arrext$\hss
$\mathord\rightarrow$}}\limits} 
\def\leftmono{\mathrel{\m@th\hbox to
14.6\P@{$\mathord\leftarrow$\hss$\arrext$\hskip-2.1\P@\lasyb\char'50%
}}\limits} 
\mathchardef\arrext="0200       

\def\settypes(#1,#2,#3){\arrowtypea#1 \arrowtypeb#2 \arrowtypec#3}
\def\settoheight#1#2{\setbox\@tempboxa\hbox{#2}#1\ht\@tempboxa\relax}%
\def\settodepth#1#2{\setbox\@tempboxa\hbox{#2}#1\dp\@tempboxa\relax}%
\def\settokens[#1`#2`#3`#4]{%
     \def\tokena{#1}\def\tokenb{#2}\def\tokenc{#3}\def\tokend{#4}}
\def\setsqparms[#1`#2`#3`#4;#5`#6]{%
\arrowtypea #1
\arrowtypeb #2
\arrowtypec #3
\arrowtyped #4
\width #5
\height #6
}
\def\setpos(#1,#2){\xpos=#1 \ypos#2}

\def\settriparms[#1`#2`#3;#4]{\settripairparms[#1`#2`#3`1`1;#4]}%

\def\settripairparms[#1`#2`#3`#4`#5;#6]{%
\arrowtypea #1
\arrowtypeb #2
\arrowtypec #3
\arrowtyped #4
\arrowtypee #5
\width #6
\height #6
}

\def\resetparms{\settripairparms[1`1`1`1`1;500]\width 500}

\resetparms

\def\mvector(#1,#2)#3{
\put(0,0){\vector(#1,#2){#3}}%
\put(0,0){\vector(#1,#2){26}}%
}
\def\evector(#1,#2)#3{{
\arrowlength #3
\put(0,0){\vector(#1,#2){\arrowlength}}%
\advance \arrowlength by-30
\put(0,0){\vector(#1,#2){\arrowlength}}%
}}

\def\horsize#1#2{%
\settowidth{\tempdimen}{$#2$}%
#1=\tempdimen
\divide #1 by\unitlength
}

\def\vertsize#1#2{%
\settoheight{\tempdimen}{$#2$}%
#1=\tempdimen
\settodepth{\tempdimen}{$#2$}%
\advance #1 by\tempdimen
\divide #1 by\unitlength
}

\def\putvector(#1,#2)(#3,#4)#5#6{{%
\ifnum3<\arrowtype
\putdashvector(#1,#2)(#3,#4)#5\arrowtype
\else
\ifnum\arrowtype<-3
\putdashvector(#1,#2)(#3,#4)#5\arrowtype
\else
\xpos=#1
\ypos=#2
\run=#3
\rise=#4
\arrowlength=#5
\ifnum \arrowtype<0
    \ifnum \run=0
    \advance \ypos by-\arrowlength
    \else
    \tempcounta \arrowlength
    \multiply \tempcounta by\rise
    \divide \tempcounta by\run
    \ifnum\run>0
        \advance \xpos by\arrowlength
        \advance \ypos by\tempcounta
    \else
        \advance \xpos by-\arrowlength
        \advance \ypos by-\tempcounta
    \fi
    \fi
    \multiply \arrowtype by-1
    \multiply \rise by-1
    \multiply \run by-1
\fi
\ifcase \arrowtype
\or \put(\xpos,\ypos){\vector(\run,\rise){\arrowlength}}%
\or \put(\xpos,\ypos){\mvector(\run,\rise)\arrowlength}%
\or \put(\xpos,\ypos){\evector(\run,\rise){\arrowlength}}%
\fi\fi\fi
}}

\def\putsplitvector(#1,#2)#3#4{
\xpos #1
\ypos #2
\arrowtype #4
\halflength #3
\arrowlength #3
\gap 140
\advance \halflength by-\gap
\divide \halflength by2
\ifnum\arrowtype>0
   \ifcase \arrowtype
   \or \put(\xpos,\ypos){\line(0,-1){\halflength}}%
       \advance\ypos by-\halflength
       \advance\ypos by-\gap
       \put(\xpos,\ypos){\vector(0,-1){\halflength}}%
   \or \put(\xpos,\ypos){\line(0,-1)\halflength}%
       \put(\xpos,\ypos){\vector(0,-1)3}%
       \advance\ypos by-\halflength
       \advance\ypos by-\gap
       \put(\xpos,\ypos){\vector(0,-1){\halflength}}%
   \or \put(\xpos,\ypos){\line(0,-1)\halflength}%
       \advance\ypos by-\halflength
       \advance\ypos by-\gap
       \put(\xpos,\ypos){\evector(0,-1){\halflength}}%
   \fi
\else \arrowtype=-\arrowtype
   \ifcase\arrowtype
   \or \advance \ypos by-\arrowlength
       \put(\xpos,\ypos){\line(0,1){\halflength}}%
       \advance\ypos by\halflength
       \advance\ypos by\gap
       \put(\xpos,\ypos){\vector(0,1){\halflength}}%
   \or \advance \ypos by-\arrowlength
       \put(\xpos,\ypos){\line(0,1)\halflength}%
       \put(\xpos,\ypos){\vector(0,1)3}%
       \advance\ypos by\halflength
       \advance\ypos by\gap
       \put(\xpos,\ypos){\vector(0,1){\halflength}}%
   \or \advance \ypos by-\arrowlength
       \put(\xpos,\ypos){\line(0,1)\halflength}%
       \advance\ypos by\halflength
       \advance\ypos by\gap
       \put(\xpos,\ypos){\evector(0,1){\halflength}}%
   \fi
\fi
}

\def\putmorphism(#1)(#2,#3)[#4`#5`#6]#7#8#9{{%
\run #2
\rise #3
\ifnum\rise=0
  \puthmorphism(#1)[#4`#5`#6]{#7}{#8}#9%
\else\ifnum\run=0
  \putvmorphism(#1)[#4`#5`#6]{#7}{#8}#9%
\else
\setpos(#1)%
\arrowlength #7
\arrowtype #8
\ifnum\run=0
\else\ifnum\rise=0
\else
\ifnum\run>0
    \coefa=1
\else
   \coefa=-1
\fi
\ifnum\arrowtype>0
   \coefb=0
   \coefc=-1
\else
   \coefb=\coefa
   \coefc=1
   \arrowtype=-\arrowtype
\fi
\width=2
\multiply \width by\run
\divide \width by\rise
\ifnum \width<0  \width=-\width\fi
\advance\width by60
\if l#9 \width=-\width\fi
\putbox(\xpos,\ypos){#4}
{\multiply \coefa by\arrowlength
\advance\xpos by\coefa
\multiply \coefa by\rise
\divide \coefa by\run
\advance \ypos by\coefa
\putbox(\xpos,\ypos){#5} }%
{\multiply \coefa by\arrowlength
\divide \coefa by2
\advance \xpos by\coefa
\advance \xpos by\width
\multiply \coefa by\rise
\divide \coefa by\run
\advance \ypos by\coefa
\if l#9%
   \putrbox(\xpos,\ypos){#6}%
\else\if r#9%
   \putlbox(\xpos,\ypos){#6}%
\fi\fi }%
{\multiply \rise by-\coefc
\multiply \run by-\coefc
\multiply \coefb by\arrowlength
\advance \xpos by\coefb
\multiply \coefb by\rise
\divide \coefb by\run
\advance \ypos by\coefb
\multiply \coefc by70
\advance \ypos by\coefc
\multiply \coefc by\run
\divide \coefc by\rise
\advance \xpos by\coefc
\multiply \coefa by140
\multiply \coefa by\run
\divide \coefa by\rise
\advance \arrowlength by\coefa
\ifcase\arrowtype
\or \put(\xpos,\ypos){\vector(\run,\rise){\arrowlength}}%
\or \put(\xpos,\ypos){\mvector(\run,\rise){\arrowlength}}%
\or \put(\xpos,\ypos){\evector(\run,\rise){\arrowlength}}%
\fi}\fi\fi\fi\fi}}

\newcount\numbdashes \newcount\lengthdash \newcount\increment

\def\howmanydashes{
\numbdashes=\arrowlength \lengthdash=40
\divide\numbdashes by \lengthdash
\lengthdash=\arrowlength
\divide\lengthdash by \numbdashes
\increment=\lengthdash
\multiply\lengthdash by 3
\divide\lengthdash by 5
}

\def\putdashvector(#1)(#2,#3)#4#5{%
\ifnum#3=0 \putdashhvector(#1){#4}#5
\else
\ifnum#2=0
\putdashvvector(#1){#4}#5\fi\fi}

\def\putdashhvector(#1,#2)#3#4{{%
\arrowlength=#3 \howmanydashes
\multiput(#1,#2)(\increment,0){\numbdashes}%
{\vrule height .4pt width \lengthdash\unitlength}
\arrowtype=#4 \xpos=#1
\ifnum\arrowtype<0 \advance\arrowtype by 7 \fi
\ifcase\arrowtype
\or \advance\xpos by 10
    \put(\xpos,#2){\vector(-1,0){\lengthdash}}
    \advance\xpos by 40
    \put(\xpos,#2){\vector(-1,0){\lengthdash}}
\or \advance \xpos by 10
    \put(\xpos,#2){\vector(-1,0){\lengthdash}}
    \advance\xpos by  \arrowlength
    \advance\xpos by  -50
    \put(\xpos,#2){\vector(-1,0){\lengthdash}}
\or \advance\xpos by 10
    \put(\xpos,#2){\vector(-1,0){\lengthdash}}
\or \advance\xpos by \arrowlength
    \advance\xpos by -\lengthdash
    \put(\xpos,#2){\vector(1,0){\lengthdash}}
\or {\advance\xpos by 10
    \put(\xpos,#2){\vector(1,0){\lengthdash}}}
    \advance\xpos by \arrowlength
    \advance\xpos by -\lengthdash
    \put(\xpos,#2){\vector(1,0){\lengthdash}}
\or \advance\xpos by \arrowlength
    \advance\xpos by -\lengthdash
    \put(\xpos,#2){\vector(1,0){\lengthdash}}
    \advance\xpos by -40
    \put(\xpos,#2){\vector(1,0){\lengthdash}}
   \fi
}}

\def\putdashvvector(#1,#2)#3#4{{%
\arrowlength=#3 \howmanydashes
\ypos=#2 \advance\ypos by -\arrowlength
\multiput(#1,#2)(0,\increment){\numbdashes}%
    {\vrule width .4pt height \lengthdash\unitlength}
\arrowtype=#4 \ypos=#2
\ifnum\arrowtype<0 \advance\arrowtype by 7 \fi
\ifcase\arrowtype
\or \advance\ypos by \arrowlength \advance\ypos by -40
    \put(#1,\ypos){\vector(0,1){\lengthdash}}
    \advance\ypos by -40
    \put(#1,\ypos){\vector(0,1){\lengthdash}}
\or \advance\ypos by 10
    \put(#1,\ypos){\vector(0,1){\lengthdash}}
    \advance\ypos by \arrowlength \advance\ypos by -40
    \put(#1,\ypos){\vector(0,1){\lengthdash}}
\or \advance\ypos by \arrowlength \advance\ypos by -40
    \put(#1,\ypos){\vector(0,1){\lengthdash}}
\or \advance\ypos by 10
    \put(#1,\ypos){\vector(0,-1){\lengthdash}}
\or \advance\ypos by 10
    \put(#1,\ypos){\vector(0,-1){\lengthdash}}
    \advance\ypos by \arrowlength \advance\ypos by -40
    \put(#1,\ypos){\vector(0,-1){\lengthdash}}
\or \advance\ypos by 10
    \put(#1,\ypos){\vector(0,-1){\lengthdash}}
    \advance\ypos by 40
    \put(#1,\ypos){\vector(0,-1){\lengthdash}}
\fi
}}

\def\puthmorphism(#1,#2)[#3`#4`#5]#6#7#8{{%
\xpos #1
\ypos #2
\width #6
\arrowlength #6
\arrowtype=#7
\putbox(\xpos,\ypos){#3\vphantom{#4}}%
{\advance \xpos by\arrowlength
\putbox(\xpos,\ypos){\vphantom{#3}#4}}%
\horsize{\tempcounta}{#3}%
\horsize{\tempcountb}{#4}%
\divide \tempcounta by2
\divide \tempcountb by2
\advance \tempcounta by30
\advance \tempcountb by30
\advance \xpos by\tempcounta
\advance \arrowlength by-\tempcounta
\advance \arrowlength by-\tempcountb
\putvector(\xpos,\ypos)(1,0)\arrowlength\arrowtype
\divide \arrowlength by2
\advance \xpos by\arrowlength
\vertsize{\tempcounta}{#5}%
\divide\tempcounta by2
\advance \tempcounta by20
\if a#8 %
   \advance \ypos by\tempcounta
   \putbox(\xpos,\ypos){#5}%
\else
   \advance \ypos by-\tempcounta
   \putbox(\xpos,\ypos){#5}%
\fi}}

\def\putvmorphism(#1,#2)[#3`#4`#5]#6#7#8{{%
\xpos #1
\ypos #2
\arrowlength #6
\arrowtype #7
\settowidth{\xlen}{$#5$}%
\putbox(\xpos,\ypos){#3}%
{\advance \ypos by-\arrowlength
\putbox(\xpos,\ypos){#4}}%
{\advance\arrowlength by-140
\advance \ypos by-70
\ifdim\xlen>0pt
   \if m#8%
      \putsplitvector(\xpos,\ypos)\arrowlength\arrowtype
   \else
   \putvector(\xpos,\ypos)(0,-1)\arrowlength\arrowtype
   \fi
\else
   \putvector(\xpos,\ypos)(0,-1)\arrowlength\arrowtype
\fi}%
\ifdim\xlen>0pt
   \divide \arrowlength by2
   \advance\ypos by-\arrowlength
   \if l#8%
      \advance \xpos by-40
      \putrbox(\xpos,\ypos){#5}%
   \else\if r#8%
      \advance \xpos by40
      \putlbox(\xpos,\ypos){#5}%
   \else
      \putbox(\xpos,\ypos){#5}%
   \fi\fi
\fi
}}

\def\putsquarep<#1>(#2)[#3;#4`#5`#6`#7]{{%
\setsqparms[#1]%
\setpos(#2)%
\settokens[#3]%
\puthmorphism(\xpos,\ypos)[\tokenc`\tokend`{#7}]{\width}{\arrowtyped}b%
\advance\ypos by \height
\puthmorphism(\xpos,\ypos)[\tokena`\tokenb`{#4}]{\width}{\arrowtypea}a%
\putvmorphism(\xpos,\ypos)[``{#5}]{\height}{\arrowtypeb}l%
\advance\xpos by \width
\putvmorphism(\xpos,\ypos)[``{#6}]{\height}{\arrowtypec}r%
}}

\def\putsquare{\@ifnextchar <{\putsquarep}{\putsquarep%
   <\arrowtypea`\arrowtypeb`\arrowtypec`\arrowtyped;\width`\height>}}
\def\square{\@ifnextchar< {\squarep}{\squarep
   <\arrowtypea`\arrowtypeb`\arrowtypec`\arrowtyped;\width`\height>}}
\def\squarep<#1>[#2`#3`#4`#5;#6`#7`#8`#9]{{
\setsqparms[#1]
\diagram
\putsquarep<\arrowtypea`\arrowtypeb`\arrowtypec`
\arrowtyped;\width`\height>
(0,0)[#2`#3`#4`{#5};#6`#7`#8`{#9}]
\enddiagram
}}                                                 
\def\putptrianglep<#1>(#2,#3)[#4`#5`#6;#7`#8`#9]{{%
\settriparms[#1]%
\xpos=#2 \ypos=#3
\advance\ypos by \height
\puthmorphism(\xpos,\ypos)[#4`#5`{#7}]{\height}{\arrowtypea}a%
\putvmorphism(\xpos,\ypos)[`#6`{#8}]{\height}{\arrowtypeb}l%
\advance\xpos by\height
\putmorphism(\xpos,\ypos)(-1,-1)[``{#9}]{\height}{\arrowtypec}r%
}}

\def\putptriangle{\@ifnextchar <{\putptrianglep}{\putptrianglep
   <\arrowtypea`\arrowtypeb`\arrowtypec;\height>}}
\def\ptriangle{\@ifnextchar <{\ptrianglep}{\ptrianglep
   <\arrowtypea`\arrowtypeb`\arrowtypec;\height>}}
\def\ptrianglep<#1>[#2`#3`#4;#5`#6`#7]{{
\settriparms[#1]
\diagram
\putptrianglep<\arrowtypea`\arrowtypeb`
\arrowtypec;\height>
(0,0)[#2`#3`#4;#5`#6`{#7}]
\enddiagram
}}                                            

\def\putqtrianglep<#1>(#2,#3)[#4`#5`#6;#7`#8`#9]{{%
\settriparms[#1]%
\xpos=#2 \ypos=#3
\advance\ypos by\height
\puthmorphism(\xpos,\ypos)[#4`#5`{#7}]{\height}{\arrowtypea}a%
\putmorphism(\xpos,\ypos)(1,-1)[``{#8}]{\height}{\arrowtypeb}l%
\advance\xpos by\height
\putvmorphism(\xpos,\ypos)[`#6`{#9}]{\height}{\arrowtypec}r%
}}

\def\putqtriangle{\@ifnextchar <{\putqtrianglep}{\putqtrianglep
   <\arrowtypea`\arrowtypeb`\arrowtypec;\height>}}
\def\qtriangle{\@ifnextchar <{\qtrianglep}{\qtrianglep
   <\arrowtypea`\arrowtypeb`\arrowtypec;\height>}}
\def\qtrianglep<#1>[#2`#3`#4;#5`#6`#7]{{
\settriparms[#1]
\width=\height                                
\diagram
\putqtrianglep<\arrowtypea`\arrowtypeb`
\arrowtypec;\height>
(0,0)[#2`#3`#4;#5`#6`{#7}]
\enddiagram
}}

\def\putdtrianglep<#1>(#2,#3)[#4`#5`#6;#7`#8`#9]{{%
\settriparms[#1]%
\xpos=#2 \ypos=#3
\puthmorphism(\xpos,\ypos)[#5`#6`{#9}]{\height}{\arrowtypec}b%
\advance\xpos by \height \advance\ypos by\height
\putmorphism(\xpos,\ypos)(-1,-1)[``{#7}]{\height}{\arrowtypea}l%
\putvmorphism(\xpos,\ypos)[#4``{#8}]{\height}{\arrowtypeb}r%
}}

\def\putdtriangle{\@ifnextchar <{\putdtrianglep}{\putdtrianglep
   <\arrowtypea`\arrowtypeb`\arrowtypec;\height>}}
\def\dtriangle{\@ifnextchar <{\dtrianglep}{\dtrianglep
   <\arrowtypea`\arrowtypeb`\arrowtypec;\height>}}
\def\dtrianglep<#1>[#2`#3`#4;#5`#6`#7]{{
\settriparms[#1]
\width=\height                                
\diagram
\putdtrianglep<\arrowtypea`\arrowtypeb`
\arrowtypec;\height>
(0,0)[#2`#3`#4;#5`#6`{#7}]
\enddiagram
}}

\def\putbtrianglep<#1>(#2,#3)[#4`#5`#6;#7`#8`#9]{{%
\settriparms[#1]%
\xpos=#2 \ypos=#3
\puthmorphism(\xpos,\ypos)[#5`#6`{#9}]{\height}{\arrowtypec}b%
\advance\ypos by\height
\putmorphism(\xpos,\ypos)(1,-1)[``{#8}]{\height}{\arrowtypeb}r%
\putvmorphism(\xpos,\ypos)[#4``{#7}]{\height}{\arrowtypea}l%
}}

\def\putbtriangle{\@ifnextchar <{\putbtrianglep}{\putbtrianglep
   <\arrowtypea`\arrowtypeb`\arrowtypec;\height>}}
\def\btriangle{\@ifnextchar <{\btrianglep}{\btrianglep
   <\arrowtypea`\arrowtypeb`\arrowtypec;\height>}}
\def\btrianglep<#1>[#2`#3`#4;#5`#6`#7]{{
\settriparms[#1]
\width=\height                               
\diagram
\putbtrianglep<\arrowtypea`\arrowtypeb`
\arrowtypec;\height>
(0,0)[#2`#3`#4;#5`#6`{#7}]
\enddiagram
}}

\def\putAtrianglep<#1>(#2,#3)[#4`#5`#6;#7`#8`#9]{{%
\settriparms[#1]%
\xpos=#2 \ypos=#3
{\multiply \height by2
\puthmorphism(\xpos,\ypos)[#5`#6`{#9}]{\height}{\arrowtypec}b}%
\advance\xpos by\height \advance\ypos by\height
\putmorphism(\xpos,\ypos)(-1,-1)[#4``{#7}]{\height}{\arrowtypea}l%
\putmorphism(\xpos,\ypos)(1,-1)[``{#8}]{\height}{\arrowtypeb}r%
}}

\def\putAtriangle{\@ifnextchar <{\putAtrianglep}{\putAtrianglep
   <\arrowtypea`\arrowtypeb`\arrowtypec;\height>}}
\def\Atriangle{\@ifnextchar <{\Atrianglep}{\Atrianglep
   <\arrowtypea`\arrowtypeb`\arrowtypec;\height>}}
\def\Atrianglep<#1>[#2`#3`#4;#5`#6`#7]{{
\settriparms[#1]
\width=\height                                     
\diagram
\putAtrianglep<\arrowtypea`\arrowtypeb`
\arrowtypec;\height>
(0,0)[#2`#3`#4;#5`#6`{#7}]
\enddiagram
}}

\def\putAtrianglepairp<#1>(#2)[#3;#4`#5`#6`#7`#8]{{%
\settripairparms[#1]%
\setpos(#2)%
\settokens[#3]%
\puthmorphism(\xpos,\ypos)[\tokenb`\tokenc`{#7}]{\height}{\arrowtyped}b%
\advance\xpos by\height
\puthmorphism(\xpos,\ypos)[\phantom{\tokenc}`\tokend`{#8}]%
{\height}{\arrowtypee}b%
\advance\ypos by\height
\putmorphism(\xpos,\ypos)(-1,-1)[\tokena``{#4}]{\height}{\arrowtypea}l%
\putvmorphism(\xpos,\ypos)[``{#5}]{\height}{\arrowtypeb}m%
\putmorphism(\xpos,\ypos)(1,-1)[``{#6}]{\height}{\arrowtypec}r%
}}

\def\putAtrianglepair{\@ifnextchar <{\putAtrianglepairp}{\putAtrianglepairp%
   <\arrowtypea`\arrowtypeb`\arrowtypec`\arrowtyped`\arrowtypee;\height>}}
\def\Atrianglepair{\@ifnextchar <{\Atrianglepairp}{\Atrianglepairp%
   <\arrowtypea`\arrowtypeb`\arrowtypec`\arrowtyped`\arrowtypee;\height>}}

\def\Atrianglepairp<#1>[#2;#3`#4`#5`#6`#7]{{
\settripairparms[#1]
\settokens[#2]
\width=\height                                
\diagram
\putAtrianglepairp                            
<\arrowtypea`\arrowtypeb`\arrowtypec`
\arrowtyped`\arrowtypee;\height>
(0,0)[{#2};#3`#4`#5`#6`{#7}]
\enddiagram
}}

\def\putVtrianglep<#1>(#2,#3)[#4`#5`#6;#7`#8`#9]{{%
\settriparms[#1]%
\xpos=#2 \ypos=#3
\advance\ypos by\height
{\multiply\height by2
\puthmorphism(\xpos,\ypos)[#4`#5`{#7}]{\height}{\arrowtypea}a}%
\putmorphism(\xpos,\ypos)(1,-1)[`#6`{#8}]{\height}{\arrowtypeb}l%
\advance\xpos by\height
\advance\xpos by\height
\putmorphism(\xpos,\ypos)(-1,-1)[``{#9}]{\height}{\arrowtypec}r%
}}

\def\putVtriangle{\@ifnextchar <{\putVtrianglep}{\putVtrianglep
   <\arrowtypea`\arrowtypeb`\arrowtypec;\height>}}
\def\Vtriangle{\@ifnextchar <{\Vtrianglep}{\Vtrianglep
   <\arrowtypea`\arrowtypeb`\arrowtypec;\height>}}
\def\Vtrianglep<#1>[#2`#3`#4;#5`#6`#7]{{
\settriparms[#1]
\width=\height                                 
\diagram
\putVtrianglep<\arrowtypea`\arrowtypeb`
\arrowtypec;\height>
(0,0)[#2`#3`#4;#5`#6`{#7}]
\enddiagram
}}

\def\putVtrianglepairp<#1>(#2)[#3;#4`#5`#6`#7`#8]{{
\settripairparms[#1]%
\setpos(#2)%
\settokens[#3]%
\advance\ypos by\height
\putmorphism(\xpos,\ypos)(1,-1)[`\tokend`{#6}]{\height}{\arrowtypec}l%
\puthmorphism(\xpos,\ypos)[\tokena`\tokenb`{#4}]{\height}{\arrowtypea}a%
\advance\xpos by\height
\puthmorphism(\xpos,\ypos)[\phantom{\tokenb}`\tokenc`{#5}]%
{\height}{\arrowtypeb}a%
\putvmorphism(\xpos,\ypos)[``{#7}]{\height}{\arrowtyped}m%
\advance\xpos by\height
\putmorphism(\xpos,\ypos)(-1,-1)[``{#8}]{\height}{\arrowtypee}r%
}}

\def\putVtrianglepair{\@ifnextchar <{\putVtrianglepairp}{\putVtrianglepairp%
    <\arrowtypea`\arrowtypeb`\arrowtypec`\arrowtyped`\arrowtypee;\height>}}
\def\Vtrianglepair{\@ifnextchar <{\Vtrianglepairp}{\Vtrianglepairp%
    <\arrowtypea`\arrowtypeb`\arrowtypec`\arrowtyped`\arrowtypee;\height>}}
\def\Vtrianglepairp<#1>[#2;#3`#4`#5`#6`#7]{{
\settripairparms[#1]
\settokens[#2]
\diagram
\putVtrianglepairp                             
<\arrowtypea`\arrowtypeb`\arrowtypec`
\arrowtyped`\arrowtypee;\height>
(0,0)[{#2};#3`#4`#5`#6`{#7}]
\enddiagram
}}

\def\putCtrianglep<#1>(#2,#3)[#4`#5`#6;#7`#8`#9]{{%
\settriparms[#1]%
\xpos=#2 \ypos=#3
\advance\ypos by\height
\putmorphism(\xpos,\ypos)(1,-1)[``{#9}]{\height}{\arrowtypec}l%
\advance\xpos by\height
\advance\ypos by\height
\putmorphism(\xpos,\ypos)(-1,-1)[#4`#5`{#7}]{\height}{\arrowtypea}l%
{\multiply\height by 2
\putvmorphism(\xpos,\ypos)[`#6`{#8}]{\height}{\arrowtypeb}r}%
}}

\def\putCtriangle{\@ifnextchar <{\putCtrianglep}{\putCtrianglep
    <\arrowtypea`\arrowtypeb`\arrowtypec;\height>}}
\def\Ctriangle{\@ifnextchar <{\Ctrianglep}{\Ctrianglep
    <\arrowtypea`\arrowtypeb`\arrowtypec;\height>}}
\def\Ctrianglep<#1>[#2`#3`#4;#5`#6`#7]{{
\settriparms[#1]
\width=\height                               
\diagram
\putCtrianglep<\arrowtypea`\arrowtypeb`
\arrowtypec;\height>
(0,0)[#2`#3`#4;#5`#6`{#7}]
\enddiagram
}}                                           
\def\putDtrianglep<#1>(#2,#3)[#4`#5`#6;#7`#8`#9]{{%
\settriparms[#1]%
\xpos=#2 \ypos=#3
\advance\xpos by\height \advance\ypos by\height
\putmorphism(\xpos,\ypos)(-1,-1)[``{#9}]{\height}{\arrowtypec}r%
\advance\xpos by-\height \advance\ypos by\height
\putmorphism(\xpos,\ypos)(1,-1)[`#5`{#8}]{\height}{\arrowtypeb}r%
{\multiply\height by 2
\putvmorphism(\xpos,\ypos)[#4`#6`{#7}]{\height}{\arrowtypea}l}%
}}

\def\putDtriangle{\@ifnextchar <{\putDtrianglep}{\putDtrianglep
    <\arrowtypea`\arrowtypeb`\arrowtypec;\height>}}
\def\Dtriangle{\@ifnextchar <{\Dtrianglep}{\Dtrianglep
   <\arrowtypea`\arrowtypeb`\arrowtypec;\height>}}
\def\Dtrianglep<#1>[#2`#3`#4;#5`#6`#7]{{
\settriparms[#1]
\width=\height                              
\diagram
\putDtrianglep<\arrowtypea`\arrowtypeb`
\arrowtypec;\height>
(0,0)[#2`#3`#4;#5`#6`{#7}]
\enddiagram
}}                                          
\def\setrecparms[#1`#2]{\width=#1 \height=#2}%

\def\recursep<#1`#2>[#3;#4`#5`#6`#7`#8]{{%
\width=#1 \height=#2
\settokens[#3]
\settowidth{\tempdimen}{$\tokena$}
\ifdim\tempdimen=0pt
  \savebox{\tempboxa}{\hbox{$\tokenb$}}%
  \savebox{\tempboxb}{\hbox{$\tokend$}}%
  \savebox{\tempboxc}{\hbox{$#6$}}%
\else
  \savebox{\tempboxa}{\hbox{$\hbox{$\tokena$}\times\hbox{$\tokenb$}$}}%
  \savebox{\tempboxb}{\hbox{$\hbox{$\tokena$}\times\hbox{$\tokend$}$}}%
  \savebox{\tempboxc}{\hbox{$\hbox{$\tokena$}\times\hbox{$#6$}$}}%
\fi
\ypos=\height
\divide\ypos by 2
\xpos=\ypos
\advance\xpos by \width
\bfig
\putCtrianglep<-1`1`1;\ypos>(0,0)[`\tokenc`;#5`#6`{#7}]%
\puthmorphism(\ypos,0)[\tokend`\usebox{\tempboxb}`{#8}]{\width}{-1}b%
\puthmorphism(\ypos,\height)[\tokenb`\usebox{\tempboxa}`{#4}]{\width}{-1}a%
\advance\ypos by \width
\putvmorphism(\ypos,\height)[``\usebox{\tempboxc}]{\height}1r%
\efig
}}

\def\recurse{\@ifnextchar <{\recursep}{\recursep<\width`\height>}}

\def\puttwohmorphisms(#1,#2)[#3`#4;#5`#6]#7#8#9{{%
\puthmorphism(#1,#2)[#3`#4`]{#7}0a
\ypos=#2
\advance\ypos by 20
\puthmorphism(#1,\ypos)[\phantom{#3}`\phantom{#4}`#5]{#7}{#8}a
\advance\ypos by -40
\puthmorphism(#1,\ypos)[\phantom{#3}`\phantom{#4}`#6]{#7}{#9}b
}}

\def\puttwovmorphisms(#1,#2)[#3`#4;#5`#6]#7#8#9{{%
\putvmorphism(#1,#2)[#3`#4`]{#7}0a
\xpos=#1
\advance\xpos by -20
\putvmorphism(\xpos,#2)[\phantom{#3}`\phantom{#4}`#5]{#7}{#8}l
\advance\xpos by 40
\putvmorphism(\xpos,#2)[\phantom{#3}`\phantom{#4}`#6]{#7}{#9}r
}}

\def\puthcoequalizer(#1)[#2`#3`#4;#5`#6`#7]#8#9{{%
\setpos(#1)%
\puttwohmorphisms(\xpos,\ypos)[#2`#3;#5`#6]{#8}11%
\advance\xpos by #8
\puthmorphism(\xpos,\ypos)[\phantom{#3}`#4`#7]{#8}1{#9}
}}

\def\putvcoequalizer(#1)[#2`#3`#4;#5`#6`#7]#8#9{{%
\setpos(#1)%
\puttwovmorphisms(\xpos,\ypos)[#2`#3;#5`#6]{#8}11%
\advance\ypos by -#8
\putvmorphism(\xpos,\ypos)[\phantom{#3}`#4`#7]{#8}1{#9}
}}

\def\putthreehmorphisms(#1)[#2`#3;#4`#5`#6]#7(#8)#9{{%
\setpos(#1) \settypes(#8)
\if a#9 %
     \vertsize{\tempcounta}{#5}%
     \vertsize{\tempcountb}{#6}%
     \ifnum \tempcounta<\tempcountb \tempcounta=\tempcountb \fi
\else
     \vertsize{\tempcounta}{#4}%
     \vertsize{\tempcountb}{#5}%
     \ifnum \tempcounta<\tempcountb \tempcounta=\tempcountb \fi
\fi
\advance \tempcounta by 60
\puthmorphism(\xpos,\ypos)[#2`#3`#5]{#7}{\arrowtypeb}{#9}
\advance\ypos by \tempcounta
\puthmorphism(\xpos,\ypos)[\phantom{#2}`\phantom{#3}`#4]{#7}{\arrowtypea}{#9}
\advance\ypos by -\tempcounta \advance\ypos by -\tempcounta
\puthmorphism(\xpos,\ypos)[\phantom{#2}`\phantom{#3}`#6]{#7}{\arrowtypec}{#9}
}}

\def\setarrowtoks[#1`#2`#3`#4`#5`#6]{%
\def\toka{#1}
\def\tokb{#2}
\def\tokc{#3}
\def\tokd{#4}
\def\toke{#5}
\def\tokf{#6}
}
\def\hex{\@ifnextchar <{\hexp}{\hexp<1000`400>}}
\def\hexp<#1`#2>[#3`#4`#5`#6`#7`#8;#9]{%
\setarrowtoks[#9]
\yext=#2 \advance \yext by #2
\xext=#1 \advance\xext by \yext
\bfig
\putCtriangle<-1`0`1;#2>(0,0)[`#5`;\tokb``\tokd]
\xext=#1 \yext=#2 \advance \yext by #2
\putsquare<1`0`0`1;\xext`\yext>(#2,0)[#3`#4`#7`#8;\toka```\tokf]
\advance \xext by #2
\putDtriangle<0`1`-1;#2>(\xext,0)[`#6`;`\tokc`\toke]
\efig
}

\def\skp{\\}

\newtheorem*{conjecture}{Conjecture}
\newtheorem{myconjecture}{Conjecture}
\newtheorem{myproposition}{Proposition}
\newtheorem{myproposition2}{Proposition}
\newtheorem*{theorem}{Theorem}
\newtheorem{mytheorem}{Theorem}
\newtheorem{mytheorem2}{Theorem}
\newtheorem*{corollary}{Corollary}

\author{Carina Curto}

\title{
Matrix model superpotentials and Calabi--Yau spaces: an A-D-E classification
}

\advisor{David R. Morrison}

\department{Mathematics}
\date{2005}

\member{M. Ronen Plesser}
\member{Paul S. Aspinwall}
\member{William L. Pardon}

\begin{document}
\bibliographystyle{abbrv} 

\Copyright

\maketitle

\subject{Mathematics}
\makeabstract

\abstract
We use F. Ferrari's methods relating matrix models to Calabi-Yau
spaces in order to explain Intriligator and
Wecht's ADE classification of $\N=1$ superconformal theories which
arise as RG fixed points of $\N = 1$ SQCD theories with adjoints. The connection between
matrix models and $\N = 1$ gauge theories can be seen as evidence for the
Dijkgraaf--Vafa conjecture.  We find that ADE superpotentials in
the Intriligator--Wecht classification exactly match matrix model
superpotentials obtained from Calabi-Yau's with corresponding ADE
singularities. Moreover, in the additional $\Hat{O}, \Hat{A}, \Hat{D}$ and $\Hat{E}$ cases
we find new singular geometries.  These `hat' geometries are closely related to their ADE counterparts, 
but feature non-isolated singularities. As a byproduct, we give simple descriptions for small resolutions of Gorenstein threefold singularities in terms of transition functions between just two coordinate charts.  
To obtain these results we develop techniques for performing small resolutions and small blow-downs, including an
algorithm for blowing down exceptional $\PP^1$'s.  In particular, we conjecture that small resolutions
for isolated Gorenstein threefold singularities can be obtained by deforming matrix factorizations
for simple surface singularities -- and prove this in the length 1 and length 2 cases.

\newpage

\vspace*{2in}

\begin{flushright}
A mi t\'io, Pompilio Zigrino, por despertar mi inter\'es en la f\'isica.\\
To my professor, Vincent Rodgers, for ten years of pedagogy and support.\\
A mis padres, In\'es y Ra\'ul, por darme el mundo y mucho m\'as.
\end{flushright}
\acknowledgements
First and most importantly, I would like to thank my advisor, Dave Morrison -- it has 
been amazing to work with you.  

Thank you Ronen for introducing me to strings.  
Thanks to Ilarion and Sven for learning with me and for teaching me so many beautiful things.  
Thank you Anda and Janice -- it would have been ten times harder without you!  
Thanks to Carly and Pedro for turning the whole world into playdough.  
And a grateful thanks to Ryan and Kris, who had to live with me.  To Vladimir it goes
without say.

Gracias a Julio Cort\'azar, Che Guevara, Frida Kahlo, Subcomandante Marcos, y Evita Per\'on
por ser figuras facinantes y controversiales.  Me han distraido cuando m\'as lo he necesitado, 
y me mantienen conciente del continente.

Finally, I would like to acknowledge Shakira, Belle \& Sebastian, Elliott Smith, Tori Amos, 
Gal Costa, Manu Chao, and all of the Russian music I listened to while writing.

This work is supported by NSF and VIGRE graduate fellowships.

\tableofcontents

\listoftables	

\listoffigures	

%
%
%
\chapters

\begin{onecolumn}
\chapter{Introduction}\label{ch:intro}
\renewcommand{\a}{\alpha}
\renewcommand{\c}{\gamma}
\renewcommand{\d}{\delta}
\newcommand{\e}{\varepsilon}

Duality has long played an important role in string theory.  Dualities between the five major classes of string theories led Witten to propose in 1995 the idea of an underlying 11-dimensional M-theory \cite{witten2}.  Earlier, the discovery of mirror symmetry as a duality of CFT's \cite{plesser} motivated the proposed duality between type IIA and type IIB string theories, and gave birth to purely geometric versions of mirror symmetry with rich implications for the geometry of Calabi-Yau manifolds \cite{hori,cox}.  By relating physical quantities (correlators, partition functions, spectra) between different theories with geometric input, dualities have uncovered many unexpected patterns in geometry. This has led to surprising conjectures (such as mirror symmetry and T-duality) which not only have important implications for physics, but are interesting and meaningful in a purely geometric light.

Recently, there has been a tremendous amount of work surrounding dualities which relate string theories to other classes of theories.  Maldacena's 1997 AdS/CFT correspondence is perhaps the most famous example of such a duality \cite{maldacena}.  The connection between Chern-Simons gauge theory and string theory was first introduced by Witten in 1992 \cite{witten}. In 1999, Gopakumar and Vafa initiated a program to study the relationship between large N limits of Chern-Simons theory (gauge theory) and type IIA topological string theory (geometry) by using ideas originally proposed by 't Hooft in the 1970's \cite{gopakumar}.  The resulting gauge theory/geometry correspondence led to a conjecture about extremal transitions, often referred to as the ``geometric transition conjecture.''  In the case of conifold singularities, this is more or less understood.  The conifold singularity can be resolved in two very different ways:  (1) with a traditional blow up in algebraic geometry, in which the singular point gets replaced by an exceptional $\PP^1$, or (2) by a deformation of the algebraic equation which replaces the singular point with an $S^3$ whose size is controlled by the deformation parameter (see Figure~\ref{fig:bigpicture}).  The physical degrees of freedom associated to D5 branes wrapping the $\PP^1$ correspond to a 3-form flux through the $S^3$.  The geometric statement is that one can interpolate between the two kinds of resolutions.

In 2002, Dijkgraaf and Vafa expanded this program and proposed new dualities between type IIB topological strings on Calabi-Yau threefolds and matrix models \cite{vafa1,vafa2}.  Due to the symmetry between type IIA and type IIB string theories, this may be viewed as ``mirror'' to the Gopakumar-Vafa conjecture.  By studying the conifold case, they found strong evidence for the matching of the string theory partition function with that of a matrix model whose potential is closely related to the geometry in question.  In particular, a dual version of special geometry in Calabi-Yau threefolds is seen in the eigenvalue dynamics of the associated matrix model \cite{vafa1}. The proposed string theory/matrix model duality has led to an explosion of research on matrix models, a topic which had been dormant since the early 1990's, when it was studied in the context of 2D gravity \cite{ginsparg}.  The connection between string theory and matrix models is of very tangible practical importance, since many quantities which are difficult to compute in string theory are much easier to handle on the matrix model side.

\begin{figure}[h]
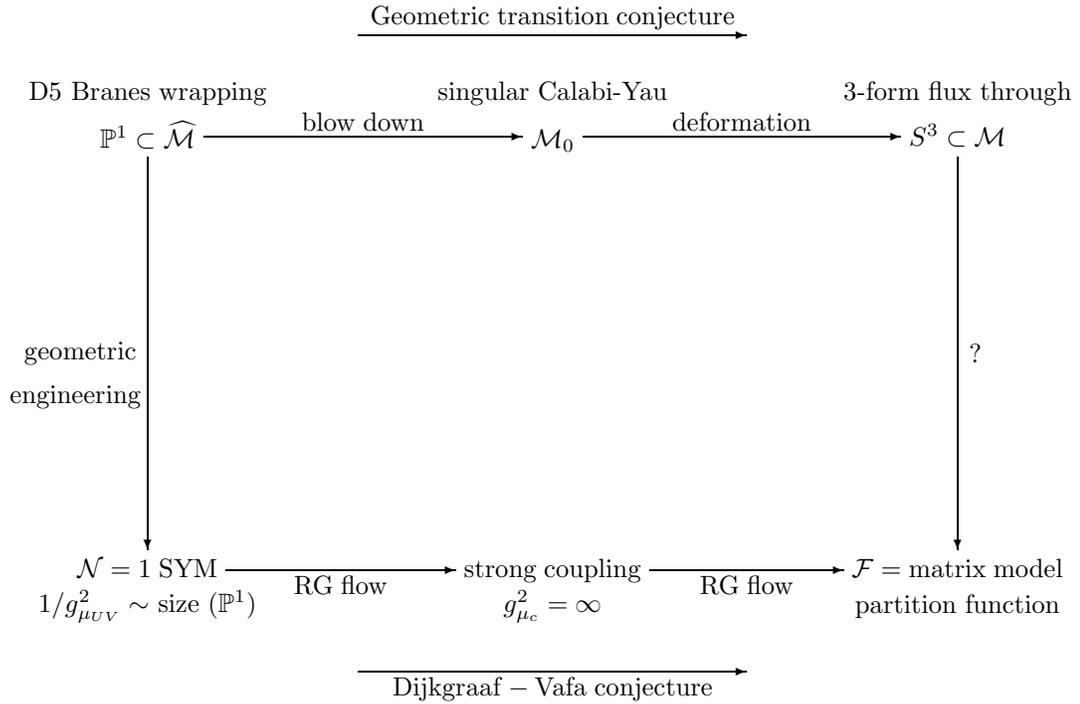

\begin{footnotesize}
$$\bfig
\putsquare<1`1`0`1;1400`1500>(0,250)[\PP^1 \subset \Hat{\M}`\M_0`\mathcal{N}=1 \mathrm{\;SYM}`\mathrm{strong\;coupling};
\mathrm{blow\;down}`\mathrm{geometric}``\mathrm{RG\;flow}]
\putsquare<1`0`1`1;1400`1500>(1400,250)[\phantom{\M_0}`S^3\subset\M`\phantom{\mathrm{strong\;coupling}}`\mathcal{F}=\mathrm{matrix\;model};\mathrm{deformation}``?`\mathrm{RG\;flow}]
\putmorphism(700,2100)(1,0)[``\mathrm{Geometric\;transition\;conjecture}]{1400}1a
\putmorphism(700,-100)(1,0)[``\mathrm{Dijkgraaf-Vafa\;conjecture}]{1400}1b
\put(0,1900){\makebox(0,0){D5 Branes wrapping}}
\put(1400,1900){\makebox(0,0){singular Calabi-Yau}}
\put(2800,1900){\makebox(0,0){3-form flux through}}
\put(0,125){\makebox(0,0){$1/g_{\mu_{UV}}^2 \sim$ size $(\PP^1)$}}
\put(1400,125){\makebox(0,0){$g_{\mu_c}^2 = \infty$}} 
\put(2800,125){\makebox(0,0){partition function}}
\put(-250,850){\makebox(0,0){engineering}}
\efig
$$
\end{footnotesize}
\caption{The big picture}
\label{fig:bigpicture}
\end{figure}

Inspired by these developments (summarized in Figure~\ref{fig:bigpicture}), in 2003 F. Ferrari was led to propose a direct connection between matrix models and the Calabi-Yau spaces of their dual string theories \cite{ferrari}.  It is well known that the solution to a 1-matrix model can be characterized geometrically, in terms of a hyperelliptic curve.  The potential for the matrix model serves as direct input into the algebraic equation for the curve, and the vacuum solutions (distributions of eigenvalues) can be obtained from the geometry of the curve and correspond to branch cuts on the Riemann surface.  The work of Vafa and collaborators on the strongly coupled dynamics of four-dimensional $\mathcal{N}=1$ supersymmetric gauge theories \cite{gopakumar,cachazo2,cachazo,vafa3} suggests that for multi-matrix models, higher-dimensional Calabi-Yau spaces might be useful. Ferrari pursues this idea in \cite{ferrari}, finding evidence that certain multi-matrix models can, indeed, be directly characterized in terms of higher-dimensional (non-compact) Calabi-Yau spaces.

By thinking of the matrix model potential $W(x_1,...,x_M)$ as providing constraints on the deformation space of an exceptional $\PP^1$ within a smooth (resolved) Calabi-Yau $\Hat{\M}$, Ferrari outlines a precise prescription for constructing such smooth geometries directly from the potential.  Specifically, the resolved geometry $\Hat{\M}$ is given by transition functions
\begin{equation}\label{eq:transfns} 
\beta = 1/\gamma, \skop v_1=\gamma^{-n}w_1, \skop
v_2=\gamma^{-m}w_2+\partial_{w_1}E(\gamma,w_1),
\end{equation}
between two coordinate patches $(\gamma,w_1,w_2)$ and $(\beta,v_1,v_2)$, where $\beta$ and $\gamma$ are
stereographic coordinates over an exceptional $\PP^1$. 
The perturbation comes from the ``geometric potential'' $E(\gamma,w)$, which is related to the matrix
model potential $W$ via
\begin{equation}\label{eq:superpotential}
W(x_1,...,x_M)=\dfrac{1}{2\pi i}\oint_{C_0}
\gamma^{-M-1}E(\gamma,\sum_{i=1}^{M}x_i\gamma^{i-1})\mathrm{d}\gamma,
\end{equation}
where $M = n+1 = -m-1.$ We explain this construction in detail in Chapter~\ref{ch:ferrari}.

In the absence of the perturbation term $\partial_{w_1}E(\gamma,w_1)$, the transition functions~\eqref{eq:transfns} simply describe an $\OO(M-1) \oplus \OO(-M-1)$ bundle over a $\PP^1$.  The matrix model superpotential $W$ encodes the constraints on the sections $x_1,...,x_M$ of the bundle due to the presence of $\partial_{w_1}E(\gamma,w_1)$.
Note that this procedure is also invertible.  In other words, given a matrix model superpotential
$W(x_1,...,x_M)$, one can find a corresponding geometric potential $E(\gamma,w_1)$.  As we will discuss in
Chapter~\ref{ch:prelim}, not all perturbation terms $\gamma^j w_1^k$ contribute to the superpotential~\eqref{eq:superpotential}, so there may be many choices of geometric potential for a given $W$. 
Nevertheless, the associated geometry $\Hat{\M}$ is unique \cite[page 634]{ferrari}.  

In 2000, S. Katz had already shown how to codify constraints on versal deformation spaces of curves in terms of a potential function,\footnote{For a rigorous derivation of the D-brane superpotential, see \cite{aspinwall}.} in cases such as~\eqref{eq:transfns} where the constraints are integrable \cite{katz}.  In 2001, F. Cachazo, S. Katz and C. Vafa constructed $\mathcal{N}=1$ supersymmetric gauge theories corresponding to D5 branes wrapping 2-cycles of ADE fibered threefolds \cite{cachazo}. Ferrari studies the same kinds of geometries in a different context, by interpreting the associated potential as belonging to a matrix model, and proposing that the Calabi-Yau geometry encodes all relevant information about the matrix theory.\footnote{Specifically, it is the triple of spaces $\Hat{\M},\M_0,$ and $\M$ that are conjectured to encode the matrix model quantities; the blow-down map $\pi:\Hat{\M} \longrightarrow \M_0$ is the most difficult step towards performing the matrix model computations \cite{ferrari}.}  He is able to verify this in a few examples, and computes known resolvents of matrix models in terms of periods in the associated geometries.  The matching results, as well as Ferrari's solution of a previously unsolved matrix model, suggest that not only can matrix models simplify computations in string theory, but associated geometries from string theory can simplify computations in matrix models.

Many questions immediately arise from Ferrari's construction.  In particular, the matrix model resolvents are not directly encoded in the resolved geometry $\Hat{\M}$, but require knowing the corresponding singular geometry $\M_0$ obtained by blowing down the exceptional $\PP^1$.  It is not clear how to do this blow-down in general. It is also not obvious that a geometry constructed from a matrix model potential in this manner will indeed contain a $\PP^1$ which can be blown down to become an isolated singularity.\footnote{We will see later that the `hat' potentials from the Intriligator-Wecht classification lead to geometries where an entire family of $\PP^1$'s is blown down to reveal a space $\M_0$ with non-isolated singularities.  It is interesting to wonder what the corresponding ``geometric transition conjecture'' should be for these cases.}  Just which matrix models can be ``geometrically engineered'' in this fashion?  What are the corresponding geometries?  Can different matrix model potentials correspond to the same geometry?  If so, what common features of those models does the geometry encode?  Ferrari asks many such questions at the end of his paper \cite{ferrari}, and also wonders whether or not it might be possible to devise an algorithm which will automatically construct the blow-down given the initial resolved space.

Previously established results in algebraic geometry such as Laufer's Theorem \cite{laufer} and the classification of Gorenstein threefold singularities by S. Katz and D. Morrison \cite{morrison} provide a partial answer to these questions.

\begin{theorem}[Laufer 1979]\label{thm:laufer}
Let $\M_0$ be an analytic space of dimension $D \geq 3$ with an isolated singularity at $p$.  Suppose there exists a non-zero holomorphic $D$-form $\Omega$ on $\M_0-\{p\}$.\footnote{This is the Calabi-Yau condition.} Let $\pi:\Hat{\M} \longrightarrow \M_0$ be a resolution of $\M_0$.  Suppose that the exceptional set $A=\pi^{-1}(p)$ is one-dimensional and irreducible. Then $A \cong \PP^1$ and $D=3$.
Moreover, the normal bundle of $\PP^1$ in $\Hat{\M}$ must be either $\mathcal{N} = \OO(-1)\oplus\OO(-1),
\OO\oplus\OO(-2),$ or $\OO(1)\oplus\OO(-3)$. 
\end{theorem}

\noindent Laufer's theorem immediately tells us that we can restrict our search of possible geometries to dimension 3, and that there are only three candidates for the normal bundle to our exceptional $\PP^1$.  In Ferrari's construction, the bundles
$\OO(-1)\oplus\OO(-1), \OO\oplus\OO(-2),$ and $\OO(1)\oplus\OO(-3)$ correspond to zero-, one-, and two-matrix models, respectively.\footnote{This is because these bundles have zero, one, and two independent global sections, respectively.} Following the Dijkgraaf-Vafa correspondence, this puts a limit of 2 adjoint fields on the associated gauge theory, which is precisely the requirement for asymptotic freedom in $\N=1$ supersymmetric gauge theories.  This happy coincidence is perhaps our first indication that the Calabi-Yau geometry encodes information about the RG flow of its corresponding matrix model or gauge theory. 

The condition $M \leq 2$ for the normal bundle $\OO(M-1)\oplus\OO(-M-1)$ 
in Laufer's theorem is equivalent to asymptotic freedom, and reflects the fact that only for asymptotically
free theories can we expect the $\PP^1$ to be exceptional.  In considering matrix model potentials with 
$M \geq 3$ fields, Ferrari points out that the normal bundle to the $\PP^1$ changes with the addition
of the perturbation $\partial_{w_1}E(\gamma,w_1)$, and makes the following conjecture \cite[page 636]{ferrari}:

\begin{conjecture}[RG Flow, Ferrari 2003]  
Consider the perturbed geometry for $m=-n-2$ and associated superpotential $W$.  Let $\N$ be the normal bundle of a $\PP^1$ that sits at a given critical point of $W$.
Let $r$ be the corank of the Hessian of $W$ at the critical point.  Then
$\N = \OO(r-1)\oplus\OO(-r-1)$.
\end{conjecture}

\noindent Ferrari proves the conjecture for $n=1$, and limits himself to two-matrix models
($M=n+1=2$) in the rest of his paper.  Our first result gives evidence in support of the RG flow conjecture
in a more general setting.

\begin{myproposition} For $-M \leq r \leq M$, the addition of the perturbation term $\partial_{w_1}E(\gamma,w_1)=\gamma^{r+1}w_1$ in the transition functions
$$\beta = \gamma^{-1}, \skop v_1 = \gamma^{-M+1}w_1, \skop v_2 =
\gamma^{M+1}w_2+\gamma^{r+1}w_1,$$ changes the bundle from
$\OO(M-1)\oplus\OO(-M-1)$ to $\OO(r-1)\oplus\OO(-r-1).$
In particular, the $M$--matrix model potential
$$W(x_1,...,x_M) = \dfrac{1}{2}\sum_{i=1}^{M-r}x_ix_{M-r+1-i}, \skop (r \geq 0)$$
is geometrically equivalent \footnote{We will call two potentials {\em geometrically equivalent} if they yield the same geometry under Ferrari's construction.} 
to the $r$--matrix model potential
$$W(x_1,...,x_r) = 0.$$
\end{myproposition}

\noindent The proof is given in Section~\ref{sec:bundlechange}.  The fact that the associated superpotential is purely quadratic is satisfying since for quadratic potentials we can often ``integrate out'' fields, giving a field--theoretic intuition for why the geometry associated to an $M$--matrix model can be equivalent to that of an $r$--matrix model, with $r < M$.

Laufer's theorem constrains the dimension, the exceptional set and its normal bundle--but what are the possible singularity types?  In the surface case (complex dimension two), the classification of simple singularities is a classic result.\footnote{An excellent reference for this and other results in singularity theory is \cite{arnold2}.  For a more applications-oriented treatment (with many cute pictures!) see Arnol$'$d's 1991 book \cite{arnold}.
For 15 characterizations of rational double points, see \cite{durfee}.}    As hypersurfaces in $\CC^3$, the distinct geometries are given by: 
\begin{equation}\label{dim2}
\begin{array}{ccccc}
A_k &: & x^2 + y^2 + z^{k+1} &=& 0\\
D_{k+2} &: & x^2 + y^2z + z^{k+1} &=& 0\\
E_6 &: & x^2 + y^3 + z^4 &=& 0\\
E_7 &: & x^2 + y^3 + yz^3 &=& 0\\
E_8 &: & x^2 + y^3 + z^5 &=& 0
\end{array}
\end{equation}




In 1992, Katz and Morrison answered this question in dimension 3 when they characterized the full set of Gorenstein threefold singularities with irreducible small resolutions using invariant theory \cite{morrison}.  In order to do the
classification, Katz and Morrison find it useful to think of threefolds as deformations of surfaces, where the deformation parameter $t$ takes on the role of the extra dimension.  The equations for the singularities can thus
be written in so-called {\em preferred versal form}, as given in Table~\ref{table:versal}.
The coefficients $\alpha_i, \delta_i, \gamma_i,$ and $\varepsilon_i$ are given by invariant polynomials, and are implicity functions of the deformation parameter $t$.

\begin{table}[h]
\begin{center}
$\begin{array}{|c|r@{\  + \ }l|} \hline
 & \multicolumn{2}{c|}{} \\
\mathrm{S} & \multicolumn{2}{c|}{\text{Preferred Versal Form}}  \\
 & \multicolumn{2}{c|}{} \\ \hline
 & \multicolumn{2}{c|}{} \\
A_{n-1}   &  - X Y + Z^n  &   \sum_{i=2}^{n}\a_i Z^{n-i} \\
(n \ge 2) & \multicolumn{2}{c|}{} \\
  & \multicolumn{2}{c|}{} \\
D_n  &  \multicolumn{1}{r@{\  - \ }}{ X^2 + Y^2 Z - Z^{n-1} } &
 \sum_{i=1}^{n-1}\d_{2i} Z^{n-i-1}  + 2 \c_n Y  \\
(n \ge 3)  & \multicolumn{2}{c|}{} \\
  & \multicolumn{2}{c|}{} \\
E_4   &  - X Y + Z^5  &   \e_2 Z^3 + \e_3 Z^2 + \e_4 Z + \e_5 \\
  & \multicolumn{2}{c|}{} \\
E_5  &  \multicolumn{1}{r@{\  - \ }}{  X^2 + Y^2 Z - Z^4 } &
 \e_2 Z^3 - \e_4 Z^2  + 2 \e_5 Y - \e_6 Z - \e_8  \\
  & \multicolumn{2}{c|}{} \\
E_6  &  - X^2 - X Z^2 + Y^3  &
      \e_2 Y Z^2 + \e_5 Y Z + \e_6 Z^2 + \e_8 Y \\
  &  &
      \e_9 Z + \e_{12} \\
  & \multicolumn{2}{c|}{} \\
E_7  &  - X^2 - Y^3 + 16 Y Z^3  &
    \e_2 Y^2 Z + \e_6 Y^2 + \e_8 Y Z + \e_{10} Z^2  \\
  &  &
       \e_{12} Y
     + \e_{14} Z + \e_{18} \\
  & \multicolumn{2}{c|}{} \\
E_8  &  - X^2 + Y^3 - Z^5 &
      \e_2 Y Z^3 + \e_8 Y Z^2 + \e_{12} Z^3 + \e_{14} Y Z \\
  &  &
      \e_{18} Z^2 + \e_{20} Y + \e_{24} Z
     + \e_{30} \\
 & \multicolumn{2}{c|}{} \\ \hline
\end{array}$
\end{center}
\caption[Gorenstein threefold singularities in preferred versal form]{Gorenstein threefold singularities in preferred versal form \cite[page 465]{morrison}.}
\label{table:versal}
\end{table}

\noindent We will also find this representation of the singular threefolds useful, both in performing small resolutions explicitly (Chapter~\ref{ch:smallres}), and in identifying what kinds of singular geometries we get upon blowing down resolved geometries.

In contrast to what one might expect,\footnote{By taking hyperplane sections, one may get surface singularities corresponding to any of the ADE Dynkin diagrams.  A priori, this could indicate that there is an infinite number of families of the threefold singularities with irreducible small resolutions.  What Katz and Morrison discovered is that only a finite number of Dynkin diagrams can arise from ``generic'' hyperplane sections.} there are only a finite number of families of Gorenstein threefold singularities with irreducible small resolutions.  They are distinguished by the Koll\'ar ``length'' invariant,\footnote{\textbf{Defn:} Let $\pi: Y \flecha X$ irreducible small resolution of an isolated threefold singularity $p \in X$.  Let $C = \pi^{-1}(p)$ be the exceptional set.  The {\em length} of $p$ is the length at the generic point of the scheme supported on $C$, with structure sheaf $\OO_Y/\pi^{-1}(m_p,x)$.} and are resolved via small resolution of the appropriate length node in the corresponding Dynkin diagram.  The precise statement of Katz and Morrison's results are given by the following theorem and corollary \cite[page 456]{morrison}:

\begin{theorem}[Katz \& Morrison 1992] The generic hyperplane section of an isolated Gorenstein threefold singularity which has an irreducible small resolution defines one of the primitive partial resolution graphs in Figure~\ref{figure1}.  Conversely, given any such primitive partial resolution graph, there exists an
irreducible small resolution $Y\to X$ whose general hyperplane
section is described by that partial resolution graph.
\end{theorem}

\begin{corollary}  The general hyperplane section of $X$ is uniquely determined
by the length of the singular point $P$.
\end{corollary}

\setlength{\unitlength}{1 true in}

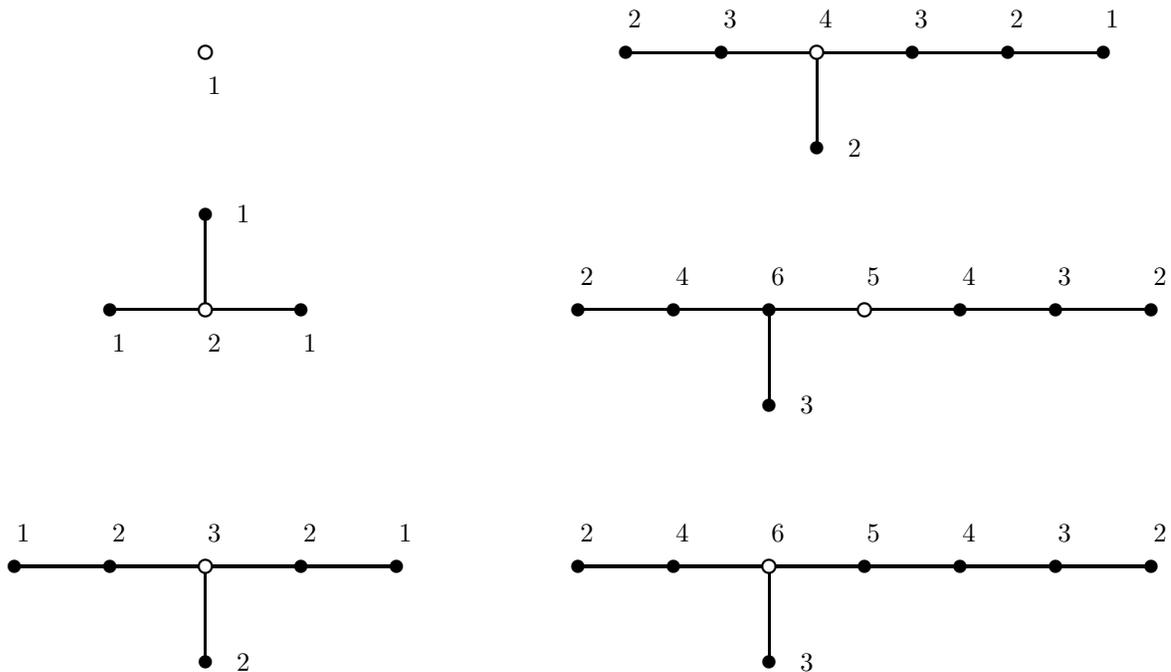
\begin{figure}[t]
\begin{picture}(2,1)(1.9,.5)
\thicklines
\put(2.9,1){\circle{.075}}
\put(2.775,.7){\makebox(.25,.25){\footnotesize 1}}
\end{picture}
\hspace*{\fill}
\begin{picture}(3,1)(1.65,.5)
\thicklines
\put(1.9,1){\circle*{.075}}
\put(1.9,1){\line(1,0){.5}}
\put(2.4,1){\circle*{.075}}
\put(2.4,1){\line(1,0){.4625}}
\put(2.9,1){\circle{.075}}
\put(2.9,.9625){\line(0,-1){.4625}}
\put(2.9,.5){\circle*{.075}}
\put(2.9375,1){\line(1,0){.4625}}
\put(3.4,1){\circle*{.075}}
\put(3.4,1){\line(1,0){.5}}
\put(3.9,1){\circle*{.075}}
\put(3.9,1){\line(1,0){.5}}
\put(4.4,1){\circle*{.075}}
\put(1.775,1.05){\makebox(.25,.25){\footnotesize 2}}
\put(2.275,1.05){\makebox(.25,.25){\footnotesize 3}}
\put(2.775,1.05){\makebox(.25,.25){\footnotesize 4}}
\put(2.925,.375){\makebox(.25,.25){\footnotesize 2}}
\put(3.275,1.05){\makebox(.25,.25){\footnotesize 3}}
\put(3.775,1.05){\makebox(.25,.25){\footnotesize 2}}
\put(4.275,1.05){\makebox(.25,.25){\footnotesize 1}}
\end{picture}

\bigskip

\bigskip

\begin{picture}(2,1)(1.9,.5)
\thicklines
\put(2.4,1){\circle*{.075}}
\put(2.4,1){\line(1,0){.4625}}
\put(2.9,1){\circle{.075}}
\put(2.9,1.0375){\line(0,1){.4625}}
\put(2.9,1.5){\circle*{.075}}
\put(2.9375,1){\line(1,0){.4625}}
\put(3.4,1){\circle*{.075}}
\put(2.275,.7){\makebox(.25,.25){\footnotesize 1}}
\put(2.775,.7){\makebox(.25,.25){\footnotesize 2}}
\put(2.925,1.375){\makebox(.25,.25){\footnotesize 1}}
\put(3.275,.7){\makebox(.25,.25){\footnotesize 1}}
\end{picture}
\hspace*{\fill}
\begin{picture}(3,1)(1.9,.5)
\thicklines
\put(1.9,1){\circle*{.075}}
\put(1.9,1){\line(1,0){.5}}
\put(2.4,1){\circle*{.075}}
\put(2.4,1){\line(1,0){.5}}
\put(2.9,1){\circle*{.075}}
\put(2.9,1){\line(0,-1){.5}}
\put(2.9,.5){\circle*{.075}}
\put(2.9,1){\line(1,0){.4625}}
\put(3.4,1){\circle{.075}}
\put(3.4375,1){\line(1,0){.4625}}
\put(3.9,1){\circle*{.075}}
\put(3.9,1){\line(1,0){.5}}
\put(4.4,1){\circle*{.075}}
\put(4.4,1){\line(1,0){.5}}
\put(4.9,1){\circle*{.075}}
\put(1.775,1.05){\makebox(.25,.25){\footnotesize 2}}
\put(2.275,1.05){\makebox(.25,.25){\footnotesize 4}}
\put(2.775,1.05){\makebox(.25,.25){\footnotesize 6}}
\put(2.925,.375){\makebox(.25,.25){\footnotesize 3}}
\put(3.275,1.05){\makebox(.25,.25){\footnotesize 5}}
\put(3.775,1.05){\makebox(.25,.25){\footnotesize 4}}
\put(4.275,1.05){\makebox(.25,.25){\footnotesize 3}}
\put(4.775,1.05){\makebox(.25,.25){\footnotesize 2}}
\end{picture}

\bigskip

\bigskip

\begin{picture}(2,1)(1.9,.5)
\thicklines
\put(1.9,1){\circle*{.075}}
\put(1.9,1){\line(1,0){.5}}
\put(2.4,1){\circle*{.075}}
\put(2.4,1){\line(1,0){.4625}}
\put(2.9,1){\circle{.075}}
\put(2.9,.9625){\line(0,-1){.4625}}
\put(2.9,.5){\circle*{.075}}
\put(2.9375,1){\line(1,0){.4625}}
\put(3.4,1){\circle*{.075}}
\put(3.4,1){\line(1,0){.5}}
\put(3.9,1){\circle*{.075}}
\put(1.775,1.05){\makebox(.25,.25){\footnotesize 1}}
\put(2.275,1.05){\makebox(.25,.25){\footnotesize 2}}
\put(2.775,1.05){\makebox(.25,.25){\footnotesize 3}}
\put(2.925,.375){\makebox(.25,.25){\footnotesize 2}}
\put(3.275,1.05){\makebox(.25,.25){\footnotesize 2}}
\put(3.775,1.05){\makebox(.25,.25){\footnotesize 1}}
\end{picture}
\hspace*{\fill}
\begin{picture}(3,1)(1.9,.5)
\thicklines
\put(1.9,1){\circle*{.075}}
\put(1.9,1){\line(1,0){.5}}
\put(2.4,1){\circle*{.075}}
\put(2.4,1){\line(1,0){.4625}}
\put(2.9,1){\circle{.075}}
\put(2.9,.9625){\line(0,-1){.4625}}
\put(2.9,.5){\circle*{.075}}
\put(2.9375,1){\line(1,0){.4625}}
\put(3.4,1){\circle*{.075}}
\put(3.4,1){\line(1,0){.5}}
\put(3.9,1){\circle*{.075}}
\put(3.9,1){\line(1,0){.5}}
\put(4.4,1){\circle*{.075}}
\put(4.4,1){\line(1,0){.5}}
\put(4.9,1){\circle*{.075}}
\put(1.775,1.05){\makebox(.25,.25){\footnotesize 2}}
\put(2.275,1.05){\makebox(.25,.25){\footnotesize 4}}
\put(2.775,1.05){\makebox(.25,.25){\footnotesize 6}}
\put(2.925,.375){\makebox(.25,.25){\footnotesize 3}}
\put(3.275,1.05){\makebox(.25,.25){\footnotesize 5}}
\put(3.775,1.05){\makebox(.25,.25){\footnotesize 4}}
\put(4.275,1.05){\makebox(.25,.25){\footnotesize 3}}
\put(4.775,1.05){\makebox(.25,.25){\footnotesize 2}}
\end{picture}

\caption{The 6 types of Gorenstein threefold singularities}
\label{figure1}
\end{figure}

We thus know that there are only a finite number of families of distinct geometries with the desired properties for Ferrari's construction, and they correspond to isolated threefold singularities with small resolutions. While much is known about the resolution of these singularities (they are obtained by blowing up divisors associated to nodes of the appropriate length in the corresponding Dynkin diagram), it is not easy to perform the small blowup explicitly.

\setlength{\unitlength}{1 true in}

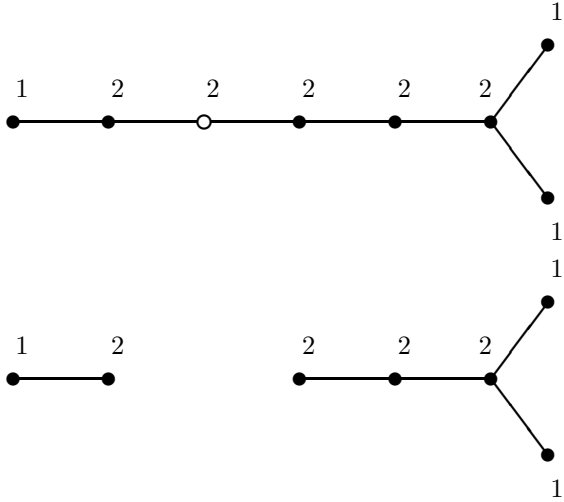
\begin{figure}[h]

\begin{picture}(3,1)(1.65,.5)
\thicklines
\put(1.9,1){\circle*{.075}}
\put(1.9,1){\line(1,0){.5}}
\put(2.4,1){\circle*{.075}}
\put(2.4,1){\line(1,0){.4625}}
\put(2.9,1){\circle{.075}}
\put(2.9375,1){\line(1,0){.4625}}
\put(3.4,1){\circle*{.075}}
\put(3.4,1){\line(1,0){.5}}
\put(3.9,1){\circle*{.075}}
\put(3.9,1){\line(1,0){.5}}
\put(4.4,1){\circle*{.075}}
\put(4.4,1){\line(3,4){.3}}
\put(4.7,1.4){\circle*{.075}}
\put(4.4,1){\line(3,-4){.3}}
\put(4.7,.6){\circle*{.075}}
\put(1.775,1.05){\makebox(.25,.25){\footnotesize 1}}
\put(2.275,1.05){\makebox(.25,.25){\footnotesize 2}}
\put(2.775,1.05){\makebox(.25,.25){\footnotesize 2}}
\put(3.275,1.05){\makebox(.25,.25){\footnotesize 2}}
\put(3.775,1.05){\makebox(.25,.25){\footnotesize 2}}
\put(4.2,1.05){\makebox(.25,.25){\footnotesize 2}}
\put(4.575,1.45){\makebox(.25,.25){\footnotesize 1}}
\put(4.575,.3){\makebox(.25,.25){\footnotesize 1}}
\end{picture}

\bigskip
\bigskip

\begin{picture}(3,1)(1.65,.5)
\thicklines
\put(1.9,1){\circle*{.075}}
\put(1.9,1){\line(1,0){.5}}
\put(2.4,1){\circle*{.075}}
\put(3.4,1){\circle*{.075}}
\put(3.4,1){\line(1,0){.5}}
\put(3.9,1){\circle*{.075}}
\put(3.9,1){\line(1,0){.5}}
\put(4.4,1){\circle*{.075}}
\put(4.4,1){\line(3,4){.3}}
\put(4.7,1.4){\circle*{.075}}
\put(4.4,1){\line(3,-4){.3}}
\put(4.7,.6){\circle*{.075}}
\put(1.775,1.05){\makebox(.25,.25){\footnotesize 1}}
\put(2.275,1.05){\makebox(.25,.25){\footnotesize 2}}
\put(3.275,1.05){\makebox(.25,.25){\footnotesize 2}}
\put(3.775,1.05){\makebox(.25,.25){\footnotesize 2}}
\put(4.2,1.05){\makebox(.25,.25){\footnotesize 2}}
\put(4.575,1.45){\makebox(.25,.25){\footnotesize 1}}
\put(4.575,.3){\makebox(.25,.25){\footnotesize 1}}
\end{picture}
\label{fig:D8blowup}
\caption[The blowup of a node leaves residual surface singularities]{The blow up of the open length 2 node in the top Dynkin diagram causes the $D_8$ surface singularity to split into two lower-order surface singularities.  The residual singularities are of types $A_2$ and $D_5$, as suggested by the remaining parts of the diagram.}
\end{figure}

In the case of simple surface singularities, the blowup of a single node in the associated Dynkin diagram can be obtained via an appropriate matrix factorization of the singular equation \cite{verdier,yoshino}.  Given a 
polynomial $f$ defining a hypersurface, a pair of square matrices $(\phi,\psi)$ such that
$$\phi \cdot \psi = f \cdot \mathbf{1}\;\;\; \and \;\;\; \psi \cdot \phi = f \cdot \mathbf{1}$$
is called a {\em matrix factorization} of $f$.  A particular ADE singularity will have a different (indecomposable) matrix factorization for each node in the associated Dynkin diagram.  
The matrix factorization picks out the ideal corresponding to the node; with this the blowup is straightforward.  The resulting surface usually contains residual singularities, which are characterized by the remaining parts of the original Dynkin diagram once the blown up node has been removed (see Figure 1.3).

Because small resolutions of Gorenstein threefold singularities are obtained by blowing up nodes in the Dynkin diagrams of corresponding hyperplane section surface singularities, we are led to make the following conjecture:

\begin{myconjecture} Small resolutions for Gorenstein threefold singularities can be obtained by deforming matrix factorizations for corresponding surface singularities.
\end{myconjecture}

\noindent We prove this in the length 1 and length 2 cases, yielding the following theorem:

\begin{mytheorem} For Gorenstein threefold singularities of length 1 and length 2, with hyperplane section an $A_n$ or $D_{n+2}$ surface singularity, the small resolution is obtained by deforming the matrix factorization for a node of the same length in the corresponding Dynkin diagram.
\end{mytheorem}

\noindent The proof is given in Chapter~\ref{ch:smallres}.  We are thus able to explicitly perform the small resolutions in the $A_n$ and $D_{n+2}$ cases.  However, even when the blow up has been found, it can be highly non-trivial to find the right coordinates for expressing the answer.  The way in which we identify the superpotential from the resolved geometry requires that this geometry be represented in terms of transition functions over just two coordinate charts covering the exceptional $\PP^1$.\footnote{The techniques in \cite{aspinwall} for computing the superpotential are not constrained in this way; perhaps those methods can be used to overcome this difficulty.}  In the length 2 case ($D_{n+2}$), this difficulty is already apparent, and performing the small resolution stops short of uncovering the corresponding matrix model potential.

The major obstacle in identifying which matrix model corresponds to each of the candidate singular geometries from \cite{morrison} is the absence of a simple description in the form of~\eqref{eq:transfns} for their small resolutions.  This frustration is also expressed in \cite{cachazo}, where the same geometries are used to construct $\mathcal{N}=1$ ADE quiver theories.\footnote{``...the gauge theory description suggests a rather simple global geometric description of the blown up $\PP^1$ for all cases.  However such a mathematical construction is not currently known in the full generality suggested by the gauge theory.  Instead only some explicit blown up geometries are known in detail...''\cite[page 35]{cachazo}}  In the case where the normal bundle to the exceptional $\PP^1$ is $\OO(1)\oplus\OO(-3)$, only Laufer's example \cite{laufer} and its extension by Pinkham and Morrison 
\cite[page 368]{pinkham} was known.  For us, the resolution to this problem came from a timely, albeit surprising, source.

In September, 2003 (just four days after Ferrari's \cite{ferrari} came out on the preprint server!), Intriligator and Wecht posted their results on RG fixed points of $\mathcal{N}=1$ SQCD with adjoints \cite{wecht}.  Using ``a-maximization'' and doing a purely field theoretic analysis, they classify all
relevant adjoint superpotential deformations for 4d $\N=1$ SQCD with
$\N_f$ fundamentals and $\N_a=2$ adjoint matter chiral superfields,
$X$ and $Y$.   The possible RG fixed points, together with the map of possible flows between fixed points, are summarized in Figure 1.4.

\begin{figure}[h]
\label{fig:RGflow}
\begin{picture}(6,3)(0,0)
\put(.8,1.2){\makebox{
$\begin{array}{c|c}
\mathrm{type} & W(X,Y)\\
\hline
&\\
\Hat{O} & 0\\
\Hat{A} & \Tr Y^2\\
\Hat{D} & \Tr XY^2\\
\Hat{E} & \Tr Y^3\\
A_k & \Tr(X^{k+1}+Y^2)\\
D_{k+2} & \Tr(X^{k+1}+XY^2)\\
E_6 & \Tr(Y^3+X^4)\\
E_7 & \Tr(Y^3+YX^3)\\
E_8 & \Tr(Y^3+X^5)
\end{array}$}}
\put(3.7,0){\makebox{\includegraphics[scale=.7]{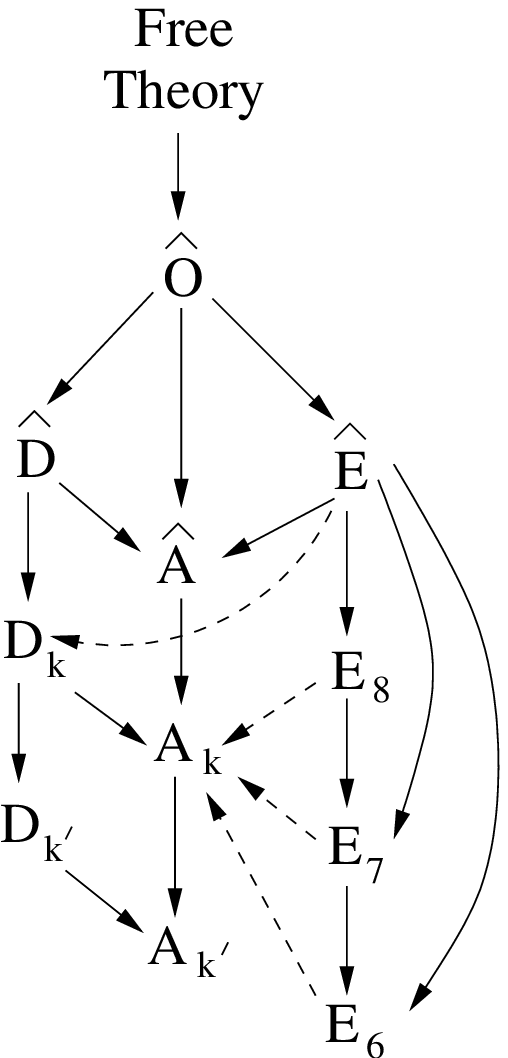}}}
\end{picture}
\caption[Intriligator--Wecht Classification of RG Fixed Points]{Intriligator--Wecht Classification 
of RG Fixed Points.  The diagram on the right shows the map of possible
flows between fixed points.  Dotted lines indicate flow to a particular value of $k$.  Note
that $k' < k$. \cite[pages 3-4]{wecht}}
\end{figure}

Due to the form of the polynomials, Intriligator and Wecht name the relevant superpotential deformations according to the famous ADE classification of singularities in dimensions \nolinebreak 1 and 2~\eqref{dim2}.  There is no geometry in their analysis, however, and they seem surprised to uncover a connection to these singularity types.\footnote{``On the face of it, this has no obvious connection to any of the other known ways in which Arnold's singularities have appeared in mathematics or physics.''\cite[page 3]{wecht}}

Naively, one may speculate that this is the answer!\footnote{In particular, if the Dijkgraaf--Vafa conjecture holds, we should expect any classification of $\N=1$ gauge theories to have a matrix model counterpart.  Verifying such a correspondence thus provides a non-trivial consistency check on the proposed string theory/matrix model duality.}  We make the following conjecture:

\begin{myconjecture}  The superpotentials in Intriligator and Wecht's ADE classification for $\N=1$ gauge theories (equivalently, the polynomials defining simple curve singularities), are precisely the matrix model potentials which yield small resolutions of Gorenstein threefold singularities under Ferrari's construction~\eqref{eq:transfns}.
\end{myconjecture}

\noindent Armed with this new conjecture, we may now run the classification program backwards.  Starting from the resolved space $\Hat{\M}$ given by transition functions over the exceptional $\PP^1$, we can verify the correspondence by simply performing the blow down and confirming that the resulting geometry has the right singularity type.  In particular, the matrix model superpotentials (if correct) give us simple descriptions for the small resolutions of Gorenstein threefold singularities in terms of transition functions as in~\eqref{eq:transfns}.  Like other geometric insights stemming from dualities in string theory, such a result is of independent mathematical interest.

This still leaves us with some major challenges.  As Ferrari points out, there is no known systematic way of performing the blow downs, and our first task is to devise an algorithm which will do so.

\begin{myproposition}\label{prop2}  There is an algorithm for blowing down the exceptional $\PP^1$'s when the resolved geometry is given by simple transition functions~\eqref{eq:transfns} as in Ferrari's construction.
\end{myproposition}

\noindent We prove Proposition~\ref{prop2} by giving an explicit algorithm in Chapter~\ref{ch:algorithm}.  The algorithm can be implemented by computer,\footnote{See the appendix for Maple code.} and searches for global holomorphic functions which can be used to construct the blow down.  Any global holomorphic function on the resolved geometry $\Hat{\M}$ is necessarily constant on the exceptional $\PP^1$, so these functions are natural candidates for coordinates on the blown down geometry $\M_0$, since the $\PP^1$ must collapse to a point.  We are unable to guarantee, however, that all (independent) global holomorphic functions are found by our algorithm.\footnote{The issue is really about how long to run the code.  Given a particular list of input monomials, the algorithm will find all independent global holomorphic functions that can be built using elements from this list.  However, the list of
possible monomials is, of course, infinite.}
The resulting singular space $\M_0$, whose defining equations are obtained by finding relations among the global holomorphic functions, must be checked.  We can verify that we do, indeed, recover the original smooth space by inverting the blow down and performing the small resolution of the singular point. 

We find that this program works perfectly in the $A_k$ (length 1) and $D_{k+2}$ (length 2) cases, lending credence to the idea that the Intriligator--Wecht classification is, indeed, the right answer!  In the exceptional cases, however, a few mysteries arise.  We are only able to find the blow down for the Intriligator--Wecht superpotential $E_7$, and the resulting singular space has a length 3 singularity.  We summarize these results in the following theorem:

\begin{mytheorem}\label{thm2} 
Consider the 2-matrix model potentials $W(x,y)$ in the table below, with corresponding resolved geometries $\Hat{\M}$
$$ \beta = 1/\gamma, \skop v_1=\gamma^{-1}w_1, \skop
v_2=\gamma^{3}w_2+\partial_{w_1}E(\gamma,w_1),$$
given by perturbation terms $\partial_{w_1}E(\gamma,w_1)$.  Blowing down the exceptional $\PP^1$ in each $\Hat{\M}$ yields the following singular geometries $\M_0$:
\begin{table}[h]
$$\begin{array}{c|c|c|c}
\mathrm{type} & W(x,y) & \partial_{w_1}E(\gamma,w_1) & \mathrm{singular\;\;geometry\;\;} \M_0\\
\hline
&&&\\
A_k & \dfrac{1}{k+1}x^{k+1}+\12 y^2 & \gamma^2 w_1^k + w_1
& XY-T(Z^k-T) = 0\\
& & &\\
D_{k+2} & \dfrac{1}{k+1}x^{k+1}+xy^2 & \gamma^2 w_1^k+w_1^2 
& X^2 - Y^2Z + T(Z^{k/2} - T)^2 = 0, \:\:\;\; k \; \even\\
& & & X^2 - Y^2Z - T(Z^k - T^2) = 0, \:\:\;\; k \; \odd\\
& & &\\
E_7 & \dfrac{1}{3}y^3+yx^3 & \gamma^{-1}w_1^2+\gamma w_1^3
& X^2 - Y^3 + Z^5 + 3TYZ^2 + T^3Z = 0.\\
&&&\\
\hline
\end{array}$$
\caption{Superpotentials corresponding to length 1, length 2, and length 3 singularities.}
\label{table:thm2}
\end{table}
\end{mytheorem}

By comparing the above singular geometries with the equations in
preferred versal form (Table~\ref{table:versal}), we immediately identify the $A_k$ and $D_{k+2}$ superpotentials 
as corresponding to length 1 and length 2 threefold singularities, respectively.  For the $E_7$ potential, we first
note that the polynomial $X^2 - Y^3 + Z^5 + 3TYZ^2 + T^3Z$ is in preferred versal form for $E_8$ (with $\varepsilon_8 = -3T$ and $\varepsilon_{24}=-T^3$).\footnote{$T=0$ yields a hyperplane section with $E_8$ surface singularity.}  On the other hand, using Proposition 4 in the proof of the Katz--Morrison classification \cite[pages 499-500]{morrison}, we see that the presence of the monomial $T^3Z$ constrains the threefold singularity type to length 3.  We thus
have the following corollary:

\begin{corollary}
The resolved geometries $\Hat{\M}$ given by the 2-matrix model potentials $A_k$, $D_{k+2}$, and $E_7$ in Table~\ref{table:thm2} are small resolutions for length 1, length 2, and length 3 singularities, respectively.
\end{corollary}

Although simple descriptions of the form~\eqref{eq:transfns} were previously known for small resolutions of length 1 and length 2 Gorenstein threefold singularities (Laufer's example \cite{laufer} in the length 2 case), it is striking that in no other case such a concrete representation for the blowup was known.  Theorem 2, together with its Corollary, show a length 3 example where the small resolution also has an extremely simple form:
$$ \beta = 1/\gamma, \skop v_1=\gamma^{-1}w_1, \skop
v_2=\gamma^{3}w_2+\gamma^{-1}w_1^2+\gamma w_1^3.$$

Chapters~\ref{ch:length1},~\ref{ch:length2}, and~\ref{ch:length3} comprise the proof of Theorem~\ref{thm2}. Missing are examples of length 4, length 5, and length 6 singularity types. For the moment, we are skeptical about whether or not these are describable using geometries that are simple enough to fit into Ferrari's framework.

In some sense the Intriligator--Wecht classification does not contain enough superpotentials; only length 1, 2, and 3 singularities appear to be included.  On the other hand, there are too many: the additional superpotentials  $\Hat{O},\Hat{A},\Hat{D}$ and $\Hat{E}$ have no candidate geometries corresponding to the Katz--Morrison classification of Gorenstein threefold singularities!  What kind of geometries do these new cases correspond to?  And what (if any) is their relation to the original ADE classification?  Using Ferrari's framework and our new algorithmic blow down methods we are able to identify the geometries corresponding to these extra `hat' cases.  We summarize the results in the following theorem:

\begin{mytheorem}\label{thm3} The singular geometries corresponding to the $\Hat{O},\Hat{A},\Hat{D}$ and $\Hat{E}$ cases in the Intriligator--Wecht classification of superpotentials are:
\begin{table}[h]
$$\begin{array}{c|c|c|c}
\mathrm{type} & W(x,y) & \partial_{w_1}E(\gamma,w_1) & \mathrm{singular\;\;geometry\;\;} \M_0\\
\hline
&&&\\
\Hat{O} & 0 & 0 & \CC^3/\ZZ_3\\
& & &\\
\Hat{A} & \12 y^2 & \gamma^2 w_1 & \CC \times \CC^2/\ZZ_2\\
& & &\\
\Hat{D} & xy^2 & \gamma w_1^2 & X^2 + Y^2Z - T^3 = 0\\
& & &\\
\Hat{E} & \dfrac{1}{3} y^3 & \gamma^2 w_1^2 & \mathrm{Spec} (\CC[a,b,u,v]/\ZZ_2)/(b^4 - u^2 - av).\\
&&&\\
\hline
\end{array}$$
\caption{Geometries for the `Hat' cases.}
\end{table}
\end{mytheorem}

The proof of Theorem~\ref{thm3} is the content of Chapter~\ref{ch:hatcases}.  We find that the resolved geometries have full families of $\PP^1$'s which are blown down, and the resulting singular spaces have interesting relations to the ADE cases.  The $\Hat{A}$ geometry is a curve of $A_1$ singularities, while the equation for $\Hat{D}$ looks like the equations for $D_{k+2}$ where the $k$-dependent terms have been dropped.  The identification of new, related geometries obtained by combining Ferrari's framework with the Intriligator--Wecht classification turns out to be one of the most interesting parts of our story.  The presence of these extra geometries may have implications for the relevant string dualities; perhaps the geometric transition conjecture can be expanded beyond isolated singularities.

It is surprising that even in the $\Hat{O}$ case, with $W(x,y) = 0$ superpotential, the geometry is highly non-trivial.  In fact, we find that it is the $A_1, A_k,$ and $\Hat{A}$ cases which are, in some sense, the ``simplest.''  Although the descriptions for the resolved geometries $\Hat{\M}$ in these cases make it appear as though the normal bundles to the exceptional $\PP^1$'s are all $\OO(1)\oplus\OO(-3)$ (as required by a 2-matrix model potential $W(x,y)$), these geometries can all be described with fewer fields.  A straightforward application of Proposition 1 shows that the $\Hat{A}$ case is equivalent to a one-matrix model with $W(x)=0$ superpotential.  Similarly, we will see in Chapter~\ref{ch:length1}
that 
$$W(x,y) = \dfrac{1}{k+1}x^{k+1} + \dfrac{1}{2}y^2 \skop \and \skop W(x) = \dfrac{1}{k+1}x^{k+1}$$
are geometrically equivalent, so the $A_k$ (length 1) cases are also seen to correspond to 1-matrix models, where
the $y$ field has been ``integrated out.''  This is a relief because we know that the exceptional $\PP^1$ after blowing up an $A_k$ singularity should have normal bundle $\OO \oplus \OO(-2)$.  When $k=1$, Proposition 1 further reduces the geometry to that of a 0-matrix model (no superpotential possible!), showing that $A_1$ is the most trivial case, with normal bundle $\OO(-1)\oplus\OO(-1)$.\footnote{In contrast, the normal bundle in the $\Hat{O}$ case is truly $\OO(1)\oplus\OO(-3)$, showing that this geometry requires a 2-matrix model description.}  These results are summarized in Table 1.4, and can be understood as evidence for Ferrari's RG Conjecture.
(Compare with Intriligator and Wecht's map of possible RG flows in Figure 1.4.)

\begin{table}[h]
$$\begin{array}{c|c|c|c}
\mathrm{type} & 2-\mathrm{matrix\;model} & 1-\mathrm{matrix\;model} & 0-\mathrm{matrix\;model}\\
\hline
& & &\\
\Hat{O} & W(x,y) = 0 & &\\
& & &\\
\Hat{A} & W(x,y) = \12 y^2 & W(x) = 0 & \\
& & &\\
A_k & W(x,y) = \dfrac{1}{k+1}x^{k+1}+\12 y^2 & W(x) = \dfrac{1}{k+1}x^{k+1} &\\
& & &\\
A_1 & W(x,y) = \dfrac{1}{2}x^{2}+\12 y^2 & W(x) = \dfrac{1}{2}x^{2} & W = 0\\
&&&\\
\hline
\end{array}$$
\label{table:Acases}
\caption{Geometrically equivalent superpotentials.}
\end{table}

Our analysis also indicates that the names are well-chosen: the $\Hat{A}$ and $\Hat{D}$ geometries are closely related to their $A_{k}$ and $D_{k+2}$ counterparts, and the same might be true for the $\Hat{E}$ case.  The relationship between the $\Hat{A},\Hat{D}$ and $\Hat{E}$ geometries and the ADE singularities is worth exploring for purely geometric reasons.  To summarize, string dualities have told us to enlarge the class of geometries considered in \cite{morrison}, and have pointed us to closely related `limiting cases' of these geometries with non-isolated singularities.

The organization of this thesis is as follows.  Chapter~\ref{ch:ferrari} reviews Ferrari's construction of Calabi-Yau's from matrix model superpotentials.  Chapter~\ref{ch:prelim} contains preliminary computations which are used throughout this work, as well as the proof of Proposition \nolinebreak 1.  Chapter~\ref{ch:node} reviews how to blow up a node in the Dynkin diagram corresponding to a simple surface singularity.  In the $A_n$ and $D_{n+2}$ cases, we rederive matrix factorizations and explicitly show how these are used to obtain the blowups.  In Chapter~\ref{ch:smallres} we prove Theorem \nolinebreak 1, and explore small resolutions of Gorenstein threefold singularities as a ``bottom-up'' approach for identifying corresponding matrix model superpotentials. In Chapter~\ref{ch:algorithm} we develop a ``top-down'' approach to the problem: the Intriligator--Wecht classification provides us with ``guesses'' for the resolved geometries, and an algorithm is described which performs the blow-downs systematically (thereby proving Proposition 2).   Chapters~\ref{ch:length1},~\ref{ch:length2}, and~\ref{ch:length3} present our results in the length 1, 2, and 3 cases, and collectively form the proof of Theorem \nolinebreak 2.  Chapter~\ref{ch:othercases} summarizes our partial knowledge in the case of the remaining $E_6$ and $E_8$ Intriligator--Wecht superpotentials.  Chapter~\ref{ch:hatcases} shows our results in finding the new geometries associated to the $\Hat{O},\Hat{A},\Hat{D}$ and $\Hat{E}$ cases, proving Theorem \nolinebreak 3.  Finally, in Chapter~\ref{ch:future} we discuss open questions and future directions for research.  Appendix A contains our Maple code implemention of the blow-down algorithm.

\end{onecolumn}

\begin{onecolumn}
\chapter{Ferrari's construction}\label{ch:ferrari}
\noindent In this chapter we review some basics of matrix models, and then present Ferrari's construction of Calabi-Yau spaces from matrix model superpotentials in detail.  We illustrate the construction in several examples.

\section{A bit about matrix models.}

\noindent Here we briefly sketch the rudiments of the hermitian one-matrix model, and show the emergence
of the hyperelliptic curve which geometrically encodes the solution.  The general idea that one may be able to 
``solve'' the matrix model by means of an algebraic variety motivates Ferrari to propose Calabi-Yau spaces as higher-dimensional analogues which can encode the solutions of multi-matrix models.  An excellent reference for matrix models in physics is \cite{ginsparg}.

The partition function for a hermitian one-matrix model has the form
$$e^{Z} = \int \d M e^{-\mathrm{tr} W(M)} = \int \prod_{i=1}^N \d \lambda_i \Delta^2(\lambda) e^{-\sum_i W(\lambda_i)},$$
where $M$ is an $N \times N$ hermitian matrix, $W$ a polynomial potential, and the $\lambda_i$'s are the $N$ eigenvalues of $M$.  In going from the first expression to the second, we diagonalize the matrices and pick up the Vandermonde determinant
$$\Delta(\lambda) = \prod_{i<j}(\lambda_i - \lambda_j).$$
To ``solve'' the matrix model means to find a distribution of eigenvalues which minimizes the effective action
$$S(\lambda) = \dfrac{1}{g_s}\sum_i W(\lambda_i) - 2 \sum_{i<j} \log(\lambda_i-\lambda_j).$$
The log term in the effective action comes from exponentiating the contribution from the Vandermonde determinant
$$\Delta(\lambda) = e^{\sum_{i<j}\log(\lambda_i-\lambda_j)},$$
which acts as a Coulomb repulsion between the eigenvalues.

In the large $N$ limit the solution to the matrix model is given by a continuum of eigenvalues, 
with density $\rho$.  The eigenvalues will fill a domain on the real axis; the different components of the domain are known as cuts.  Because of this, the solution to the 1-matrix model can be expressed in terms of an
algebraic equation called the spectral curve, a hyperelliptic curve in the $(x,y)$ plane,
$$y^2 - W'(x)^2 + f(x) = 0,$$
where $f(x)$ depends on the potential.\skp

\section{The geometric framework}

\noindent Here we review Ferrari's construction of non-compact Calabi-Yau's from
matrix model superpotentials, following section 3 of \cite{ferrari}.

The main idea behind the geometric setup is that deformations of the exceptional $\PP^1$ in a resolved geometry $\Hat{\M}$ correspond to adjoint fields in the gauge theory \cite{cachazo}.  Alternatively, the deformation space for a $\PP^1$ wrapped by D-branes can be thought of in terms of matrix models.  The number $N$ of D-branes wrapping the $\PP^1$ gives the size of the matrices ($N \times N$), while the number $M$ of independent sections of the $\PP^1$ normal bundle gives the number of matrices (an $M$--matrix model).

\begin{figure}[h]
\includegraphics[scale=.8]{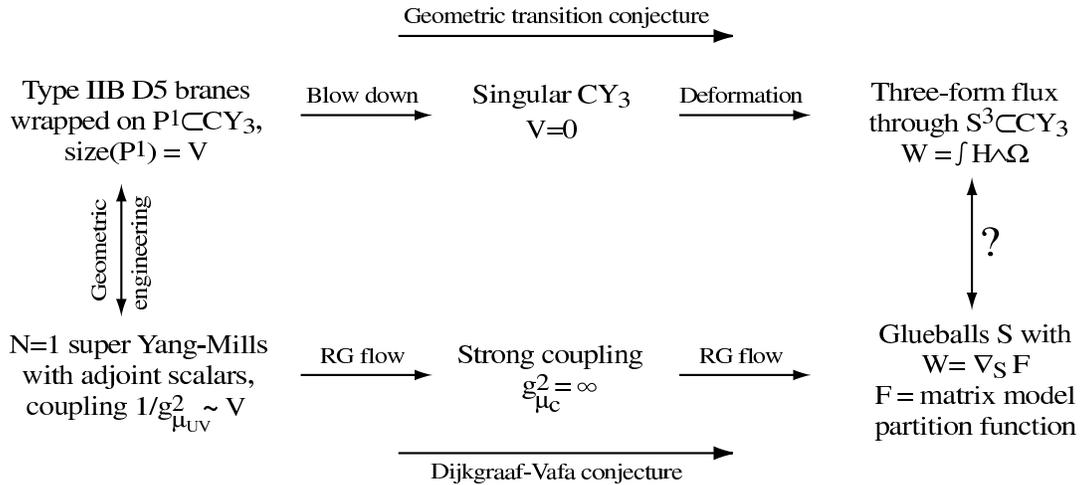}
\label{fig:ferraridiagram}
\caption{Ferrari's diagram \cite[page 631]{ferrari}.}
\end{figure}

Inspired by the string theory dualities (summarized in Figure 2.1), Ferrari develops a recipe to go straight from the matrix model to a Calabi-Yau space.  If the dualities hold, all of the matrix model quantities should be computable from the corresponding geometry.  In this way, Ferrari's prescription provides a non-trivial consistency check on the Gopakumar-Vafa and Dijkgraaf-Vafa conjectures.  Moreover, such a Calabi-Yau space provides a natural higher-dimensional analogue for the spectral (hyperelliptic) curve which encodes the solution to the hermitian one-matrix model.

\begin{figure}[h]
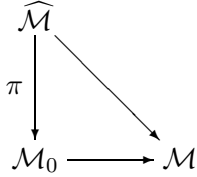

$$\btriangle[\Hat{\M}`\M_0`\M;\pi``]$$
\label{fig:spaces}
\caption{The triple of Calabi-Yau spaces corresponding to a single matrix model.}
\end{figure}

For each theory, there are three relevant Calabi-Yau spaces:
the resolved Calabi-Yau $\Hat{\M}$, the singular Calabi-Yau $\M_0$, and 
the smooth deformed space $\M$ (see Figure~\ref{fig:spaces}).
In short, Ferrari's game consists of the following steps:
\begin{enumerate}
\item Start with an $M$-matrix model matrix model superpotential 
      $W(x_1,...,x_M)$ and construct a smooth Calabi-Yau $\Hat{\M}$.  
      The details of this construction are presented below.
\item Identify the exceptional $\PP^1$ in the resolved space $\Hat{\M}$.
\item Blow down the exceptional $\PP^1$ to get the singular $\M_0$.  The blow down map is 
      $\pi:\Hat{\M} \flecha \M_0$.
\item Perturb the algebraic equation for $\M_0$ to get the smooth deformed space $\M$.
\item From the triple of geometries, compute matrix model quantities (resolvents).
\item Use standard matrix model techniques (loop equations) to check answers in cases where
      the matrix model solution is known.
\end{enumerate}

\noindent In this framework, the matrix model superpotential is encoded in the transition functions defining
the resolved geometry $\Hat{\M}$.  Ferrari shows that a wide variety of matrix model superpotentials
arise in this fashion, and that matrix model resolvents can be computed directly from the geometry.
In other words, the solution to the matrix model is encoded in the corresponding triple of Calabi-Yau's.

The bottleneck to this program is Step 3, the construction of the blow-down map.  While Ferrari's
ad-hoc methods for constructing the blow-down are successful in his particular examples, he does not know
how to construct the blow down in general.  Moreover, it seems the calculation of the blow down map $\pi$ is essentially equivalent to solving the associated matrix model, and hence it would be very useful to have an
algorithm which computes $\pi$ \cite[page 655]{ferrari}.


In this thesis, we are mainly concerned with Steps 1-3.  Our goal is to show how to compute the blow-down map in a large class of examples, and therefore to understand better which singular geometries $\M_0$ arise from matrix models in Ferrari's framework.  In the future, it would be nice to implement the deformation to $\M$ and also to compute matrix model quantities for our examples (Steps 4-6).  For the present, this is beyond our scope.\skp

\subsection{Step 1: Construction of resolved Calabi-Yau}
 
\noindent We now turn to Step 1, the construction of the ``upstairs'' resolved space $\Hat{\M}$ given the matrix
model superpotential $W(x_1,...,x_k)$.  $\Hat{\M}$ is given be transition functions between just two coordinate charts over an exceptional $\PP^1$: $(\beta,v_1,v_2)$ in the first chart, and $(\gamma,w_1,w_2)$ in the second chart.  $\beta$ and $\gamma$ should thought of as stereographic coordinates for the $\PP^1$, with $\beta = \gamma^{-1}$.  The other
coordinates $v_1,v_2$ and $w_1,w_2$ span the normal directions to the $\PP^1$, and have non-trivial transition functions.

We first discuss the case where $W = 0$, in which the Calabi-Yau is the total space of a vector bundle over
the exceptional $\PP^1$.  We then show how a simple deformation of the transition functions leads to constraints on
the sections of the bundle.  The independent sections $x_1,...,x_k$ correspond to matrix degrees of freedom ($k$ independent sections for a $k$-matrix model).  The constraints can be encoded in a potential $W(x_1,...,x_k)$.  When $W$ is non-zero, our geometry $\Hat{\M}$ is no longer a vector bundle -- if the total space were a vector bundle, the sections $x_1,...,x_k$ would be allowed to move freely and therefore satisfy no constraints.  We shall refer to geometries with $W \neq 0$ as ``deformed'' or ``constrained'' bundles.\skp

\noindent \textbf{Pure $\OO(n)\oplus\OO(m)$ bundle}\\
\noindent Consider the following $\Hat{\MM}$ geometry for $n \geq 0$ and $m < 0$:
$$\beta = \gamma^{-1}, \skop v_1=\gamma^{-n}w_1, \skop v_2=\gamma^{-m}w_2.$$
There is an $(n+1)$-dimensional family of $\PP^1$'s that sit at
$$w_1(\gamma)=\sum_{i=1}^{n+1} x_i \gamma^{i-1}, \skop w_2(\gamma)=0.$$
We have no freedom in the $w_2$ coordinate because $m < 0$ precludes $v_2(\beta)$
from being holomorphic whenever $w_2(\gamma)$ is.
$w_1(\gamma)$ and $w_2(\gamma)$ define globally holomorphic sections, and
in the $\beta$ coordinate patch become
$$v_1(\beta)=\sum_{i=1}^{n+1} x_i \beta^{n-i+1},\skop v_2(\beta)=0.$$
The parameters $x_i$ are precisely the fields in the associated superpotential,
and they span the versal deformation space of the $\PP^1$'s. In this case there are no
constraints on the $x_i$'s, which corresponds to the fact that the superpotential is
$$W(x_1,...,x_{n+1})=0.$$
The geometry $\Hat{\M}$ is the total space of a vector bundle, which we might refer to as a ``free'' bundle because
it is not constrained.\skp

\noindent \textbf{Deformed bundle; enter superpotential}\\
\noindent Now consider the deformed geometry (with $n \geq 0$ and $m < 0$):
$$\beta = \gamma^{-1}, \skop v_1=\gamma^{-n}w_1, \skop
v_2=\gamma^{-m}w_2+\partial_{w_1}E(\gamma,w_1),$$
where $E(\gamma,w_1)$ is a function of two complex variables which can be Laurent
expanded in terms of entire functions $E_i$,
$$E(\gamma,w_1) = \sum_{i=-\infty}^{\infty} E_i(w_1)\gamma^i.$$
We call $E$ the ``geometric potential.''
The most general holomorphic section $(w_1(\gamma),w_2(\gamma))$ of
the normal bundle $\N$ to the $\PP^1$'s still has
$$w_1(\gamma)=\sum_{i=1}^{n+1} x_i \gamma^{i-1}, \skop
v_1(\beta)=\sum_{i=1}^{n+1} x_i \beta^{n-i+1},$$
but in order to ensure that $v_2(\beta)$ is holomorphic, the $x_i$'s will
have to satisfy some constraints.  Since a holomorphic $w_2(\gamma)$ can
only cancel poles in $\beta^{-j}$ for $j \geq |m|$, the $x_i$'s must
satisfy $|m|-1$ constraints in order to cancel remaining lower-order poles introduced
by the perturbation.  Hence the versal deformation space of the $\PP^1$ is spanned
by $n+1$ parameters $x_i$ satisfying $|m|-1 = -m-1$ constraints.  

For the $\PP^1$ to be isolated we need $n+1 = -m-1$, and we denote this quantity
(the number of fields) by $M$.  The constraints are integrable, and
equivalent to the extremization $\d W=0$ of the corresponding superpotential
$W(x_1,...,x_M).$  The $\PP^1$'s then sit at the critical points of the
superpotential, in the sense that for critical values of the $x_i$'s,
the pair $(w_1(\gamma), w_2(\gamma))$ will be a global holomorphic section
defining a $\PP^1$.\skp

\noindent \textbf{Summary: General transition functions for $\Hat{\M}$}\\
The resolved geometry $\Hat{\M}$ is described by two coordinate patches $(\gamma,w_1,w_2)$ and
$(\beta,v_1,v_2)$, with transition functions
$$\beta = 1/\gamma, \skop v_1=\gamma^{-n}w_1, \skop
v_2=\gamma^{-m}w_2+\partial_{w_1}E(\gamma,w_1).$$
In the absence of the $\partial_{w_1}E(\gamma,w_1)$ term, this would simply be an $\OO(n) \oplus \OO(m)$
bundle over the $\PP^1$ parametrized by the stereographic coordinates $\gamma$ and $\beta$.
The perturbation comes from the ``geometric potential'' $E(\gamma,w)$, which can be expanded as
$$E(\gamma,w) = \sum_{i=-\infty}^{+\infty} E_i(w)\gamma^i.$$
\skp

\noindent \textbf{The superpotential}\\
The matrix model superpotential encodes the constraints on the sections $x_1,...,x_M$ due to
the presence of the perturbation term $\partial_{w_1}E(\gamma,w_1)$ in the defining transition functions
for $\Hat{\M}$.  It can be obtained directly from the geometric potential via
$$W(x_1,...,x_M)=\dfrac{1}{2\pi i}\oint_{C_0}
\gamma^{-M-1}E(\gamma,\sum_{i=1}^{M}x_i\gamma^{i-1})\d\gamma,$$
where $$M = n+1 = -m-1.$$
The contour integral is meant as a bookkeeping device; $C_0$ should be taken as
a small loop encircling the origin.  The integral is a compact notation used
by Ferrari to encode all of the constraints at once.  The general method for going from
transition function perturbation (geometric potential) to superpotential was first presented in \cite{katz}.

This procedure is also invertible.  In other words, given a matrix model superpotential
$W(x_1,...,x_M)$ one can find a corresponding geometric potential $E(\gamma,w_1)$, 
and therefore construct the associated geometry.  $E$ is not in general unique; from the expression
for $W$ one can see that terms can always be added to the geometric potential which will not contribute
to the residue of the integrand, and hence will not affect the superpotential.  Such terms have no
effect on the geometry, however.

Going from $W$ to $E$ is essentially the implementation of Step 1
in Ferrari's game.  We now turn to Step 2, which is to locate the exceptional $\PP^1$'s.\skp

\subsection{Step 2: Locating the $\PP^1$'s}

\noindent The first task in constructing the blow-down maps is figuring out 
where the $\PP^1$'s that we want to blow down are located.
We will mostly be interested in the $M=2$ case,
$$\beta=\gamma^{-1},\skop v_1=\gamma^{-1}w_1, \skop v_2=\gamma^{3}w_2+\partial_{w_1}E(\gamma,w_1),$$
where we always have
$$w_1(\gamma)=x+\gamma y, \skop v_1(\beta)=\beta x + y,$$
with $x$ and $y$ critical points of $W(x,y)$ at the $\PP^1$'s.  Depending
on the form of the perturbation $\partial_{w_1}E(\gamma,w_1)$, $w_2(\gamma)$
will be chosen to cancel poles of order $\geq 3$.  The requirement that
$v_2(\beta)$ be holomorphic will fix $x$ and $y$ values to be the same as
for the critical points of $W(x,y)$.\skp

\subsection{Step 3: Finding the blow-down map in Ferrari's examples}\label{step3}

\noindent As previously mentioned, Ferrari has no systematic way of constructing the blow-down map
$$\pi:\Hat{\M} \flecha \M_0.$$  He successfully finds the blow-down in several examples, however, through clever
but ad-hoc methods.  Here we show how Ferrari finds the blow-down by discussing the main examples from his paper \cite{ferrari}.\skp

\noindent {\bf Example 1} \cite[pages 636-7]{ferrari}\\
\noindent The resolved geometry $\Hat{\M}$ given by transition functions
$$\beta = \gamma^{-1}, \skop v_1 = w_1, \skop v_2 = \gamma^2 w_2 + \gamma P(w_1),$$
engineers the one-matrix model with potential $W(x)$ such that $W'(x) = P(x)$.  The $\PP^1$'s
sit at
$$w_1(\gamma) = x = v_1(\beta), \skop w_2(\gamma) = 0 = v_2(\beta),$$
with $P(x) = 0$.  In this case, global holomorphic functions are not difficult to find:
\begin{eqnarray*}
\pi_1 &=& w_1 = v_1,\\
\pi_2 &=& 2\gamma w_2 + P(w_1) = 2\beta v_2 - P(v_1),\\
\pi_3 &=& w_2 - \gamma P(w_1) - \gamma^2 w_2 = -v_2 - \beta P(v_1) + \beta^2 v_2,\\
\pi_4 &=& w_2 + \gamma P(w_1) + \gamma^2 w_2 = v_2 - \beta P(v_1) + \beta^2 v_2.
\end{eqnarray*}
These functions collapse the isolated $\PP^1$'s to points, and provide a birational
isomorphism outside the set of $\PP^1$'s.  The singular blown down geometry is uncovered
by observing that the $\pi_i$'s satisfy the algebraic equation
$$\M_0: \pi_4^2 = \pi_3^2 + \pi_2^2 - P^2(\pi_1).$$
Taking hyperplane sections $\pi_1 = x_c$, for critical points of the potential $P(x_c) = 0,$
we see that we have $A_1$ surface singularities which are simultaneously resolved upon inverting
the blow-down.\skp

\noindent {\bf Example 2} \cite[pages 637-8]{ferrari}\\
\noindent Consider the superpotential
$$W(x,y) = V(x) + U(y) - xy,$$
with resolved geometry $\Hat{\M}$ of the form
$$\beta = \gamma^{-1}, \skop v_1 = \gamma^{-1} w_1, \skop v_2 = \gamma^3 w_2 + \gamma Q(w_1/\gamma)
+ \gamma^2 P(w_1) - \gamma w_1,$$
where $U'(y) = Q(y)$ and $V'(x) = P(x)$.  The blow-down is constructed by trial and error.  Starting
with the ansatz
$$\pi_1 = -\beta v_2 + ... = -\gamma^2 w_2 - Q(w_1/\gamma) + w_1 - \gamma P(w_1) + ...,$$
we can move the non-holomorphic term $Q(w_1/\gamma)$ to the left-hand-side, yielding
$$\pi_1 = -\beta v_2 + Q(v_1) = -\gamma^2 w_2 + w_1 - \gamma P(w_1).$$
Similarly we get 
$$\pi_2 = \gamma w_2 + P(w_1) = v_1 - \beta Q(v_1) + \beta^2 v_2.$$
$\pi_3$ and $\pi_4$ require a bit more creativity, and suitable functions are
\begin{eqnarray*}
\pi_3 &=& -\beta^3 v_2 + \beta^2 Q(v_1) + \beta P(\pi_1) - \beta v_1 = \dfrac{P(\pi_1) - P(w_1)}{\gamma}-w_2,\\
\pi_4 &=& \gamma^3 w_2 + \gamma^2 P(w_1) + \gamma Q(\pi_2) - \gamma w_1 = 
\dfrac{Q(\pi_2) - Q(v_1)}{\beta} + v_2.
\end{eqnarray*}
The singular Calabi--Yau is
$$\M_0: \pi_3\pi_4 = (P(\pi_1) - \pi_2)(Q(\pi_2) - \pi_1).$$
As in Example 1, this is a length 1 threefold singularity.  Why does the normal bundle to the $\PP^1$ in the 
resolved geometry look like $\OO(1) \oplus \OO(-3)$?  We know that for $A_1$ singularities it should be $\OO(-1)\oplus\OO(-1)$, and for $A_k$ ($k > 1$) the normal bundle is $\OO\oplus\OO(-2)$.  We will discuss this
puzzle in detail in Section~\ref{example2}.\skp

\noindent {\bf Example 3} \cite[pages 638-9]{ferrari}\\
\noindent This is a generalization of Laufer's example \cite{laufer}, which can also be found in \cite{cachazo}.
The superpotential is
$$W(x,y) = xy^2 + V(x) + y U(y^2),$$
and the resolved geometry $\Hat{\M}$ has the form\footnote{Note that this is almost identical to
the previous example, except for the last perturbation term $w_1^2$.  The geometry in this case,
however, is quite different and corresponds to a length 2 threefold singularity.}
$$\beta = \gamma^{-1}, \skop v_1 = \gamma^{-1} w_1, \skop v_2 = \gamma^3 w_2 + \gamma Q(w_1/\gamma)
+ \gamma^2 P(w_1) + w_1^2,$$
where
\begin{eqnarray*}
V'(x) &=& P(x) = -F_1(x^2) - xF_2(x^2),\\
Q(y) &=& U(y^2) + 2y^2 U'(y^2) = -G(y^2).
\end{eqnarray*}
Again, the blow-down map is constructed using the tricks in Examples 1 and 2, although in this
case $\pi_1,...,\pi_4$ have a rather more complicated form (see \cite[page 639]{ferrari} for details).
The singular geometry is as usual determined by finding a relation among the $\pi_i$'s, 
$$\M_0: \pi_2\left[ (\pi_2 - F_1(\pi_1))^2 - \pi_1 F_2^2(\pi_1)\right] = 
\pi_4^2 - \pi_1\pi_3^2 - G(\pi_2) \left[ (\pi_2 - F_1(\pi_1))\pi_3 + \pi_4 F_2(\pi_1)\right].$$
Despite the complicated form, special cases reveal that this threefold has length 2 
(type $D_{n+2}$) singularities, as in Laufer's original example.

\section{Example:  $A_k$}\label{sec:A_k}

We now illustrate Steps 1-3 in a simple example.
Consider the matrix model potential $$W(x)=\dfrac{1}{k+1}x^{k+1}.$$
Since there is only 1 field, 
the resolved geometry $\Hat{\M}$ is given by transition functions
$$\beta = \gamma^{-1}, \skop v_1 = w_1, \skop v_2 = \gamma^2 w_2 + \gamma w_1^k,$$
for an $\OO\oplus\OO(-2)$ bundle over the exceptional $\PP^1$.
To locate the $\PP^1$, we first note that 
$$w_1(\gamma)=x=v_1(\beta)$$ 
are the only holomorphic sections for the $\OO$ line bundle.  Substituting this into
the transition function for $v_2$ yields
$$v_2(\beta) = \beta^{-2} w_2(\gamma) + \beta^{-1}x^k,$$
which is only holomorphic if $x^k = w_2(\gamma) = 0$.  Therefore we have a single $\PP^1$ located at
$$w_1(\gamma)=w_2(\gamma)=0, \skop v_1(\beta)=v_2(\beta)=0.$$
Note that the position of the $\PP^1$ corresponds exactly to the critical point of the superpotential:
$$\d W = x^k dx = 0 \skop \implies \skop x^k = 0.$$

\noindent \textbf{The blow-down}\\  
\noindent To find the blow down map $\pi$, we must look for global holomorphic functions (which will necessarily be constant on the $\PP^1$).  We can immediately write down
\begin{eqnarray*}
\pi_1 &=& v_1 = w_1,\\
\pi_2 &=& v_2 = \gamma^2 w_2 + \gamma w_1^k,
\end{eqnarray*}
which are independent.  Moreover, notice the combination
$\beta v_2 - v_1^k = \gamma w_2$.  This gives us
\begin{eqnarray*}
\pi_3 &=& \beta v_2 - v_1^k = \gamma w_2,\\
\pi_4 &=& \beta^2 v_2 - \beta v_1^k = w_2.
\end{eqnarray*}
Since $\beta = \pi_4/\pi_3$, we have immediately from the definition of $\pi_3$ the relation
$$\M_0:\;\; \pi_3^2 = \pi_4\pi_2-\pi_3\pi_1^k.$$
This corresponds to an $A_k$ (length 1) singularity!\skp

\noindent \textbf{The blowup}\\ 
\noindent We check our computation by inverting the blow-down.  If we define
$$v_3 = \beta v_2 - v_1^k, \skop \mathrm{and} \skop w_3 = \gamma w_2 + w_1^k,$$ 
we can write
$$\begin{array}{ccccc}
\pi_1 &=& v_1 &=& w_1,\\
\pi_2 &=& v_2 &=& \gamma w_3,\\
\pi_3 &=& v_3 &=& \gamma w_2,\\
\pi_4 &=& \beta v_3 &=& w_2. 
\end{array}$$
In particular
$$\beta = \pi_4/\pi_3 = \gamma^{-1}.$$
This suggests that to recover the small resolution $\Hat{\M}$ we should blow up 
$$\pi_3 = \pi_4 = 0 \skop \mathrm{in} \skop \M_0:\;\; \pi_3^2 = \pi_4\pi_2-\pi_3\pi_1^k.$$
Denoting the $\PP^1$ coordinates by $[\beta:\gamma]$ and imposing
the relation $$\beta \pi_3 = \gamma \pi_4$$ in the blowup, we find in each chart 
$$\begin{array}{c|c}
(\gamma = 1)& (\beta = 1)\\
\hline
&\\
\pi_4 = \beta \pi_3 & \pi_3 = \gamma \pi_4\\
\pi_3 = \beta \pi_2 - \pi_1^k & \gamma^2\pi_4 = \pi_2 - \gamma\pi_1^k\\
&\\
(\beta,\pi_1,\pi_2) & (\gamma,\pi_1,\pi_4)
\end{array}$$

\noindent The transition functions between the two charts are easily found to be
$$\beta = \gamma^{-1}, \skop \pi_1 = \pi_1, \skop \pi_2 = \gamma^2 \pi_4 + \gamma \pi_1^k.$$
Identifying with the original coordinates, we find
$$\beta = \gamma^{-1}, \skop v_1 = w_1, \skop v_2 = \gamma^2 w_2 + \gamma w_1^k.$$
This is exactly what we started with!

Note that for $k=1$, we have a bundle-changing superpotential (see Section~\ref{sec:bundlechange}), and the normal bundle to the exceptional $\PP^1$ is $\OO(-1)\oplus \OO(-1)$ instead of $\OO \oplus \OO(-2)$.

\end{onecolumn}

\begin{onecolumn}
\chapter{Preliminary computations}\label{ch:prelim}
In this chapter we do some preliminary calculations in order to investigate more carefully
the interplay between perturbation terms, superpotentials, and bundle structure.  Due to 
Laufer's theorem (see Introduction), we focus on the one- and two-matrix model cases.  After
computing which perturbation terms come from superpotentials, we turn our attention to a puzzle that arises
from Ferrari's example 2: a $\PP^1$ with what seems to be an $\OO(1)\oplus\OO(-3)$ normal bundle is blown
down to reveal a geometry which should have come from blowing down a $\PP^1$ with normal bundle
$\OO \oplus \OO(-2)$.  We resolve this apparent contradiction in Section~\ref{sec:bundlechange} by investigating terms in the superpotential which not only constrain the sections of the normal bundle, but 
actually change the local structure.  This leads us to prove Proposition \nolinebreak 1:

\begin{myproposition2} For $-M \leq r \leq M$, the addition of the perturbation term $\partial_{w_1}E(\gamma,w_1)=\gamma^{r+1}w_1$ in the transition functions
$$\beta = \gamma^{-1}, \skop v_1 = \gamma^{-M+1}w_1, \skop v_2 =
\gamma^{M+1}w_2+\gamma^{r+1}w_1,$$ changes the bundle from
$\OO(M-1)\oplus\OO(-M-1)$ to $\OO(r-1)\oplus\OO(-r-1).$
In particular, the $M$--matrix model potential
$$W(x_1,...,x_M) = \dfrac{1}{2}\sum_{i=1}^{M-r}x_ix_{M-r+1-i}, \skop (r \geq 0)$$
is geometrically equivalent \footnote{We will call two potentials {\em geometrically equivalent} if they yield the same geometry under Ferrari's construction.} 
to the $r$--matrix model potential
$$W(x_1,...,x_r) = 0.$$
\end{myproposition2}

\noindent We finish by presenting a resolution of the ideal sheaf, which introduces new coordinates
which will be useful for understanding some later examples.

\section{Perturbations which contribute to superpotential}

\noindent Consider the transition functions for a (putative) $\OO(M-1)\oplus\OO(-M-1)$ bundle with
perturbation:
$$\beta = \gamma^{-1},\skop v_1=\gamma^{1-M}w_1, \skop v_2=\gamma^{M+1}w_2+\partial_{w_1}E(\gamma,w_1).$$
Not all possible perturbation terms
$$\partial_{w_1}E(\gamma,w_1) = \gamma^n w_1^m, \skop n \in \ZZ, \;\;\; m \in \ZZ_{\geq 0}$$
will contribute to the superpotential.
Here we compute the superpotential associated to such perturbations
\begin{eqnarray*}
W(x_1,...,x_M)&=&\dfrac{1}{2\pi i}\oint_{C_0}
\gamma^{-M-1}E(\gamma,\sum_{i=1}^{M}x_i\gamma^{i-1})\d\gamma\\
&=& \dfrac{1}{2\pi i}\oint_{C_0}
    \dfrac{\d\gamma}{m+1}\gamma^{n-M-1}(\sum_{i=1}^{M}x_i\gamma^{i-1})^{m+1}
\end{eqnarray*}
to find exactly which terms contribute in
the 1-matrix and 2-matrix models.  This will provide a useful reference
for future computations.\skp

\noindent {\bf Linear perturbations}\\
\noindent In the case of linear perturbations ($m=1$), the terms
$$\partial_{w_1}E(\gamma,w_1) = \gamma^n w_1$$
yield purely quadratic superpotentials:
\begin{eqnarray*}
W(x_1,...,x_M)&=&\dfrac{1}{2\pi i}\oint_{C_0}
\12 \gamma^{n-M-1}(\sum_{i=1}^{M}x_i\gamma^{i-1})^2\d\gamma\\
&=& \12 \sum_{i=1}^{M}x_i x_{M-n+2-i}
\end{eqnarray*}
We will see later that these are the terms which can change bundle structure.\skp

\noindent {\bf 1--matrix model}\\
\noindent The expression for the potential is
\begin{eqnarray*}
W(x)&=&\dfrac{1}{2\pi i}\oint_{C_0}
    \dfrac{\d\gamma}{m+1}\gamma^{n-2}x^{m+1} \\
    &=&\left\{\begin{array}{cc} 0 & n \neq 1 \\ \dfrac{1}{m+1}x^{m+1} & n=1 \end{array}\right.
\end{eqnarray*}
The only perturbation terms
$$\partial_{w_1}E(\gamma,w_1) = \gamma^n w_1^m$$
which contribute satisfy $$n=1.$$
In particular, a superpotential of the form
$$W(x) = \dfrac{1}{m+1}x^{m+1}$$
corresponds to a perturbation term of the form
$$\partial_{w_1}E(\gamma,w_1) = \gamma w_1^m.$$

\noindent {\bf 2--matrix model}\\
\noindent The expression for the potential is
\begin{eqnarray*}
W(x,y)&=&\dfrac{1}{2\pi i}\oint_{C_0}
    \dfrac{\d\gamma}{m+1}\gamma^{n-3}(x+\gamma y)^{m+1}  \\
      &=& \dfrac{1}{m+1}{m+1 \choose 2-n} x^{m+n-1}y^{2-n}
\end{eqnarray*}
The only perturbation terms
$$\partial_{w_1}E(\gamma,w_1) = \gamma^n w_1^m$$
which contribute satisfy $0 \leq 2-n \leq m+1,$ equivalently
$$1-m \leq n \leq 2.$$
In particular, a superpotential of the form
$$W(x,y) = \dfrac{1}{j+k}{j+k \choose k} x^j y^k$$
corresponds to a perturbation term of the form
$$P(j,k) \equiv \partial_{w_1}E(\gamma,w_1) = \gamma^{2-k}w_1^{j+k-1}.$$

\noindent{\bf Remark}\\
\noindent The change of variables
$$ x \leftrightarrow y, $$
in the superpotential is equivalent to interchanging $$j \leftrightarrow k.$$
A quick computation shows that
$$ P(k,j) = \gamma^{2-j}w_1^{j+k-1} = \gamma^{k-j} P(j,k).$$
In other words,
$$W(x,y) \mapsto W(y,x) \;\;\; \Longrightarrow \;\;\; P(j,k) \mapsto  \gamma^{k-j} P(j,k).$$

\section{Ferrari's example 2}\label{example2}

Recall Ferrari's second example from Section~\ref{step3}, where from $\Hat{\M}$ transition functions
$$\beta = \gamma^{-1}, \skop v_1 = \gamma^{-1} w_1, \skop v_2 = \gamma^3 w_2 + \gamma Q(w_1/\gamma)
+ \gamma^2 P(w_1) - \gamma w_1,$$
he performs the blow-down to get singular geometry $\M_0$:
$$ XY = (P(t)-Z)(Q(Z)-t).$$
We will check the result by blowing up the Weil divisor
$$Y = P(t)-Z = 0,$$
to verify that the resolved space is, indeed, what Ferrari started with.
We introduce coordinates $[\gamma:\beta]$ on the exceptional $\PP^1$, and
the relation
$$\gamma Y = \beta (P(t)-Z).$$
The blowup is summarized in two charts
$$\begin{array}{c|c}
(\gamma,Y,t) &  (\beta, Z, X)\\
\hline 
&\\
\gamma Y = P(t) -Z & Y = \beta (P(t) -Z)\\
&\\
X = \gamma (Q(Z)-t) & \beta X = Q(Z)-t
\end{array}$$
with transition functions
$$\begin{array}{ccc|ccc}
&&&&&\\
\beta &=& 1/\gamma & \gamma &=& 1/\beta\\
&&&&&\\
Z &=& -\gamma Y +P(t) & Y &=& -\beta Z + \beta P(t)\\
&&&&&\\
X &=& -\gamma t + \gamma Q(Z) & t &=& -\beta X + Q(Z).\\
&&&&&
\end{array}$$
$$ $$

\noindent{\bf The puzzle}\\
\noindent Our blowup clearly yields a geometry with normal bundle $\OO(-1)\oplus\OO(-1)$ over the
exceptional $\PP^1$!  How is it that Ferrari represents such a space by transition functions for an $\OO(1)\oplus\OO(-3)$ bundle?  In order to investigate this question, we first need to identify what
the change of coordinates is between the two representations for the resolved geometry.\skp

\noindent{\bf Identification with Ferrari's coordinates}\\
\noindent We begin the identification by noting that
$$\pi_1 = t, \skop \pi_2 = Z, \skop \pi_3 = Y, \skop \pi_4= X.$$
From Ferrari (3.40) \cite[page 638]{ferrari}, we can explicitly compute
\begin{eqnarray*}
\gamma &=& \frac{\pi_4}{Q(\pi_2)-\pi_1} = \gamma,\\
w_1 &=& \pi_1 + \frac{\pi_2\pi_4}{Q(\pi_2)-\pi_1} = t - \gamma^2 Y + \gamma P(t),\\
w_2 &=& \frac{(\pi_2 - P(w_1))\pi_3}{P(\pi_1)-\pi_2} = -Y + \gamma^{-1}(P(t)-P(w_1)).
\end{eqnarray*}
Using our transition functions along with Ferrari's (3.30) \cite[page 637]{ferrari}
$$\beta = \gamma^{-1}, \skop v_1 = \gamma^{-1} w_1, \skop v_2 = \gamma^3 w_2 + \gamma Q(w_1/\gamma)
+ \gamma^2 P(w_1) - \gamma w_1,$$
we complete the identifications
$$\begin{array}{ccc|ccc}
& &(\gamma,Y,t) \flecha (\gamma,w_1,w_2) & &&(\beta, Z, X) \flecha (\beta,v_1,v_2)\\
\hline 
& & & & &\\
\gamma &=& \gamma & \beta &=& \beta\\
w_1 &=& t-\gamma^2Y + \gamma P(t) & v_1 &=& Z - \beta^2 X + \beta Q(Z)\\
w_2 &=& -Y+\gamma^{-1}(P(t)-P(w_1)) & v_2 &=& X - \beta^{-1}(Q(Z)-Q(v_1)).\\
\end{array}$$

\noindent{\bf Jacobian determinants}\\
\noindent To check that these are valid coordinate changes, we compute the Jacobians in each chart.  In each case,
the straightforward identification of the $\PP^1$ coordinate means that we only need the remaining $2 \times 2$ minor
to check that the determinant does not vanish.

In the $\gamma$ chart we have
$$\begin{array}{cccccc}
\dfrac{\partial w_1}{\partial Y} &=& -\gamma^2, & \;\;\dfrac{\partial w_2}{\partial Y} &=& -1+\gamma P'(w_1),\\
&&&&&\\
\dfrac{\partial w_1}{\partial t} &=& 1 + \gamma P'(t), & \;\;\dfrac{\partial w_2}{\partial t} &=& 
\gamma^{-1}(P'(t)-P'(w_1)) - P'(w_1)P'(t).
\end{array}$$
The determinant is
\begin{eqnarray*}
\det &=& \dfrac{\partial w_1}{\partial Y}\dfrac{\partial w_2}{\partial t}-
\dfrac{\partial w_2}{\partial Y}\dfrac{\partial w_1}{\partial t} \\
&=& -\gamma (P'(t)-P'(w_1)) +\gamma^2 P'(w_1)P'(t) + (1-\gamma P'(w_1))(1 + \gamma P'(t))\\
&=& 1.
\end{eqnarray*}
In the $\beta$ chart we have
$$\begin{array}{cccccc}
\dfrac{\partial v_1}{\partial Z} &=& 1+\beta Q'(Z), & \;\;\dfrac{\partial v_2}{\partial Z} &=& 
\beta^{-1}(Q'(v_1)-Q'(Z))+Q'(v_1)Q'(Z),\\
&&&&&\\
\dfrac{\partial v_1}{\partial X} &=& -\beta^2, & \;\;\dfrac{\partial v_2}{\partial X} &=& 
1-\beta Q'(v_1).
\end{array}$$
The determinant is
\begin{eqnarray*}
\det &=& \dfrac{\partial v_1}{\partial Z}\dfrac{\partial v_2}{\partial X}-
\dfrac{\partial v_2}{\partial Z}\dfrac{\partial v_1}{\partial X} \\
&=& (1+\beta Q'(Z))(1-\beta Q'(v_1)) +\beta (Q'(v_1)-Q'(Z)) + \beta^2 Q'(v_1)Q'(Z)\\
&=& 1.
\end{eqnarray*}

\noindent This shows that we have a valid change of coordinates.\skp

\noindent {\bf Resolution of the puzzle}\\
\noindent How can it be that via a valid change of coordinates, one can go from a normal bundle
which is $\OO(-1)\oplus\OO(-1)$ to one which is $\OO(1)\oplus\OO(-3)$?  The resolution
of this apparent paradox lies in the fact that Ferrari never had an $\OO(1)\oplus\OO(-3)$
normal bundle to begin with.  One of his superpotential terms actually {\em changes the bundle structure}, and
despite appearances, his original $\Hat{\M}$ transition functions also correspond to an
 $\OO(-1)\oplus\OO(-1)$ bundle.

Consider the linearized version of the transition functions,
$$\beta = \gamma^{-1}, \skop v_1 = \gamma^{-1} w_1, \skop v_2 = \gamma^3 w_2 - \gamma w_1.$$
The $\gamma w_1$ term disrupts the $\OO(1)\oplus\OO(-3)$ structure, and suggests the following change
of coordinates
$$ \til{w_1}=w_1-\gamma^2w_2, \skop \til{v_1}=v_1+\beta^2v_2.$$
With these new coordinates we find the transition functions to be
$$\beta=1/\gamma,\skop \til{v_1}=\gamma w_2, \skop v_2 = -\gamma\til{w_1},$$
which clearly exhibits an $\OO(-1)\oplus\OO(-1)$ structure.

With an $\OO(-1)\oplus\OO(-1)$ bundle structure, $M = 0$.  We therefore can't
have any nontrivial superpotential, even if there's a geometric potential.
However, by changing coordinates Ferrari achieved a setup where $M = 2$, and
the superpotential was
$$W(x,y) = V(x)+U(y)-xy,$$
with $V'(x)=P(x)$ and $U'(y)=Q(y).$
In particular, it is the $xy$ term in the superpotential which leads to 
the change in the normal bundle structure.  (Recall that we're assuming from
the beginning that $x^2|P(x)$ and $x^2|Q(x)$.  Otherwise our $\OO(-1)\oplus\OO(-1)$
identification is invalid, as the polynomials contribute to the linearized transition functions.)

In this example, Ferrari is implicitly using the superpotential term to artificially
introduce fields $x$ and $y$ into a resolved geometry that should correspond to $M=0$ fields.  
The 2-matrix model theory he gets, therefore, should be related to the ``0-matrix model'' 
corresponding to the geometry with an $\OO(-1)\oplus\OO(-1)$ bundle.

\section{Superpotentials which change bundle structure}\label{sec:bundlechange}

What we see in Ferrari's example 2 is actually part of a more general phenomenon, which we were led
to uncover in trying to understand precisely this puzzle.  In this section we describe a family of superpotentials which change the underlying bundle structure, thus proving Proposition 1.

Consider an $\OO(n)\oplus\OO(m)$ normal bundle over a $\PP^1$ with geometric potential
$$E(\gamma,w_1) = \frac{1}{2}\gamma^k w_1^2, \skop k \in \ZZ.$$
The perturbed transition functions are
$$\beta = 1/\gamma, \skop v_1 = \gamma^{-n} w_1, \skop v_2 = \gamma^{-m}w_2+\gamma^k w_1.$$
If $(n,m,k)$ satisfies
$$ -n \leq k \leq -m,$$
we can perform the following change of coordinates:
$$\begin{array}{cccccc}
\til{w_1} &=& w_1+\gamma^{-m-k}w_2, &\til{v_1} &=& v_2,\\
\til{w_2} &=& w_2, & \til{v_2} &=& -v_1+\beta^{n+k}v_2.
\end{array}$$
Notice that
\begin{eqnarray*}
\til{v_1} &=&  \gamma^{-m}w_2+\gamma^k w_1 = \gamma^k \til{w_1},\\
\til{v_2} &=& -\gamma^{-n} w_1 + \gamma^{-n-k}(\gamma^{-m}w_2+\gamma^k w_1) = \gamma^{-n-m-k}\til{w_2},
\end{eqnarray*}
and so the new transition functions are
$$\beta = 1/\gamma, \skop \til{v_1} = \gamma^{k}\til{w_1}, \skop \til{v_2} = \gamma^{-n-m-k}\til{w_2}.$$
The geometric potential has changed our $\OO(n) \oplus \OO(m)$ bundle into
an $\OO(-k)\oplus\OO(n+m+k)$ bundle, with no superpotential.
The corresponding superpotential can be computed for $n+m = -2$, and depends on
the number of fields $M = n+1 = -m -1$:

\begin{eqnarray*}
W(x_1,...,x_M) &=& \frac{1}{2\pi i} \oint_{C_0} \gamma^{-M-1} E(\gamma, \sum_{i=1}^M x_i\gamma^{i-1})\d\gamma\\
&=& \frac{1}{2}\sum_{i=1}^M x_i x_{M-k-i+2}
\end{eqnarray*}
Note that all of these bundle-changing superpotentials are purely quadratic!
In the cases of interest, where $n+m=-2,$ the condition for the change of coordinates to be valid becomes
$$-n \leq k \leq n+2.$$
For allowed pairs $(n,k)$ we can thus get $$\OO(n)\oplus\OO(-n-2) \flecha \OO(-k)\oplus\OO(k-2)$$
by means of the perturbation.  Alternatively, we can think of these examples as ``true'' $\OO(-k)\oplus\OO(k-2)$
bundles which can be rewritten to ``look like'' $\OO(n)\oplus\OO(-n-2)$ plus a superpotential term.
\skp

\noindent {\bf RG conjecture}\\
\noindent In order to make contact with Ferrari's RG conjecture (see Introduction), 
we change notation a bit from the previous discussion:
$$M = n+1, \skop r = k-1.$$
The perturbation term in the following transition functions
$$\beta = \gamma^{-1}, \skop w_1' = \gamma^{-M+1}w_1, \skop w_2' =
\gamma^{M+1}w_2+\gamma^{r+1}w_1,$$ changes the bundle
$$\OO(M-1)\oplus\OO(-M-1) \flecha \OO(r-1)\oplus\OO(-r-1),
\skop \for \skop -M \leq r \leq M.$$ The change of coordinates:
$$\begin{array}{cccccc}
v_1 &=& w_1+\gamma^{M-r}w_2, &v_1' &=& w_2',\\
v_2 &=& w_2, &v_2' &=& -w_1'+\beta^{M+r}w_2',
\end{array}$$
yields new transition functions
$$\beta = \gamma^{-1}, \skop v_1' = \gamma^{r+1}v_1, \skop v_2' = \gamma^{1-r}v_2.$$

\noindent {\bf The superpotential}\\
\noindent The superpotential corresponding to the perturbation
$$\partial_{w_1}E(\gamma,w_1)=\gamma^{r+1}w_1$$
is given by
$$W_r(x_1,...,x_M) = \left\{\begin{array}{c}
\displaystyle{\frac{1}{2}\sum_{i=1}^{M-r}x_ix_{M-r+1-i}}, \skop
\for \skop r \geq 0,\\
\\
\displaystyle{\frac{1}{2}\sum_{i=1-r}^{M}x_ix_{M-r+1-i}}, \skop \for
\skop r \leq 0.\end{array}\right.$$
This completes the proof of Proposition 1.

Notice that the case $r=M$ is not interesting, as the bundle remains
unchanged and the superpotential vanishes.  Moreover,
the symmetry $r \mapsto -r$ in the bundle expression interchanges
two different superpotentials, but this amounts to a simple change
of coordinates.  To see this, first note that:
$$\begin{array}{cccccc}
&r \geq 0: &W_r(x_1,...,x_M) &=
&\displaystyle{\frac{1}{2}\sum_{i=1}^{M-|r|}x_ix_{M-|r|+1-i}}&\\
&r \leq 0: &W_r(x_1,...,x_M) &=
&\displaystyle{\frac{1}{2}\sum_{i=1+|r|}^{M}x_ix_{M+|r|+1-i}}& =
\displaystyle{\frac{1}{2}\sum_{i=1}^{M-|r|}x_{i+|r|}x_{M+1-i}}\\
\end{array}$$
The direction of the coordinate shift depends on the sign of $r$:
$$\begin{array}{ccccc}
&r \geq 0:& r \mapsto -r & \mathrm{is \;\; equivalent \;\; to} & x_i
\mapsto x_{i+|r|}\\
&r \leq 0:& r \mapsto -r & \mathrm{is \;\; equivalent \;\; to} & x_i
\mapsto x_{i-|r|}
\end{array}$$
i.e. $r \mapsto -r$ on the bundle side is equivalent to a simple
coordinate change for the corresponding superpotential.

We summarize the first few examples in the following table:
$$\begin{array}{c|c|c|c|c|c|c|c}
r & -3 & -2 & -1 & 0 & 1 & 2 & 3 \\
\hline 
&&&&&&&\\
M=1 & & & & \12 x_1^2 & & &\\
M=2 & & & \12 x_2^2 & x_1 x_2 & \12 x_1^2 & &\\
M=3 & & \12 x_3^2 & x_2x_3 & x_1 x_3 + \12 x_2^2 & x_1 x_2 & \12 x_1^2 & \\
M=4 & \12 x_4^2 & x_3x_4 & x_2x_4+\12 x_3^2 & x_1x_4+x_2x_3 &
x_1x_3+ \12 x_2^2 & x_1 x_2 & \12 x_1^2
\end{array}$$

\noindent{\bf The Hessian}\\
\noindent We compute the partial derivatives of our bundle-changing
superpotentials:
$$\begin{array}{cccccccccc}
&r \geq 0:& \dfrac{\partial W_r}{\partial x_j} &=& x_{M-r+1-j},&
\dfrac{\partial^2 W_r}{\partial x_k \partial x_j} &=&
\delta_{k,M-r+1-j}& \for & 1 \leq j \leq M-r.\\
&&&&&&&&&\\
&r \leq 0:& \dfrac{\partial W_r}{\partial x_j} &=& x_{M-r+1-j},&
\dfrac{\partial^2 W_r}{\partial x_k \partial x_j} &=&
\delta_{k,M-r+1-j}& \for & 1-r \leq j \leq M.\\
\end{array}$$
In each case, there is only one $k$ for every $j$ which yields a
non-zero second-partial.  This means the Hessian matrix has at most
one non-zero entry in each row and in each column.  The corank of
the Hessian is thus easy to compute, and is equal to the number of
rows (or columns) comprised entirely of zeroes.  In both the $r \geq
0$ and $r \leq 0$ cases, the corank of the Hessian is $r$ (see the
ranges for $j$ values). This is consistent with what we expect from
Ferrari's RG conjecture.

\section{Classification of perturbation terms}
 
\noindent{\bf 1--matrix model}\\
\noindent We consider the case with $M=1$ fields, which gives an $\OO\oplus\OO(-2)$ bundle with perturbation:
$$\beta=\gamma^{-1}, \skop v_1=w_1, \skop v_2=\gamma^2 w_2+ \partial_{w_1}E(\gamma,w_1).$$
The only perturbation terms which contribute to the superpotential are those
of the form
$$ \partial_{w_1}E(\gamma,w_1)=\gamma w_1^m.$$
All other terms can be added freely to the RHS of the $v_2$ transition function without affecting
the superpotential.  The linear term $\gamma w_1$ changes the bundle structure to $\OO(-1)\oplus\OO(-1)$,
which we will see corresponds to the $A_1$ case, with potential $W(x) = \dfrac{1}{2}x^2.$\skp

\noindent{\bf 2--matrix model}\\
\noindent We consider the case with $M=2$ fields, which gives an $\OO(1)\oplus\OO(-3)$ bundle with perturbation:
$$\beta=\gamma^{-1}, \skop v_1=\gamma^{-1}w_1, \skop v_2=\gamma^3 w_2+ \partial_{w_1}E(\gamma,w_1).$$
The only perturbation terms which contribute to the superpotential are those
of the form
$$ \partial_{w_1}E(\gamma,w_1)=\gamma^n w_1^m, \skop \mathrm{s.t.} \;\;\; -m < n < 3.$$
All other terms can be added freely to the RHS of the $v_2$ transition function without affecting
the superpotential.  Interchanging $x$ and $y$ in the superpotential leads to
$$\partial_{w_1}E(\gamma,w_1)=\gamma^n w_1^m \longmapsto \gamma^{3-m-n}w_1^m.$$
Again, the linear ($m=1$) perturbation terms are the ones which change the bundle structure.  We have
$w_1, \gamma w_1,$ and $\gamma^2 w_1$.  The first and third terms are related by interchanging $x$ and $y$, and
yield an $\OO \oplus \OO(-2)$ bundle with zero potential.  The $\gamma w_1$ term reduces the bundle all the way to $\OO(-1)\oplus\OO(-1)$.  The effect of these terms is summarized in Table 3.1.

\begin{table}[h]
$$\begin{array}{c|c|c|c}
\partial_{w_1}E(\gamma,w_1) & W(x,y) & \mathrm{new\;bundle} & \mathrm{type}\\
\hline
&&&\\
w_1 & \dfrac{1}{2}y^2 & \OO \oplus \OO(-2) & \Hat{A}\\
\gamma w_1 & xy & \OO(-1) \oplus \OO(-1) & A_1 \;\;(\mathrm{Ferrari's\;Example\;2})\\
\gamma^2 w_1 & \dfrac{1}{2}x^2 & \OO \oplus \OO(-2) & \Hat{A}
\end{array}$$
\label{bundle-change}
\caption{Bundle-changing perturbation terms in the 2-matrix model.}
\end{table}

\noindent{\bf Summary}\\
\noindent Table 3.2 summarizes everything we need to know about perturbation terms which
contribute to the superpotential in the one- and two-matrix model cases.  

\begin{table}[h]
$$\begin{array}{c|c|c|c}
& \mathrm{potential} & \partial_{w_1}E(\gamma,w_1) & \mathrm{bundle-changing}\\
\hline 
&&&\\
1-\mathrm{matrix} &  W(x) = \sum \dfrac{a_m}{m+1} x^{m+1} & \sum a_m \gamma w_1^m & \gamma w_1\\
&&&\\
2-\mathrm{matrix} &  W(x,y) = \sum \dfrac{a_{jk}}{j+k}{j+k \choose k} x^j y^k & 
\sum a_{jk} \gamma^{2-k} w_1^{j+k-1} & w_1, \gamma w_1, \gamma^2 w_1
\end{array}$$
\label{pert-term}
\caption{Summary of superpotential-contributing perturbation terms.}
\end{table}

\pagebreak 

\section{Resolution of the ideal sheaf}

Here we present a resolution of the ideal sheaf over the exceptional $\PP^1$ by introducing
additional coordinates $w_3$ and $v_3$ which simplify the transition functions.  In some cases,
this strategy allows us to immediately write down global holomorphic functions which help to
perform the blow-down of the $\PP^1$.  However, what we find most useful in this framework is the
identification of natural, additional $w_3, v_3$ coordinates.  In the case of the $E_7$ Intriligator--Wecht
superpotential (see Chapter~\ref{ch:length3}), these extra coordinates will help us to do the blowup and verify that the singular geometry we find is the correct one.

Starting from the transition functions
$$\beta=\gamma^{-1},\skop v_1=\gamma^{-1}w_1, \skop v_2=\gamma^{3}w_2+\partial_{w_1}E(\gamma,w_1),$$
with perturbation term
$$\partial E = \sum c_i \gamma^{2-a_i} w_1^{a_i+b_i-1}
= \sum c_i \beta^{-b_i-1}v_1^{a_i+b_i-1},$$
let
$$a = \max{a_i}, \skop b = \max{b_i}.$$
Consider the simpler transition functions
\begin{eqnarray*}
\beta &=& \gamma^{-1},\\
v_1 &=& \gamma^{-1} w_1,\\
v_2 &=& \gamma^{2-a} w_3,\\
v_3 &=& \gamma^{2-b} w_2,
\end{eqnarray*}
where the new functions $w_3$ and $v_3$ are defined by
\begin{eqnarray*}
w_3 &=& \gamma^{a+1}w_2 + \gamma^{a-2} \partial E = \gamma^{a+1}w_2 + p(\gamma,w_1)w_1,\\
v_3 &=& \beta^{b+1}v_2 - \beta^{b+1} \partial E = \beta^{b+1}v_2 - q(\beta,v_1)v_1,
\end{eqnarray*}
and the functions $p$ and $q$ are given by
\begin{eqnarray*}
p(\gamma,w_1) &=& \sum c_i \gamma^{a-a_i} w_1^{a_i+b_i-2},\\
q(\beta,v_1) &=& \sum c_i \beta^{b-b_i} v_1^{a_i+b_i-2}.
\end{eqnarray*}

In the framework of the $$\OO(1)\oplus\OO(a-2)\oplus\OO(b-2)$$ bundle, we can easily
write down global holomorphic functions as long as either $a<2$ or $b<2$.
Unfortunately this only covers the $\Hat{A}$, $\Hat{D},$ and $\Hat{E}$ cases.
However, this strategy also yields some insight in the other cases.

We have an exact sequence
$$0 \flecha \OO(-a-b) \flecha \OO(1-a)\oplus\OO(1-b)\oplus\OO(1-a-b) \flecha
\OO(2-a)\oplus\OO(-1)\oplus\OO(2-b) \flecha \OO,$$
with maps in the first chart,
$$\bRow{1}{\gamma^{a+1}}{w_1} \;\; \bMtx{w_1}{-w_3}{0}{0}{w_2}{-w_1}{-1}{p(\gamma,w_1)}{\gamma^{a+1}}
\;\; \Bcol{w_3}{w_1}{w_2},$$
and in the second chart,
$$\bRow{\beta^{b+1}}{1}{v_1} \;\; \bMtx{v_1}{-v_2}{0}{0}{v_3}{-v_1}{-\beta^{b+1}}
{q(\beta,v_1)}{1} \;\; \Bcol{v_2}{v_1}{v_3}.$$
For reference, we also record the transformation properties of the maps:
$$\bRow{\OO(b+1)}{\OO(a+1)}{\OO(1)} \;\; \bMtx{\OO(1)}{\OO(a-2)}{\OO(1+a-b)}{\OO(1+b-a)}
{\OO(b-2)}{\OO(1)}{\OO(b+1)}{\OO(a+b-2)}{\OO(a+1)}.$$
\end{onecolumn}

\begin{onecolumn}
\chapter{Blowups for surface singularities}\label{ch:node}
\def\ES{\mathrm{ES}}
\def\PT{\mathrm{PT}}
\def\eqn{\mathrm{eqn}}
\def\Rodd{\mathcal{R}_{\odd}} 
\def\Reven{\mathcal{R}_{\even}}

\noindent In this chapter we will begin by describing in general terms how to blow up a single node of a Dynkin
diagram.  This will require some sophisticated techniques, including use of the McKay correspondence in order
to identify the Weil divisor associated to each particular node.  Once the appropriate matrix factorization has
been found, the blowup is fairly straightforward.  Here we focus on blowing up nodes in $A_{n-1}$ and $D_{n+2}$ diagrams.  We will use these results in Chapter 5, where corresponding small resolutions will be obtained by deforming the matrix factorizations in the surface case.

Recall that given a polynomial $f$ defining a hypersurface, a pair of square matrices $(\phi,\psi)$ such that
$$\phi \cdot \psi = f \cdot \mathbf{1}\;\;\; \and \;\;\; \psi \cdot \phi = f \cdot \mathbf{1}$$
is called a {\em matrix factorization}\footnote{For everything you ever wanted to know about matrix factorizations, see \cite{yoshino}.} of $f$.

\section{How to blow up a single node}

\noindent The main ideas in this section come from \cite{verdier} and \cite{yoshino}.

We would like to blow up an arbitrary node in the Dynkin diagram for a simple surface singularity.  
The McKay correspondence will allow us to identify which Weil divisor corresponds to which node.
The main idea is to find the module corresponding to the representation for a particular node.  The presence of
the singularity is reflected in that this module is not locally free. Relations among the generators can be
arranged in a matrix which, together with its matrix of syzygies, gives a matrix factorization for the singular equation.  Each node yields a different (indecomposable) matrix factorization for the original equation.
Once the matrix factorization is found, the blowup is obtained by trying to make the module locally free.

Let $\Gamma \subset SU(2)$ be a finite subgroup, and $\CC^2/\Gamma$ the corresponding
singularity.  The McKay correspondence says that there is a $1-1$ correspondence
between irreducible representations of $\Gamma$ and nodes on the extended
Dynkin diagram, where the trivial representation corresponds to the added node.
The degree (dimension) of an irreducible representation $\rho$ is given by the
coefficient of that node in the largest root.  Moreover, $\rho$ and $\rhotil$
are adjacent in the diagram iff $\rhotil$ appears in the decomposition
of $\rho \otimes \rho_0$, where $\rho_0$ is the defining 2-dimensional representation
(from the induced action on $\CC^2$).

In order to find the Weil divisor for a particular node, we will look at the
corresponding representation $\rho$ of $\Gamma$.  Identifying the action of $\Gamma$
on the ring $\CC[x,y]$, we will find that the elements transforming in the $\rho$
representation form a $\CC[x,y]^{\Gamma}$--module which is not locally free.  Relations
among the generators give a matrix factorization for the original equation.  This is
then used to do the blowup.

Using matrix factorizations to blow up a node in an ADE Dynkin diagram is a natural continuation of the work of Gonzalez-Sprinberg \& Verdier \cite{verdier}.  We rederive matrix factorizations for $A_{n-1}$ and $D_{n+2}$ surface singularities in a similar, but slightly different way.
In the $A_{n-1}$ and $D_{n+2}$ cases, the finite groups $\Gamma$ are:
\begin{eqnarray*}
\ZZ_n &=& \left\langle g = \mtx{\xi}{0}{0}{\xi^{-1}} \right\rangle,
\skop \xi = e^{2\pi i/n}, \skop \xi^n = 1,\\
\BB_{4n} &=& \left\langle g = \mtx{\zeta}{0}{0}{\zeta^{-1}},
h = \mtx{0}{1}{-1}{0} \right\rangle, \skop \zeta = e^{\pi i/n},
\skop \zeta^{2n} = 1.
\end{eqnarray*}
$\Gamma \subset SU(2)$ acts on $\CC[x,y]$, and its irreducible representations
(in theory) can all be found there.  Note that $x$ and $y$ should be regarded
as coordinate functions on $\CC^2$:
$$ x = \row{1}{0}, \skop y = \row{0}{1}.$$


\noindent {\bf General recipe}\\
\noindent The general procedure for blowing up a single node in a Dynkin diagram thus consists
of the following steps:
\begin{enumerate}
\item Identify the group action for the diagram type.
\item Study the invariant theory under this action to identify the $\CC[x,y]^{\Gamma}$--module containing the irreducible representations corresponding to the chosen node.
\item Find relations among the generators of this module, arrange in a matrix.  This, together
with the matrix of syzygies, gives a matrix factorization.
\item From the matrix of relations identify the Weil divisor corresponding to the node.
\item Blow up this divisor!
\end{enumerate}
\noindent In the fully deformed cases (see Chapter 5), we will take the matrix of relations for generators of the module from the non-deformed case, and cleverly ``extend'' the entries to the full deformation.

\section{The $A_{n-1}$ story}

\subsection{Invariant theory}

\noindent $\bf \ZZ_n$ {\bf Group action}\\ 
The action of $\ZZ_n$ on $\CC[x,y]$ is characterized by
\begin{eqnarray*}
g \cdot x &=& \row{1}{0} \mtx{\xi}{0}{0}{\xi^{-1}} = \row{\xi}{0} = \xi x,\\
g \cdot y &=& \row{0}{1} \mtx{\xi}{0}{0}{\xi^{-1}} = \row{0}{\xi^{-1}} = \xi^{-1}y,
\end{eqnarray*}
from which we induce the action on monomials of higher degree:
\begin{eqnarray*}
g \cdot x^m &=& (g \cdot x)^m = \xi^m x^m,\\
g \cdot y^m &=& (g \cdot y)^m = \xi^{-m} y^m,\\
g \cdot x^k y^l &=& (g \cdot x)^k(g \cdot y)^l = \xi^{k-l} x^k y^l.
\end{eqnarray*}

\noindent $\bf \ZZ_n$ {\bf Invariant theory}\\  
Since $\xi^n = 1$, the invariant polynomials under $\ZZ_n$ are generated by:
$$x^n, \;\; xy, \;\; y^n.$$
If we define
$$X = x^n, \;\; Y = y^n, \;\; Z = xy,$$
then $X,Y$ and $Z$ are invariant under $\Gamma$.  In particular, $X, Y$ and $Z$ satisfy the
equation for an $A_{n-1}$ surface singularity:
$$XY - Z^n = 0.$$
The irreducible representations $\{\rho_j\}$ of $\ZZ_n$ are labeled by the characters
$\{\chi_j\}$, where $\chi_j(g) = \xi^j$, $j = 1...n$.  From now on we will
simply refer to $\xi^j$ as the character.  $\rho_n$ is the trivial
representation; all the irreducible representations are one-dimensional.
For instance, the defining representation $\rho_0$ splits as
$\rho_0 = \rho_1 \oplus \rho_{n-1}$.  Within the ring $\CC[x,y]$, the irreducible
representations with character $\xi^m$ are given by the
one-dimensional vector spaces spanned by the monomials
$$\xi^m: \;\;\; \{x^j y^k\} \;\;\; \st \;\; j-k \equiv m \mod n.$$
Let $\CC[x,y]^{\Gamma}$ denote the subring of invariant functions.
It is generated by $X, Y$ and $Z$, and can also be written as
$$\CC[x,y]^{\Gamma} = \CC[X,Y,Z]/(XY - Z^n).$$
By definition $\Gamma$ acts trivially on $\CC[x,y]^{\Gamma}$, so
multiplying a monomial by an element of the invariant ring will yield a new
representation with the same character. 
Consider the representation $x^j y^k$ where $j-k = m + qn.$
If $q \geq 0$ we can write
$$x^j y^k = x^{j-k}Z^k = x^m X^q Z^k \in x^m \CC[x,y]^{\Gamma}.$$
Similarly, for $q \leq -1$ we can write
$$x^j y^k = y^{k-j}Z^j =  y^{-m-qn}Z^j = y^{n-m} Y^{-q-1} Z^j
\in y^{n-m} \CC[x,y]^{\Gamma}.$$
In other words, all of the irreducible representations
with character $\xi^m$ can be found in the $\CC[x,y]^{\Gamma}$--module
$$\langle x^m, y^{n-m} \rangle \CC[x,y]^{\Gamma}.$$
Importantly, this is not a free module.  Let $s_1 = x^m,$ and $s_2 = y^{n-m}.$
Then any element in the module can be written as $a s_1 + b s_2$, where
$a,b \in \CC[x,y]^{\Gamma}$.  However, there are relations, and the above
representation is not unique (if elements in the module could be uniquely
specified by ordered pairs $(a,b)$, the module would be free).  For example,
\begin{eqnarray*}
Y s_1 &=& Z^m s_2,\\
X s_2 &=& Z^{n-m} s_1.
\end{eqnarray*}
These relations are related to a factorization of the defining equation:
$$X\cdot Y - Z^m \cdot Z^{n-m} = 0.$$
Blowing up the singularity renders the module (locally) free.

\subsection{Matrix of relations}
The relations among the generators $s_1$ and $s_2$ of the $\CC[x,y]^{\Gamma}$--module corresponding to the $m$th node
(the $\ZZ_n$ representation of character $\xi^m$) can be put into a matrix
$$\mathcal{R} =
\left(\begin{array}{cc}
Y & -Z^m \\
Z^{n-m} & -X
\end{array}\right)$$
which has a corresponding matrix of syzygies
$$\mathcal{S} = 
\left(\begin{array}{cc}
X & -Z^m \\
Z^{n-m} & -Y
\end{array}\right).$$
Together, the pair $(\mathcal{R},\mathcal{S})$ gives a matrix factorization for the $A_{n-1}$ equation:
$$\mathcal{S}\mathcal{R} = \mathcal{R}\mathcal{S} = (XY-Z^n)\mathbf{1}_{2\times 2}.$$

\subsection{Summary of results}

\noindent $\Gamma = \ZZ_{n}$, $X = \Spec \CC[x,y]^{\Gamma}$. If we let
$$X = x^n, \skop Y = y^n, \skop Z = xy,$$
then the ring of invariant polynomials can be written
$$\CC[x,y]^{\Gamma} = \CC[X,Y,Z]/(XY - Z^n),$$
which we recognize as a space with $A_{n-1}$ singularity at the origin.
The sheaf of $\OO_X$-modules we are interested in is generated by $s_1$ and $s_2$:
$$\FF = \{s_1 = x^m, s_2 = y^{n-m}\} \cdot \OO_X.$$
We can view these as providing a map
\begin{eqnarray*}
(s_1,s_2): \OO_X^{\;\oplus 2} & \flecha & \FF.\\
(f,g) & \longmapsto & f \cdot s_1 + g \cdot s_2
\end{eqnarray*}
$\FF$ corresponds to elements of $\CC[x,y]$ which transform as a
$\xi^m$ representation of $\ZZ_n$.  We find relations:
\begin{eqnarray*}
Y s_1 - Z^m s_2 &=& 0,\\
Z^{n-m} s_1 - X s_2 &=& 0.
\end{eqnarray*}
These provide a matrix factorization
$$\mathcal{R} =
\left(\begin{array}{cc}
Y & -Z^m \\
Z^{n-m} & -X
\end{array}\right), \skop
\mathcal{S} = 
\left(\begin{array}{cc}
X & -Z^m \\
Z^{n-m} & -Y
\end{array}\right),$$
and correspond to a factorization of the defining equation:
$$X \cdot Y - Z^m \cdot Z^{n-m} = 0.$$
Away from the origin the two relations are the same, and allow us to solve
for one generator in terms of the other:
$$s_2 = \frac{Y}{Z^m} s_1 = \frac{Z^{n-m}}{X} s_1.$$
This suggests that for purposes of figuring out the blow-up, it might be okay
to consider just one of the two relations.
\skp

\subsection{The blowup}

\noindent The original equation in $\CC^3$ is $$XY = Z^n.$$
We have the exact sequence:
$$\begin{array}{cccccc}
\OO_X & \maparrow{(Y,-Z^m)} & \OO_X^{\;\oplus 2} & \maparrow{(s_1, s_2)} &
\FF & \flecha 0.\\
& & & & & \\
f & \longmapsto & (fY,-fZ^m) & \longmapsto & f(Ys_1 - Z^ms_2) = 0 & \\
& & & & & \\
& & (f,g) & \longmapsto & fs_1+gs_2 &
\end{array}$$
In particular, we see that the map to $\FF$ has a one-dimensional
kernel---implying that $\FF$ has rank 1---except at
$$Y = Z^m = 0,$$
where the rank of the $\OO_X \flecha \OO_X^{\;\oplus 2}$ map drops,
the relation becomes vacuous, and the rank of $\FF$ jumps to two.
We blow up the Weil divisor
$Y = Z^m = 0$, which corresponds to the ideal
$$I = \langle Y,Z^m \rangle.$$
We introduce coordinates $[\mu, \nu]$ for $\PP^1$, and the relation
$$\mu Y = \nu Z^m.$$
We find the new relation and proper transform of the original curve in two coordinate charts:
$$\begin{array}{cc|cc}
&\mu = 1 & \nu = 1 &\\
\hline
&&&\\
& Y = \nu Z^m & \mu Y = Z^m & \rightarrow (A_{m-1})\\
(A_{n-m-1}) \leftarrow & X\nu = Z^{n-m} & X = \mu Z^{n-m} & \\
&&&\\
& (\nu, X ,Z) & (\mu,  Y , Z) &\\
\end{array}$$
The proper transform in the first chart is the equation for an
$A_{n-m-1}$ singularity, whereas the new relation in the second chart
describes an $A_{m-1}$ singularity!  We can thus see
how the blowup breaks $A_{n-1} \flecha A_{n-m-1}\;\oplus \; A_{m-1}$.
\skp


\setlength{\unitlength}{1 true in}

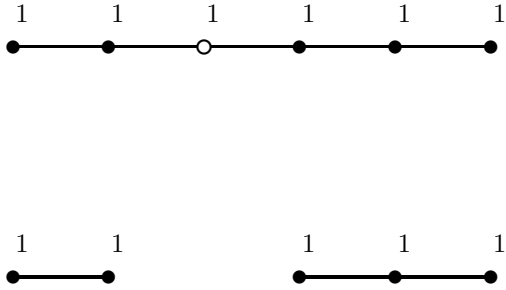
\begin{figure}[ht]

\begin{picture}(3,1)(1.65,.5)
\thicklines
\put(1.9,1){\circle*{.075}}
\put(1.9,1){\line(1,0){.5}}
\put(2.4,1){\circle*{.075}}
\put(2.4,1){\line(1,0){.4625}}
\put(2.9,1){\circle{.075}}
\put(2.9375,1){\line(1,0){.4625}}
\put(3.4,1){\circle*{.075}}
\put(3.4,1){\line(1,0){.5}}
\put(3.9,1){\circle*{.075}}
\put(3.9,1){\line(1,0){.5}}
\put(4.4,1){\circle*{.075}}
\put(1.775,1.05){\makebox(.25,.25){\footnotesize 1}}
\put(2.275,1.05){\makebox(.25,.25){\footnotesize 1}}
\put(2.775,1.05){\makebox(.25,.25){\footnotesize 1}}
\put(3.275,1.05){\makebox(.25,.25){\footnotesize 1}}
\put(3.775,1.05){\makebox(.25,.25){\footnotesize 1}}
\put(4.275,1.05){\makebox(.25,.25){\footnotesize 1}}
\end{picture}

\hspace*{\fill}

\begin{picture}(3,1)(1.65,.5)
\thicklines
\put(1.9,1){\circle*{.075}}
\put(1.9,1){\line(1,0){.5}}
\put(2.4,1){\circle*{.075}}
\put(3.4,1){\circle*{.075}}
\put(3.4,1){\line(1,0){.5}}
\put(3.9,1){\circle*{.075}}
\put(3.9,1){\line(1,0){.5}}
\put(4.4,1){\circle*{.075}}
\put(1.775,1.05){\makebox(.25,.25){\footnotesize 1}}
\put(2.275,1.05){\makebox(.25,.25){\footnotesize 1}}
\put(3.275,1.05){\makebox(.25,.25){\footnotesize 1}}
\put(3.775,1.05){\makebox(.25,.25){\footnotesize 1}}
\put(4.275,1.05){\makebox(.25,.25){\footnotesize 1}}
\end{picture}
\caption[Residual singularities from the blowup of an $A_6$ node]{An $A_6$ surface singularity splits into a pair of singularities of type $A_2$ and $A_3$ after the open node in the $A_6$ diagram is blown up.}
\end{figure}

\setlength{\unitlength}{1 pt}

\section{The $D_{n+2}$ story}

\subsection{Invariant theory}

\noindent $\bf \BB_{4n}$ {\bf Group action}\\  
\noindent The action of $\BB_{4n}$ on $\CC[x,y]$ is characterized by
\begin{eqnarray*}
g \cdot x &=& \row{1}{0} \mtx{\zeta}{0}{0}{\zeta^{-1}}
= \row{\zeta}{0} = \zeta x,\\
g \cdot y &=& \row{0}{1} \mtx{\zeta}{0}{0}{\zeta^{-1}}
= \row{0}{\zeta^{-1}} = \zeta^{-1} y,\\
h \cdot x &=& \row{1}{0} \mtx{0}{1}{-1}{0} = \row{0}{1} = y,\\
h \cdot y &=& \row{0}{1} \mtx{0}{1}{-1}{0} = \row{-1}{0} = -x.
\end{eqnarray*}

\noindent Extending the action to arbitrary monomials we have:
\begin{eqnarray*}
g \cdot x^j y^k &=& \zeta^{j-k} x^j y^k,\\
h \cdot x^j y^k &=& (-1)^k x^k y^j.
\end{eqnarray*}

\noindent $\bf \BB_{4n}$ {\bf Invariant theory}\\
\noindent Since $\zeta^{2n} = 1$, the invariant
polynomials under $\ZZ_{2n} \subset \BB_{4n}$ are generated by:
$$x^{2n}, \;\; xy, \;\; y^{2n}.$$
However, $h \cdot x^{2n} = y^{2n}$, $h \cdot y^{2n} = x^{2n}$, and
$h \cdot xy = -xy.$  If we define
$$X = \frac{i}{2}xy(x^{2n}-y^{2n}), \;\; Y = \frac{1}{2}(x^{2n}+y^{2n}),
\;\; Z = x^2y^2,$$
then $X,Y$ and $Z$ are invariant under $\Gamma = \BB_{4n}$.  
In particular, $X, Y$ and $Z$ satisfy the equation:
$$X^2 + Y^2Z - Z^{n+1} = 0,$$
which is preferred versal form for a $D_{n+2}$ singularity.

To construct irreducible representations of $\BB_{4n}$ in $\CC[x,y]$,
we start with one-dimensional representations under $\ZZ_{2n}$ and
complete the orbit under $\BB_{4n}$.  Since $h^2 = -1$, the pairs
$(x^j y^k, h \cdot x^j y^k)$ span complete and irreducible representations
of $\BB_{4n}$. Observe that when $j \neq k$ the representations are two-dimensional,
while $j = k$ yields one-dimensional representations.  Again, we label
the representations by the character evaluated on $g$.  For character
$\zeta^m + \zeta^{-m}$, $m = 1...2n$, we have:
$$\zeta^m + \zeta^{-m}: \;\;\; \{(x^j y^k, (-1)^k x^k y^j)\} \;\;\;
\st \;\; j-k \equiv m \mod 2n.$$
Let $\CC[x,y]^{\Gamma}$ denote the subring of invariant functions.
It is generated by $X, Y$ and $Z$, and can also be written as
$$\CC[x,y]^{\Gamma} = \CC[X,Y,Z]/(X^2 + Y^2Z - Z^{n+1}).$$
By definition $\Gamma$ acts trivially on $\CC[x,y]^{\Gamma}$, so
multiplying a monomial by an element of the invariant ring will yield a new
representation with the same character.
Consider the $\zeta^m + \zeta^{-m}$ representations,
$$(x^j y^k, (-1)^k x^k y^j) \;\; \where \;\; j-k = m + 2nq.$$
We split these into four cases:

\begin{eqnarray*}
q \geq 0, k \;\even &:& \;\;\;
(x^j y^k, x^k y^j) =  (x^{j-k}, y^{j-k}) Z^{k/2}\\
&& \skop \in (x^{m+2nq}, y^{m+2nq}) \CC[x,y]^{\Gamma},\\
q \geq 0, k \;\odd &:& \;\;\;
(x^j y^k, -x^k y^j) = (x^{j-k+1}y, -y^{j-k+1}x) Z^{(k-1)/2}\\
&& \skop \in (x^{m+1+2nq}y, -y^{m+1+2nq}x) \CC[x,y]^{\Gamma},\\
q \leq -1, j \;\even &:& \;\;\;
(x^j y^k, (-1)^k x^k y^j) =  (y^{k-j}, (-1)^k x^{k-j}) Z^{j/2}\\
&& \skop \in (y^{-m-2nq}, (-1)^m x^{-m-2nq}) \CC[x,y]^{\Gamma},\\
q \leq -1, j \;\odd &:& \;\;\;
(x^j y^k, (-1)^k x^k y^j) = (y^{k-j+1}x, (-1)^k x^{k-j+1}y) Z^{(j-1)/2}\\
&& \skop \in (y^{-m-2nq+1}x, (-1)^{m+1} x^{-m-2nq+1}y) \CC[x,y]^{\Gamma}.
\end{eqnarray*}

\noindent It seems at first glance that we have four infinite families of generators
for the $\CC[x,y]^{\Gamma}$-module corresponding to the $\zeta^m + \zeta^{-m}$
representations:
\begin{footnotesize}
$$s_1^{(q)} = \col{x^{m+2nq}}{y^{m+2nq}}, \;\;\;
s_2^{(q)} = \col{y^{-2nq-m}}{(-1)^m x^{-2nq-m}},$$
$$s_3^{(q)} = \col{x^{m+1+2nq}y}{-y^{m+1+2nq}x}, \;\;\;
s_4^{(q)} = \col{y^{-2nq-m+1}x}{(-1)^{m+1} x^{-2nq-m+1}y}.$$
\end{footnotesize}

\noindent Since neither $x^{2n}$ nor $y^{2n}$ are invariant functions,
it is not obvious that we don't need the generators for all values of $q$:
$q \geq 0$ for $s_1^{(q)}, s_3^{(q)}$ and $q \leq -1$ for $s_2^{(q)}, s_4^{(q)}$.
\skp

\noindent {\bf The module}\\
\noindent It seems that in the end there will only be four generators
for this module, with relations.  We conjecture is that the following four
representations (one from each of the four families) generate the module:
\begin{footnotesize}
$$s_1 = \col{x^m}{y^m}, \;\;\;
s_2 = \col{y^{2n-m}}{(-1)^m x^{2n-m}}, \;\;\;
s_3 = \col{x^{m+1}y}{-y^{m+1}x}, \;\;\;
s_4 = \col{y^{2n-m+1}x}{(-1)^{m+1} x^{2n-m+1}y}.$$
\end{footnotesize}
Note that $s_1 = s_1^{(0)}, s_2 = s_2^{(-1)}, s_3 = s_3^{(0)},$
and $s_4 = s_4^{(-1)}$.

\subsection{Matrix of relations}

\noindent {\bf Cramer's rule}\\
\noindent The module is not free, however; there are relations among the generators.
We can compute these using Cramer's rule, which states that the solution
to a system of equations
$$Av = b,$$
is given by
$$v_j = \frac{\det M_j}{\det A},$$
where $M_j$ is the matrix obtained from $A$ by replacing the $j$th column
with the vector $b$.
If we want to solve for $s_3$ and $s_4$ in terms of $s_1$ and $s_2$, our
systems of equations will look like
\begin{eqnarray*}
\row{s_1}{s_2} \col{v_1}{v_2} &=& s_3,\\
\row{s_1}{s_2} \col{v_1}{v_2} &=& s_4.
\end{eqnarray*}
In order to implement Cramer's rule, we first compute all of the pairwise
determinants $\det(s_is_j)$. By restricting to $m$ even or $m$ odd we can
also write these in terms of the invariant functions $X, Y$ and $Z$.

$$\begin{array}{c|c|c|c}
\det(s_is_j) & \any \;\; m & m \;\; \even & m \;\; \odd\\
\hline 
&&&\\
\det(s_1s_2) & (-1)^{m} x^{2n}-y^{2n} & -2iX/\sqrt{Z} & -2Y\\
\det(s_1s_3) & -2 x^{m+1} y^{m+1} & -2\sqrt{Z}Z^{m/2} & -2 Z^{(m+1)/2}\\
\det(s_1s_4) & (-1)^{m+1} x^{2n+1}y - y^{2n+1}x & -2\sqrt{Z}Y & -2iX\\
\det(s_2s_3) & (-1)^{m+1} x^{2n+1}y - y^{2n+1}x & -2\sqrt{Z}Y & -2iX\\
\det(s_2s_4) & 2(-1)^{m+1} x^{2n-m+1}y^{2n-m+1} & -2\sqrt{Z}Z^{n-m/2} & 2 Z^{n+1-(m+1)/2}\\
\det(s_3s_4) & (-1)^{m+1} x^{2n+2}y^2 + y^{2n+2}x^2 & 2i\sqrt{Z}X & 2ZY\\
\end{array}$$
For $m$ even,
\begin{eqnarray*}
b = s_3:\skop v_1 &=& \frac{\det(s_3 s_2)}{\det(s_1 s_2)} = \frac{iYZ}{X},\skop
          v_2 = \frac{\det(s_1 s_3)}{\det(s_1 s_2)} = \frac{-iZ^{m/2+1}}{X}.\\
b = s_4:\skop v_1 &=& \frac{\det(s_4 s_2)}{\det(s_1 s_2)} = \frac{iZ^{n+1-m/2}}{X},\skop
          v_2 = \frac{\det(s_1 s_4)}{\det(s_1 s_2)} = \frac{-iYZ}{X}.
\end{eqnarray*}
For $m$ odd,
\begin{eqnarray*}
b = s_3:\skop v_1 &=& \frac{\det(s_3 s_2)}{\det(s_1 s_2)} = \frac{-iX}{Y},\skop
          v_2 = \frac{\det(s_1 s_3)}{\det(s_1 s_2)} = \frac{Z^{(m+1)/2}}{Y}.\\
b = s_4:\skop v_1 &=& \frac{\det(s_4 s_2)}{\det(s_1 s_2)} = \frac{Z^{n+1-(m+1)/2}}{Y},\skop
          v_2 = \frac{\det(s_1 s_4)}{\det(s_1 s_2)} = \frac{iX}{Y}.
\end{eqnarray*}

\noindent {\bf Matrices of relations}\\
\noindent In the case of $m$ odd, we find the following matrix of relations:
$$ m \;\; \odd \skop
\left( \begin{array}{cccc}
s_1 & s_2 & s_3 & s_4\\
-iX & Z^{(m+1)/2} & -Y & 0\\
Z^{n+1-(m+1)/2} & iX & 0 & -Y\\
-YZ & 0 & -iX & Z^{(m+1)/2}\\
0 & -YZ & Z^{n+1-(m+1)/2} & iX
\end{array} \right)$$
For $m$ even, we initially find the matrix:
$$ m \;\; \even \skop
\left( \begin{array}{cccc}
s_1 & s_2 & s_3 & s_4\\
iYZ & -iZ^{m/2+1} & -X & 0\\
iZ^{n+1-m/2} & -iYZ & 0 & -X\\
-X & 0 & iY & -iZ^{m/2}\\
0 & -iX & -Z^{n-m/2} & Y
\end{array} \right)$$
Switching the top two rows with the bottom two 
$$ m \;\; \even \skop
\left( \begin{array}{cccc}
s_1 & s_2 & s_3 & s_4\\
-X & 0 & iY & -iZ^{m/2}\\
0 & -iX & -Z^{n-m/2} & Y\\
iYZ & -iZ^{m/2+1} & -X & 0\\
iZ^{n+1-m/2} & -iYZ & 0 & -X
\end{array} \right)$$
and multiplying the first and third rows by $i$, the second row by $-1$, and the fourth row
by $-i$ (we can do this because the rows just reflect relations among the generators $s_1,..,s_4$
and can thus be shuffled and scaled independently):
$$m \;\; \even \skop
\left(\begin{array}{cccc}
s_1 & s_2 & s_3 & s_4\\
-iX & 0 & -Y & Z^{m/2}\\
0& iX &Z^{n-m/2}&-Y\\
-YZ&Z^{m/2+1}&-iX&0\\
Z^{n+1-m/2}&-YZ&0&iX\\
\end{array}\right)$$
Now the $m$ even matrix matches the $m$ odd one!\skp

\noindent {\bf Matrix factorizations}\\
\noindent Now in both the even and odd cases, we can manipulate the relations to
find syzygies, modulo the defining equation $-X^2-Y^2Z+Z^{n+1}$.  There
are four syzygies which are easy to determine; these can also be 
represented as a $4\times4$ matrix which multiplies by the relations matrix
to give zero (modulo the equation).

In the case of $m$ odd, the syzygies $\mathcal{S}$ and relations $\mathcal{R}$ matrices are
$$\mathcal{S}_{\odd} = 
\left(\begin{array}{cccc}
-iX&Z^{(m+1)/2}&Y&0\\
Z^{n+1-(m+1)/2}&iX&0&Y\\
YZ&0&-iX&Z^{(m+1)/2}\\
0&YZ&Z^{n+1-(m+1)/2}&iX\\
\end{array}\right)$$
$$\mathcal{R}_{\odd} = 
\left( \begin{array}{cccc}
-iX & Z^{(m+1)/2} & -Y & 0\\
Z^{n+1-(m+1)/2} & iX & 0 & -Y\\
-YZ & 0 & -iX & Z^{(m+1)/2}\\
0 & -YZ & Z^{n+1-(m+1)/2} & iX
\end{array} \right)$$
The product of these matrices is 
$(-X^2-Y^2Z+Z^{n+1})$ times the $4\times4$
identity matrix.

In the case of $m$ even, the syzygies and relations matrices are
$$\mathcal{S}_{\even} =
\left(\begin{array}{cccc}
-iX&0&Y&Z^{m/2}\\
0&iX&Z^{n-m/2}&Y\\
YZ&Z^{m/2+1}&-iX&0\\
Z^{n+1-m/2}&YZ&0&iX\\
\end{array}\right)$$
$$\mathcal{R}_{\even} =
\left(\begin{array}{cccc}
-iX & 0 & -Y & Z^{m/2}\\
0& iX &Z^{n-m/2}&-Y\\
-YZ&Z^{m/2+1}&-iX&0\\
Z^{n+1-m/2}&-YZ&0&iX\\
\end{array}\right)$$

\noindent The product of these matrices is 
again $(-X^2-Y^2Z+Z^{n+1})$ times the $4\times4$
identity matrix.\skp

\subsection{Summary of results}

\noindent $\Gamma = \BB_{4n}$, $X = \Spec \CC[x,y]^{\Gamma}.$  If we let
$$X = \frac{i}{2}xy(x^{2n}-y^{2n}), \;\; Y = \frac{1}{2}(x^{2n}+y^{2n}),
\;\; Z = x^2y^2,$$
then the ring of invariant polynomials can be written
$$\CC[x,y]^{\Gamma} = \CC[X,Y,Z]/(X^2+Y^2Z-Z^{n+1}),$$
which we recognize as a space with $D_{n+2}$ singularity at the origin.  The
sheaf of $\OO_X$-modules of interest is:
$$\FF = \{s_1, s_2, s_3, s_4\} \cdot \OO_X,$$
where
\begin{footnotesize}
$$s_1 = \col{x^m}{y^m}, \;\;\;
s_2 = \col{y^{2n-m}}{(-1)^m x^{2n-m}}, \;\;\;
s_3 = \col{x^{m+1}y}{-y^{m+1}x}, \;\;\;
s_4 = \col{y^{2n-m+1}x}{(-1)^{m+1} x^{2n-m+1}y}.$$
\end{footnotesize}

\noindent The relations among these generators depend on whether $m$ is even or odd.
In each case we have found $4 \times 4$ matrices of relations.  As in the
$A_{n-1}$ case, we only need half of the relations to do the blowup.
In both the $m$ even and the $m$ odd cases we find that it is simplest to
choose the top two relations, which express $s_3,s_4$ in terms of $s_1,s_2$.

In the case of $m$ odd, the relations matrix is
$$\mathcal{R}_{\odd} = 
\left( \begin{array}{cccc}
s_1 & s_2 & s_3 & s_4\\
-iX & Z^{(m+1)/2} & -Y & 0\\
Z^{n+1-(m+1)/2} & iX & 0 & -Y\\
-YZ & 0 & -iX & Z^{(m+1)/2}\\
0 & -YZ & Z^{n+1-(m+1)/2} & iX
\end{array} \right)$$
In the case of $m$ even, the relations matrix is
$$\mathcal{R}_{\even} = 
\left(\begin{array}{cccc}
s_1 & s_2 & s_3 & s_4\\
-iX & 0 & -Y & Z^{m/2}\\
0& iX &Z^{n-m/2}&-Y\\
-YZ&Z^{m/2+1}&-iX&0\\
Z^{n+1-m/2}&-YZ&0&iX\\
\end{array}\right)$$

\subsection{How to use this to blow up}

\noindent Our ``relations'' matrix $\mathcal{R}$ 
is a $4\times4$ matrix which describes a map
$A^4\to A^4$, where $A$ is the coordinate ring of our space $M$ 
(a $D_{n+2}$ singularity). We have $A=\CC[X,Y,Z]/E$, where $E$ is the
equation given above.

We calculated that for the syzygies matrix $\mathcal{S}$ and the
relations matrix $\mathcal{R}$, we have
$$\mathcal{S}\mathcal{R}=E I_4.$$
This implies that in the quotient field,
$$\mathcal{R}^{-1} = \frac1E\mathcal{S}$$
but also, since $\det\mathcal{R}=E^2$, the matrix
$$E^2\mathcal{R}^{-1}=E\mathcal{S}$$
is the matrix of cofactors of $\mathcal{R}$, i.e., the matrix whose entries
are all the $3\times3$ sub-determinants of $\mathcal{R}$.  

It follows that when restricted to our space $E=0$, the matrix of cofactors
vanishes identically.  Thus, $\mathcal{R}$ generically has rank at most $2$.

We see, therefore, that $\mathcal{R}$ defines a module over $A$ which can be 
generated by $4$ elements yet has rank $2$ (unless the entries in $\mathcal{R}$
are not generic, in which case the rank could be $1$ or $0$).  To find the
smallest blowup of our space on which this module becomes locally free,
we use a variant of the blowup construction based on a Grassmannian rather
than projective space.  That is, we consider the closure of the set
$$\{(m,g)\in M \times G(2,4) \ |\ g = \text{ fiber of our module over } m \}.$$
(At points at which the rank of the cokernel of $\mathcal{R}$ is not $2$, this
definition does not apply: that is why we need the closure.)

We can explicitly describe points in $G(2,4)$ by means of Pl\"ucker 
coordinates, and as in the usual blowup, we can take these to be proportional
to certain quantities formed out of the matrix $\mathcal{R}$.  The Pl\"ucker
coordinates correspond to $2\times2$ minors of the original matrix, and we'll
end up with relations like:
$$[\alpha_1,...,\alpha_6] = [X^2-Z^{n+1}, YZ^{n+1-(m+1)/2}, iXY, ...]$$
(where the equality is evaluated in $\mathbb{P}^5$).  In addition, we need
to impose the equations satisfied by the Grassmannian $G(2,4)\subset
\mathbb{P}^5$ which in this case is a single quadratic equation.

To see the equation of $G(2,4) \subseteq \PP^5$, let $v_1,v_2,v_3,v_4$ span a 2-dimensional subspace of $\CC^2$.  Cramer's rule tells us that
$$(v_1 \wedge v_2) v_3 + v_1 (v_2 \wedge v_3) + (v_3 \wedge v_1) v_2 = 0.$$
Wedging the above equation with $v_4$ yields
$$(v_1 \wedge v_2) (v_3 \wedge v_4) + (v_1 \wedge v_4) (v_2 \wedge v_3) - 
(v_1 \wedge v_3) (v_2 \wedge v_4) = 0.$$
In terms of the Pl\"ucker coordinates
$$\gamma = v_1 \wedge v_2, \;\; \beta = v_1 \wedge v_3,\;\; \delta = v_1 \wedge v_4,\;\;
\alpha = v_2 \wedge v_3, \;\; \varepsilon = v_2 \wedge v_4,\;\; \varphi = v_3 \wedge v_4,$$
the equation satisfied by $G(2,4)$ is
\begin{eqnarray}\label{grass}
\gamma \varphi + \delta \alpha - \beta \varepsilon = 0.
\end{eqnarray}

\subsection{The blowup}

We can find additional relations among the Pl\"ucker coordinates by restricting attention to 
the 2 $\times$ 2 minors in the top two rows of the relations matrices $\Rodd$ and $\Reven$.
Together with the relation~\eqref{grass}
$$\gamma \varphi + \delta \alpha - \beta \varepsilon = 0,$$
this gives us the blowup.\skp

\noindent{\bf $\mathbf{m}$ odd}\\
\noindent The top two rows of $\Rodd$ gives
$$
\left( \begin{array}{cccc}
-iX & Z^{(m+1)/2} & -Y & 0\\
Z^{n+1-(m+1)/2} & iX & 0 & -Y
\end{array} \right)$$
and has 2 $\times$ 2 minors $[ij]$ (where $i$ and $j$ denote columns):
$$\begin{array}{ccccccccccc}
\gamma &=& [12] &=& X^2-Z^{n+1}, &\skop& \varphi &=& [34] &=& Y^2,\\
\beta &=& [13] &=& YZ^{n+1-(m+1)/2}, &\skop& \varepsilon &=& [24] &=& -YZ^{(m+1)/2},\\
\delta &=& [14] &=& iXY, &\skop& \alpha &=& [23] &=& iXY.
\end{array}$$
We find additional Pl\"ucker coordinate relations
\begin{eqnarray*}
\delta &=& \alpha,\\
\beta &=& -Z^{n-m}\varepsilon,\\
\gamma &\equiv& -Z\varphi \mod \mathrm{equation},
\end{eqnarray*}
where the last relation holds modulo the equation for the singular surface
$$X^2 + Y^2Z - Z^{n+1} = 0.$$
\noindent{\bf $\mathbf{m}$ even}\\
\noindent The top two rows of $\Reven$ gives
$$
\left(\begin{array}{cccc}
-iX & 0 & -Y & Z^{m/2}\\
0& iX &Z^{n-m/2}&-Y
\end{array}\right)$$
and has 2 $\times$ 2 minors $[ij]$ (where $i$ and $j$ denote columns):
$$\begin{array}{ccccccccccc}
\gamma &=& [12] &=& X^2, &\skop& \varphi &=& [34] &=& Y^2-Z^n,\\
\beta &=& [13] &=& -iXZ^{n-m/2}, &\skop& \varepsilon &=& [24] &=& -iXZ^{m/2},\\
\delta &=& [14] &=& iXY, &\skop& \alpha &=& [23] &=& iXY.
\end{array}$$
Again we find Pl\"ucker coordinate relations
\begin{eqnarray*}
\delta &=& \alpha,\\
\beta &=& Z^{n-m}\varepsilon,\\
\gamma &\equiv& -Z\varphi \mod \mathrm{equation}.
\end{eqnarray*}

\noindent {\bf Identification of residual singularities}\\
We can combine the two cases by writing the Pl\"ucker coordinate relations as
\begin{eqnarray*}
\delta &=& \alpha,\\
\beta &=& (-1)^m Z^{n-m}\varepsilon,\\
\gamma &\equiv& -Z\varphi \mod \mathrm{equation},
\end{eqnarray*}
The relation for the grassmanian~\eqref{grass} becomes
$$\alpha^2 - \varphi^2Z + (-1)^{m+1}\varepsilon^2Z^{n-m}=0.$$
The interesting charts for the blowup are $\varphi=1$ and $\varepsilon=1$.\skp

\noindent{ \bf $\mathbf{m}$ odd}\\
\noindent ($\varphi=1$)  In this chart we have 
$$\varepsilon = \varphi\: \varepsilon_1, \;\; \alpha = \varphi\: \alpha_1, \;\; 
\alpha_1^2-Z+\varepsilon_1^2Z^{n-m}=0.$$
To see what kind of singularity this gives, define $G(Z,\alpha_1)$ by
$$Z^{(m+1)/2} - \alpha_1^{m+1} = (Z-\alpha_1^2)G(Z,\alpha_1),$$
and notice that 
$$-Z^{(m+1)/2}=\varepsilon_1Y, \skop \and \skop Z-\alpha_1^2 = \varepsilon_1^2Z^{n-m}.$$
We are thus left with
$$-\alpha_1^{m+1} = \varepsilon_1(Y+\varepsilon_1 Z^{n-m}G(Z,\alpha_1)) \stackrel{\mathrm{def}}{=} \varepsilon_1\til{Y},$$
which is the equation for an $A_{m}$ singularity.\skp

\noindent ($\varepsilon=1$)  In this chart 
$$\alpha = \varepsilon \alpha_2, \;\; \varphi = \varepsilon \varphi_2, \;\; \alpha_2^2 - \varphi_2^2Z + Z^{n-m} = 0,$$
which we immediately recognize as a $D_{n-m+1}$ singularity.\skp

\noindent{\bf $\mathbf{m}$ even}\\
\noindent ($\varphi=1$)  In this chart we have 
$$\varepsilon = \varphi\: \varepsilon_1, \;\; \alpha = \varphi\: \alpha_1, \;\; 
\alpha_1^2-Z-\varepsilon_1^2Z^{n-m}=0.$$
To see what kind of singularity this gives, define $G(Z,\alpha_1)$ by
$$Z^{m/2} - \alpha_1^{m} = (Z-\alpha_1^2)G(Z,\alpha_1),$$
and notice that 
$$-Z^{m/2}\alpha_1=\varepsilon_1Y, \skop \and \skop Z-\alpha_1^2 = -\varepsilon_1^2Z^{n-m}.$$
We are thus left with
$$-\alpha_1^{m+1} = \varepsilon_1(Y-\alpha_1\varepsilon_1 Z^{n-m}G(Z,\alpha_1)) \stackrel{\mathrm{def}}{=} \varepsilon_1\til{Y},$$
which is the equation for an $A_{m}$ singularity.\skp

\noindent ($\varepsilon=1$)  In this chart 
$$\alpha = \varepsilon \alpha_2, \;\; \varphi = \varepsilon \varphi_2, \;\; \alpha_2^2 - \varphi_2^2Z - Z^{n-m} = 0,$$
which we immediately recognize as a $D_{n-m+1}$ singularity.

We have thus seen that in both the odd and even cases, the blowup breaks the $D_{n+2}$ surface singularity into two ``lower order'' singularities
$$D_{n+2} \flecha D_{n-m+1} \oplus A_m.$$
This is exactly what remains in the diagram once the blown up node has been removed!

\setlength{\unitlength}{1 true in}

\begin{figure}[ht]

\begin{picture}(3,1)(1.65,.5)
\thicklines
\put(1.9,1){\circle*{.075}}
\put(1.9,1){\line(1,0){.5}}
\put(2.4,1){\circle*{.075}}
\put(2.4,1){\line(1,0){.4625}}
\put(2.9,1){\circle{.075}}
\put(2.9375,1){\line(1,0){.4625}}
\put(3.4,1){\circle*{.075}}
\put(3.4,1){\line(1,0){.5}}
\put(3.9,1){\circle*{.075}}
\put(3.9,1){\line(1,0){.5}}
\put(4.4,1){\circle*{.075}}
\put(4.4,1){\line(3,4){.3}}
\put(4.7,1.4){\circle*{.075}}
\put(4.4,1){\line(3,-4){.3}}
\put(4.7,.6){\circle*{.075}}
\put(1.775,1.05){\makebox(.25,.25){\footnotesize 1}}
\put(2.275,1.05){\makebox(.25,.25){\footnotesize 2}}
\put(2.775,1.05){\makebox(.25,.25){\footnotesize 2}}
\put(3.275,1.05){\makebox(.25,.25){\footnotesize 2}}
\put(3.775,1.05){\makebox(.25,.25){\footnotesize 2}}
\put(4.2,1.05){\makebox(.25,.25){\footnotesize 2}}
\put(4.575,1.45){\makebox(.25,.25){\footnotesize 1}}
\put(4.575,.3){\makebox(.25,.25){\footnotesize 1}}
\end{picture}

\bigskip
\bigskip

\begin{picture}(3,1.5)(1.65,.25)
\thicklines
\put(1.9,1){\circle*{.075}}
\put(1.9,1){\line(1,0){.5}}
\put(2.4,1){\circle*{.075}}
\put(3.4,1){\circle*{.075}}
\put(3.4,1){\line(1,0){.5}}
\put(3.9,1){\circle*{.075}}
\put(3.9,1){\line(1,0){.5}}
\put(4.4,1){\circle*{.075}}
\put(4.4,1){\line(3,4){.3}}
\put(4.7,1.4){\circle*{.075}}
\put(4.4,1){\line(3,-4){.3}}
\put(4.7,.6){\circle*{.075}}
\put(1.775,1.05){\makebox(.25,.25){\footnotesize 1}}
\put(2.275,1.05){\makebox(.25,.25){\footnotesize 2}}
\put(3.275,1.05){\makebox(.25,.25){\footnotesize 2}}
\put(3.775,1.05){\makebox(.25,.25){\footnotesize 2}}
\put(4.2,1.05){\makebox(.25,.25){\footnotesize 2}}
\put(4.575,1.45){\makebox(.25,.25){\footnotesize 1}}
\put(4.575,.3){\makebox(.25,.25){\footnotesize 1}}
\end{picture}

\label{fig:DynkinD}
\caption[Residual singularities from the blowup of a $D_8$ node]{The blow up of the open length 2 node in the top Dynkin diagram splits the $D_8$ surface singularity into two lower-order surface singularities: $A_2$ and $D_5$.}
\end{figure}
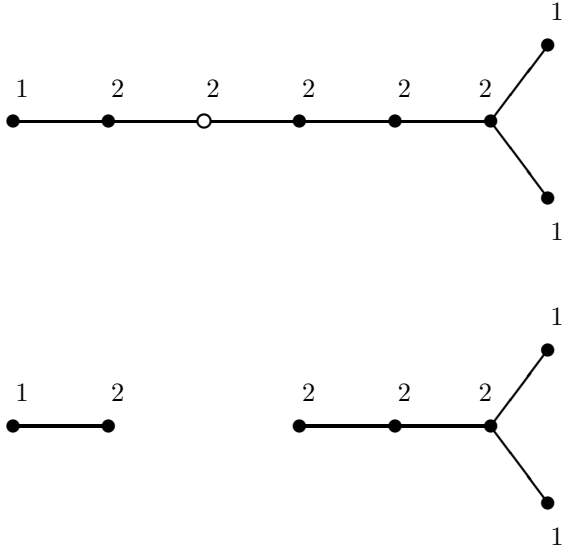

\setlength{\unitlength}{1 pt}
\end{onecolumn}

\begin{onecolumn}
\chapter{Small resolutions}\label{ch:smallres}
\def\ES{\mathrm{ES}}
\def\PT{\mathrm{PT}}
\def\eqn{\mathrm{eqn}}

\section{``Bottom-up'' approach}

\noindent In light of the algebraic geometry results of Laufer, Katz, and Morrison \cite{laufer, morrison} (see Introduction for more details), we know that there are only 6 types of geometries with isolated singularities which fit into Ferrari's framework.  These singular threefolds may have hyperplane sections pertaining to any of the ADE Dynkin diagrams, but the {\em generic hyperplane section} is determined by the {\em length} of the singular point.  The small resolution can be obtained by blowing up any node of the same length in any of the diagrams corresponding to a hyperplane section, even if it's not the generic one.

In the following ``bottom up'' approach, we will tackle the problem
of starting from one of the known singular threefolds, blowing up an arbitrary node in the associated $A_n$ or $D_{n+2}$ Dynkin diagram (this corresponds to a hyperplane section which is not generic for $n>1$), and then identifying the associated matrix model superpotential.  In this way we hope to classify all of the matrix models which can be engineered by Calabi-Yau geometries following Ferrari's construction.  The basic outline of the strategy is as follows:
\begin{itemize}
\item Start with a singular Calabi-Yau $\M_0$ with isolated Gorenstein singularity $p$ which admits an irreducible small resolution.
\item Blow up a node of a corresponding Dynkin diagram to get the small resolution $\Hat{\M}$.
\item Write the smooth space $\Hat{\M}$ in terms of transition functions in two charts over exceptional $\PP^1$.
\item Identify superpotential $W(x_1,...,x_M)$!
\end{itemize}
\noindent We find that explicitly performing the small resolution is already a challenge.  Already in the
length 2 case ($D_{n+2}$) this is highly non-trivial.  We are able to find the small resolutions by deforming
the matrix factorizations from Chapter~\ref{ch:node}.  Thus, in this section we prove Theorem 1:
\begin{mytheorem2} For Gorenstein threefold singularities of length 1 and length 2, with hyperplane section an $A_n$ or $D_{n+2}$ surface singularity, the small resolution is obtained by deforming the matrix factorization for a node of the same length in the corresponding Dynkin diagram.
\end{mytheorem2}
Perhaps surprisingly, once we are able to find the small resolution it is not easy to identify the corresponding
superpotential!  Already in the length 2 case, we find the task of identifying the right coordinates for
a description of the blowup commensurate with Ferrari's framework to be too difficult without further help.  In
the length 1 case, we can identify transition functions but only when we restrict ourselves to the $A_1$ case.

\section{The $A_{n-1}$ story}

Here we discuss the length 1 case.  The invariant theory can be found in \cite[pages 467--469]{morrison}.
The small resolution via deformation of the $A_{n-1}$ matrix factorization proves the first part of Theorem 1.

\subsection{Invariant theory}

We start with the equation for an $A_{n-1}$ singularity in preferred versal form,
$$\Phi_{A_{n-1}}: -XY + Z^n + \sum_{i=2}^{n} \alpha_i Z^{n-i}.$$
Introducing the distinguished polynomial
$$f_{n-1}(U,t) = \prod_{t=1}^{n}(U+t_i), \skop \sum_{i}t_i = 0,$$
we can rewrite the equation as
$$\Phi_{A_{n-1}}: -XY + f_{n-1}(Z,t).$$
If we are to break 
$$A_{n-1} \flecha A_{n-m-1} \oplus A_{m-1},$$
we expect the original deformed equation to factor in a way which preserves the residual
symmetries.  Rewriting $f_{n-1}$ as\footnote{Note that primes do not indicate derivatives!}
$$f_{n-1}(U,t) = \prod_{t=1}^{n-m}(U+t_i)\prod_{t=1}^{m}(U+t_i) = f_{n-m-1}'(U,t)f_{m-1}''(U,t),$$
we find the factorization
\begin{equation}\label{eq:A1}
\Phi_{A_{n-1}}: -XY + f_{n-m-1}'(Z,t)f_{m-1}''(Z,t).
\end{equation}

\subsection{Matrix factorization with deformation}

In order to ``guess'' a deformation for the matrix factorization from Chapter 4
$$\mathcal{R} =
\left(\begin{array}{cc}
Y & -Z^m \\
Z^{n-m} & -X
\end{array}\right), \skop
\mathcal{S} = 
\left(\begin{array}{cc}
X & -Z^m \\
Z^{n-m} & -Y
\end{array}\right),$$
first note that in the absence of deformation ($t=0$) the allowed invariant polynomials become
$$f_{n-m-1}'(Z) = Z^{n-m}, \skop f_{m-1}''(Z) = Z^m.$$
This leads us naturally to try the deformed matrix factorization
$$\mathcal{R_{\mathrm{def}}} =
\left(\begin{array}{cc}
Y & -f_{m-1}''(Z,t) \\
f_{n-m-1}'(Z,t) & -X
\end{array}\right), \skop
\mathcal{S_{\mathrm{def}}} = 
\left(\begin{array}{cc}
X & -f_{m-1}''(Z,t) \\
f_{n-m-1}'(Z,t) & -Y
\end{array}\right).$$
We find that this is a valid matrix factorization for the deformed equation~\eqref{eq:A1}, as
$$\mathcal{R_{\mathrm{def}}}\mathcal{S_{\mathrm{def}}} = \mathcal{S_{\mathrm{def}}}\mathcal{R_{\mathrm{def}}}
= (XY - f_{n-m-1}'(Z,t)f_{m-1}''(Z,t))\mathbf{1}_{2\times 2}.$$

\subsection{The blowup}

We blow up the Weil divisor
$Y = f_{m-1}''(Z,t) = 0$, which corresponds to the ideal
$$I = \langle Y,f_{m-1}''(Z,t)\rangle.$$
We introduce coordinates $[\mu, \nu]$ for $\PP^1$, and the relation
$$\mu Y = \nu f_{m-1}''(Z,t).$$
We find the new relation and proper transform of the original curve in two coordinate charts:
$$\begin{array}{cc|cc}
&\mu = 1 & \nu = 1 &\\
\hline
&&&\\
& Y = \nu f_{m-1}''(Z,t) & \mu Y = f_{m-1}''(Z,t) & \rightarrow (A_{m-1})\\
(A_{n-m-1}) \leftarrow & X\nu = f_{n-m-1}'(Z,t) & X = \mu f_{n-m-1}'(Z,t) & \\
&&&\\
& (\nu, X ,Z, t) & (\mu,  Y , Z, t) &\\
\end{array}$$
The proper transform in the first chart is the equation for an
$A_{n-m-1}$ singularity, whereas the new relation in the second chart
describes an $A_{m-1}$ singularity!  We can thus see
how (restricted to a hyperplane) the blowup breaks $A_{n-1} \flecha A_{n-m-1}\;\oplus \; A_{m-1}$.

\setlength{\unitlength}{1 true in}

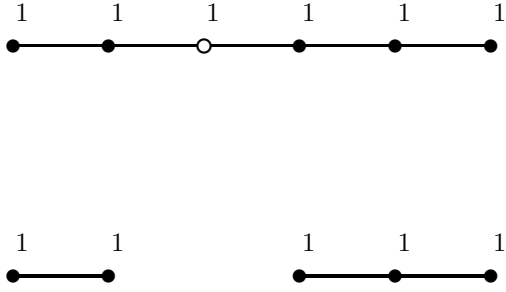
\begin{figure}[ht]

\begin{picture}(3,1)(1.65,.5)
\thicklines
\put(1.9,1){\circle*{.075}}
\put(1.9,1){\line(1,0){.5}}
\put(2.4,1){\circle*{.075}}
\put(2.4,1){\line(1,0){.4625}}
\put(2.9,1){\circle{.075}}
\put(2.9375,1){\line(1,0){.4625}}
\put(3.4,1){\circle*{.075}}
\put(3.4,1){\line(1,0){.5}}
\put(3.9,1){\circle*{.075}}
\put(3.9,1){\line(1,0){.5}}
\put(4.4,1){\circle*{.075}}
\put(1.775,1.05){\makebox(.25,.25){\footnotesize 1}}
\put(2.275,1.05){\makebox(.25,.25){\footnotesize 1}}
\put(2.775,1.05){\makebox(.25,.25){\footnotesize 1}}
\put(3.275,1.05){\makebox(.25,.25){\footnotesize 1}}
\put(3.775,1.05){\makebox(.25,.25){\footnotesize 1}}
\put(4.275,1.05){\makebox(.25,.25){\footnotesize 1}}
\end{picture}

\hspace*{\fill}

\begin{picture}(3,1)(1.65,.5)
\thicklines
\put(1.9,1){\circle*{.075}}
\put(1.9,1){\line(1,0){.5}}
\put(2.4,1){\circle*{.075}}
\put(3.4,1){\circle*{.075}}
\put(3.4,1){\line(1,0){.5}}
\put(3.9,1){\circle*{.075}}
\put(3.9,1){\line(1,0){.5}}
\put(4.4,1){\circle*{.075}}
\put(1.775,1.05){\makebox(.25,.25){\footnotesize 1}}
\put(2.275,1.05){\makebox(.25,.25){\footnotesize 1}}
\put(3.275,1.05){\makebox(.25,.25){\footnotesize 1}}
\put(3.775,1.05){\makebox(.25,.25){\footnotesize 1}}
\put(4.275,1.05){\makebox(.25,.25){\footnotesize 1}}
\end{picture}
\caption[Small resolution via blowup of a length 1 node]{An $A_6$ surface singularity splits into a pair of singularities of type $A_2$ and $A_3$ after the open node in the $A_6$ diagram is blown up.  The length 1 
singularity in the threefold can be thought of as a deformation of this surface singularity, which is in turn
given by a (non-generic) hyperplane section of the total space.  As far as the hypersurface is concerned, there are still remaining residual singularities after the blowup.  In contrast, the total space of the threefold is
smooth upon blowing up a single length 1 node in the diagram.}
\end{figure}

\setlength{\unitlength}{1 pt}

\subsection{Identification of superpotential in $A_1$ case}

We now focus on threefold singularities which are deformations of the case $A_1$ ($n=2,m=1$).  By considering
higher order deformations $t^k$ for $k>1$, we recover superpotentials for all $A_k$.\skp

\noindent {\bf Higher order deformations}\\
\noindent We start with the equation for 
an $A_1$ singularity with higher order deformations:
$$XY = Z^2-t^{2k} = (Z+t^k)(Z-t^k),$$
for $k \geq 1.$ We blow up the ideal $I = \langle Y,Z+t^k \rangle$,
corresponding to the Weil divisor $$Y = Z+t^k = 0.$$ We introduce
coordinates $[\gamma,\beta]$ for the $\PP^1$, and the relation
$$\beta Y = \gamma (Z+t^k).$$
The blowup is described in two charts:
$$\begin{array}{c|c}
&\\
\gamma\; \mathrm{chart}\;\;\;(\gamma,X,Z,t) &
\beta \; \mathrm{chart}\;\;\;(\beta,X,Y,t)\\
&\\
Y = \gamma (Z+t^k) & Z+t^k = \beta Y \\
&\\
H = \gamma X - Z+t^k = 0 & F = X - \beta^2Y +2\beta t^k = 0\\
&
\end{array}$$
The transition functions
$$\beta = \gamma^{-1},\skop X = X, \skop Y = \gamma (Z+t^k), \skop F = \gamma^{-1} H,$$
suggest that the exact sequence of normal bundles
$$0 \flecha N_{C/T} \flecha N_{C/F} \flecha N_{T/F} \flecha 0$$
for the curve $C$ within the threefold $T$, inside the fourfold
$F$, is given by\skp
$$\begin{array}{ccccccccc}
\gamma \; \mathrm{chart} & & &  & (t,X,Z) & \Bcol{0}{\gamma}{-1} &
H_0 = \gamma X - Z & &\\
& & & & & & & &\\
0 & \flecha & N_{C/T}& \maparrow{?} & \OO \oplus \OO \oplus  \OO(-1) & \flecha & \OO(1) & \flecha & 0\\
& & & & & & & &\\
\beta \; \mathrm{chart} & & &  & (t,X,Y) & \Bcol{0}{1}{-\beta^2} &
F_0 = X - \beta^2Y & &
\end{array}$$
where $H_0$ and $F_0$ are linearized versions of $H$ and $F$.
In homogeneous coordinates, the $N_{C/F} \flecha N_{T/F}$ map
becomes
$$\begin{array}{ccc}
&\Bcol{0}{\gamma}{-\beta^2}&\\
\OO \oplus \OO \oplus \OO(-1)&\flecha& \OO(1). \end{array}$$

\noindent We easily identify the following elements in the kernel:
$$[1, 0, 0], \skop [0, \beta^2, \gamma].$$
Since the corresponding matrix for the map $N_{C/T} \flecha
N_{C/F}$ has type
$$\bmtxtth{1}{0}{0}{0}{\beta^2}{\gamma} =
\bmtxtth{\OO}{*}{*}{*}{\OO(2)}{\OO(1)},$$
we can conclude that $$N_{C/T} = \OO \oplus \OO(-2).$$

\noindent The exact sequence with all the maps in each chart is thus,
\begin{footnotesize}
$$\begin{array}{ccccccccc}
\gamma \; \mathrm{chart} & &(t,r) &
\bmtxtth{1}{0}{0}{0}{1}{\gamma} & (t,X,Z) & \Bcol{0}{\gamma}{-1} &
H_0 = \gamma X - Z & &\\
& & & & & & & &\\
0 & \flecha & \OO\oplus\OO(-2) & \flecha & \OO \oplus \OO \oplus  \OO(-1) & \flecha & \OO(1) & \flecha & 0.\\
& & & & & & & &\\
\beta \; \mathrm{chart} & & (t,r') &
\bmtxtth{1}{0}{0}{0}{\beta^2}{1} & (t,X,Y) & \Bcol{0}{1}{-\beta^2}
& F_0 = X - \beta^2Y & &
\end{array}$$
\end{footnotesize}
\noindent In particular, the maps
\begin{eqnarray*}
(t,r) &\longmapsto &(t, r, \gamma r),\\
(t,r') &\longmapsto & (t, \beta^2 r', r'),
\end{eqnarray*}
imply that
$$ r = X, \skop r' = Y.$$
Using the full nonlinear transition functions, we see that
$$ Y = \gamma(Z+t^k) = \gamma(\gamma X+2t^k)= \gamma^2X+2\gamma
t^k.$$ In other words, $r,r'$ have transition function
$$ r' = \gamma^2 r +2\gamma t^k,$$
which has a superpotential correction to the $\OO(-2)$ structure.\skp

\noindent {\bf The superpotential}\\
\noindent All together, the $\OO\oplus\OO(-2)$ transition functions are
$$\beta = \gamma^{-1}, \skop t=t, \skop  r' = \gamma^2 r +
2\gamma t^k.$$ Comparing with Ferrari's setup,
$$\beta = \gamma^{-1}, \skop s'=\gamma^{-n}s, \skop
r'=\gamma^{-m}r+\partial_{s}E(\gamma,s),$$ we find that $n=0$, $m
= -2$, and $M=1$. The geometric potential is given by
$$E(\gamma,t) = \gamma\dfrac{2}{k+1} t^{k+1}.$$
The superpotential for our deformed $A_1$ case is thus,
\begin{eqnarray*}
W(x) &=& \dfrac{1}{2\pi i}\oint_{C_0}
\gamma^{-2}E(\gamma,x)\d\gamma\\
&=& \dfrac{2}{k+1} x^{k+1}.
\end{eqnarray*}
Note that the normal bundle to the $\PP^1$ is $\OO \oplus \OO(-2)$ for $k > 1$, but for $k=1$ we have a
bundle-changing superpotential, and the normal bundle is $\OO(-1) \oplus \OO(-1)$.  In particular, the different types of $A_1$ deformations yield all of the $A_k$ geometries and superpotentials.\skp

\noindent {\bf Ferrari's example 1}\\
\noindent We recover Ferrari's example if we let $$t^k \longmapsto
\dfrac{1}{2}P(t),$$ for some polynomial $P(t)$. The original
deformed equation becomes
$$XY = Z^2 - P(t)^2/4,$$
and the $r,r'$ transition function is
$$r' = \gamma^2r + \gamma P(t).$$
If $\til{W}$ is a polynomial such that $\til{W}'(t)=P(t)$, then
the geometric potential is
$$E(\gamma,t) = \gamma\til{W}(t),$$
and the superpotential is given by
\begin{eqnarray*}
W(x) &=& \dfrac{1}{2\pi i}\oint_{C_0}
\gamma^{-2}E(\gamma,x)\d\gamma\\
&=& \til{W}(x).
\end{eqnarray*}

\section{The $D_{n+2}$ story}

Now we discuss the length 2 case.  The invariant theory and the Tyurina blowup (a length 1 resolution!) can be found in \cite[pages 469--471]{morrison}.  The small resolution is obtained by deforming the matrix factorization for a length 2 node in the $D_{n+2}$ diagram, and thus completes the proof of Theorem 1.

\subsection{Invariant theory}

\noindent We start with the equation for a $D_{n+2}$ singularity in preferred versal form,
$$\Phi_{D_{n+2}}: X^2 + Y^2 Z - Z^{n+1} - \sum_{i=1}^{n+1} \delta_{2i}Z^{n+1-i} + 2 \gamma_{n+2} Y.$$
Let $t_1,...,t_{n+2}$ be coordinates on $\CC^{n+2}$, and $\WW(D_{n+2}) = (\ZZ/2)^{n+1} \ltimes
S_{n+2}$ the Weyl group of the associated lie algebra.
Then the above coefficients are invariant polynomials in $\CC^{n+2}/\WW(D_{n+2})$:
\begin{eqnarray*}
\gamma_{n+2} &:=& t_1\cdots t_{n+2} = s_{n+2},\\
\delta_{2i} &:=& \mathrm{the}\;i\mathrm{th}\;\mathrm{elementary\; symmetric\; function\; of}\;t_1^2,...,t_{n+2}^2.
\end{eqnarray*}
Define $\delta$ by
$$\delta^2 \equiv \delta_{2n+4} = \gamma_{n+2}^2.$$
The functions $\delta,\delta_2,...,\delta_{2n+2}$ generate $\CC[V]^{\WW}$.\skp

\noindent {\bf Factorization}\\
\noindent In order to factor the above $D_{n+2}$ equation, we introduce distinguished
polynomials:
$$f(U,t) = \prod_{i=1}^{n+2} (U+t_i) = U^{n+2}+ U\sum_{i=1}^{n+1} s_iU^{n+1-i}+\delta$$
$$g(Z,t) = \prod_{i=1}^{n+2} (Z+t_i^2) = Z^{n+2}+Z\sum_{i=1}^{n+1} \delta_{2i}Z^{n+1-i}+\delta^2.$$
Note that $f(0,t) = \delta$.  We can now rewrite the $D_{n+2}$ equation as
\begin{equation}\label{eq:D1}
\Phi_{D_{n+2}}: X^2 + Y^2Z - 1/Z(g(Z,t) - \delta^2) + 2\delta Y.
\end{equation}
Decomposing $f(U,t)$ into odd and even parts,
\begin{eqnarray*}
f(U,t) &=& UP(-U^2,t)+Q(-U^2,t),\\
g(-U^2,t) &=& f(U,t)f(-U,t) = (Q + UP)(Q - UP) = Q^2 - U^2P^2,\\
g(Z,t) &=& ZP(Z,t)^2 + Q(Z,t)^2.
\end{eqnarray*}
Finally, define $h(Z,t)$, $S(Z,t)$ and $G(Z,U,t)$ by
\begin{eqnarray} 
h(Z,t) &=& \frac{1}{Z}(g(Z,t) - \delta^2),\\
S(Z,t) &=& \frac{1}{Z}(Q(Z,t)- \delta),\\
UP(Z,t) + Q(Z,t) &=& (Z+U^2)G(Z,U,t) + f(U,t).\label{eq:G-defn}
\end{eqnarray}
Again we rewrite the $D_{n+2}$ equation~\eqref{eq:D1} as
\begin{eqnarray*}
\Phi_{D_{n+2}}: & & X^2 + Y^2Z - 1/Z(ZP(Z,t)^2+Q(Z,t)^2 - \delta^2) + 2\delta Y\\
&=& X^2 + Y^2Z - 1/Z(ZP(Z,t)^2 + (ZS(Z,t)+\delta)^2 - \delta^2) + 2\delta Y\\
&=& X^2 + Y^2Z - 1/Z(ZP(Z,t)^2 + ZS(Z,t)(ZS(Z,t)+2\delta)) + 2\delta Y\\
&=& X^2 + Y^2Z - P^2 - S(ZS+2\delta) + 2\delta Y.
\end{eqnarray*}

\subsection{The Tyurina blowup}

\noindent Tyurina's blow up for a $D_{n+2}$ surface singularity \cite{morrison,tyurina} results from
a very clever factorization of the deformed equation.\skp  

\noindent {\bf The miracle}\\ 
\noindent Our $D_{n+2}$ equation now factors!
$$\Phi_{D_{n+2}}: (X+P)(X-P)+(Y-S)(YZ+SZ+2\delta).$$
This is the factorization needed for the Tyurina blowup.
\skp

\noindent {\bf The blowup}\\
\noindent Our refactored equation is
$$\Phi_{D_{n+2}}: (X+P)(X-P)+(Y-S)(YZ+SZ+2\delta).$$
We blow up the top node in the Dynkin Diagram for $D_{n+2}$,
correspoding to the Weil divisor $X+P = Y-S = 0$.  This leaves us with an $A_{n+1}$ singularity,
because the remaining diagram has $n$ connected nodes.

The blow up yields
$$(X-P) = (Y-S)U, \skop \tilde{\Phi}_{D_{n+2}}: (X+P)U + (YZ+SZ+2\delta).$$
Since $X+P = (Y-S)U + 2P$, and $Q = SZ+\delta$, we can rewrite the $\tilde{\Phi}_{D_{n+2}}$ equation as
$$\tilde{\Phi}_{D_{n+2}}: (Y-S)U^2+2UP+(Y-S)Z+2(SZ+\delta) = (Y-S)(Z+U^2)+2(UP+Q).$$
Using the definition of $G$ we find
$$\tilde{\Phi}_{D_{n+2}}: (Y-S)(Z+U^2)+2(Z+U^2)G+2f(U,t) = (Z+U^2)(Y-S+2G)+2f(U,t).$$
Defining
\begin{eqnarray*}
\til{Z} &=& Z+U^2,\\
\til{Y} &=& Y-S+2G,
\end{eqnarray*}
we see that the proper transform $\tilde{\Phi}_{D_{n+2}}$ is of the form
$$\til{Y} \til{Z} = -2f_{n+2}(U,t),$$
which corresponds to an $A_{n+1}$ singularity.

Unfortunately, we see that this corresponds to blowing up a length 1
node in the associated Dynkin diagram, which will {\bf not} yield a small resolution 
for the length 2 threefold singularity type.  In the next section, we show how to blow
up length 2 node in the non-deformed case (i.e. for a surface singularity).  The results
will helps us to later do the blow up in the fully deformed case.

\setlength{\unitlength}{1 true in}

\begin{figure}[ht]

\begin{picture}(3,1)(1.65,.5)
\thicklines
\put(1.9,1){\circle*{.075}}
\put(1.9,1){\line(1,0){.5}}
\put(2.4,1){\circle*{.075}}
\put(2.4,1){\line(1,0){.5}}
\put(2.9,1){\circle*{.075}}
\put(2.9,1){\line(1,0){.5}}
\put(3.4,1){\circle*{.075}}
\put(3.4,1){\line(1,0){.5}}
\put(3.9,1){\circle*{.075}}
\put(3.9,1){\line(0,1){.4625}}
\put(3.9,1.5){\circle{.075}}
\put(3.9,1){\line(1,0){.5}}
\put(4.4,1){\circle*{.075}}
\put(1.775,1.05){\makebox(.25,.25){\footnotesize 1}}
\put(2.275,1.05){\makebox(.25,.25){\footnotesize 2}}
\put(2.775,1.05){\makebox(.25,.25){\footnotesize 2}}
\put(3.275,1.05){\makebox(.25,.25){\footnotesize 2}}
\put(3.875,1.05){\makebox(.25,.25){\footnotesize 2}}
\put(4.275,1.05){\makebox(.25,.25){\footnotesize 1}}
\put(3.775,1.55){\makebox(.25,.25){\footnotesize 1}}
\end{picture}

\hspace*{\fill}

\begin{picture}(3,1)(1.65,.5)
\thicklines
\put(1.9,1){\circle*{.075}}
\put(1.9,1){\line(1,0){.5}}
\put(2.4,1){\circle*{.075}}
\put(2.4,1){\line(1,0){.5}}
\put(2.9,1){\circle*{.075}}
\put(2.9,1){\line(1,0){.5}}
\put(3.4,1){\circle*{.075}}
\put(3.4,1){\line(1,0){.5}}
\put(3.9,1){\circle*{.075}}
\put(3.9,1){\line(1,0){.5}}
\put(4.4,1){\circle*{.075}}
\put(1.775,1.05){\makebox(.25,.25){\footnotesize 1}}
\put(2.275,1.05){\makebox(.25,.25){\footnotesize 1}}
\put(2.775,1.05){\makebox(.25,.25){\footnotesize 1}}
\put(3.275,1.05){\makebox(.25,.25){\footnotesize 1}}
\put(3.775,1.05){\makebox(.25,.25){\footnotesize 1}}
\put(4.275,1.05){\makebox(.25,.25){\footnotesize 1}}
\end{picture}

\caption[The Tyurina blowup]{The Tyurina blowup.  The special factorization of the (deformed) equation for a $D_{n+2}$ singularity allows us to blow up the open length 1 node in the top Dynkin diagram.  The resulting
space has a residual $A_{n+1}$ surface singularity.  Unfortunately, this is not the right node to blow up when resolving
a length 2 threefold singularity.  The small resolution for this singularity type can only come from blowing up
a length 2 node in the diagram for a corresponding hyperplane section.}  
\end{figure}
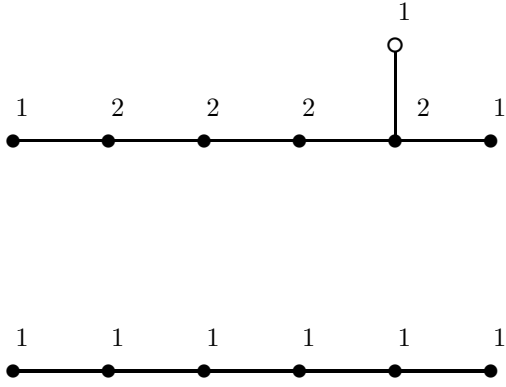

\setlength{\unitlength}{1 pt}

\subsection{Comparison with non-deformed version}

\noindent To ensure that our deformed equations match the non-deformed ones in the 
appropriate ($t=0$) limit, we explicitly compute the newly defined
polynomials in the non-deformed case:
\begin{eqnarray*}
f(U,t) &=& UP(-U^2,t)+Q(-U^2,t),\\
g(-U^2,t) &=& f(U,t)f(-U,t) = (Q + UP)(Q - UP) = Q^2 - U^2P^2,\\
g(Z,t) &=& ZP(Z,t)^2 + Q(Z,t)^2.
\end{eqnarray*}
Setting the deformation parameter $t=0$,
$$f(U,t) = U^{n+2}.$$
The computation of $P(Z),Q(Z)$ breaks up into four cases:
$$\begin{array}{cccccccc}
n = 0 &\mod 4:& P(Z)&=&0, & Q(Z)&=&-Z^{(n+2)/2},\\
n = 1 &\mod 4:& P(Z)&=&-Z^{(n+1)/2},& Q(Z) &=&0,\\
n = 2 &\mod 4:& P(Z)&=&0, & Q(Z)&=&Z^{(n+2)/2},\\
n = 3 &\mod 4:& P(Z)&=&Z^{(n+1)/2},& Q(Z) &=&0.\\
\end{array}$$
Moreover, in each of these cases we find
$$g(Z) = Z^{n+2}.$$
In particular, after doing the $D_{n+2} \flecha D_{n-m+1} \oplus
A_m$ blowup, the invariant functions $g_{n-m+1}'(Z)$ and
$f_{m+1}''(U)$ in the non-deformed case are given by:
$$\begin{array}{cccccccc}
m = 1 &\mod 4:& P''(Z)&=&0, & Q''(Z)&=&-Z^{(m+1)/2},\\
m = 2 &\mod 4:& P''(Z)&=&-Z^{m/2},& Q''(Z) &=&0,\\
m = 3 &\mod 4:& P''(Z)&=&0, & Q''(Z)&=&Z^{(m+1)/2},\\
m = 0 &\mod 4:& P''(Z)&=&Z^{m/2},& Q''(Z) &=&0,\\
\end{array}$$
and $$f''(U)=U^{m+1},\skop g'(Z)=Z^{n-m+1}.$$

\subsection{Matrix of relations for the full deformation}

\noindent {\bf Invariant theory}\\
\noindent If we are to break
$$D_{n+2} \flecha D_{n-m+1} \oplus A_m,$$
we expect the original equation with deformations to factor in a certain way.
Recalling the definition 
$$h(Z,t) = \dfrac{1}{Z}(g(Z,t)-\delta^2),$$
we see that we can write the deformed $D_{n+2}$ equation~\eqref{eq:D1} as
\begin{equation}\label{eq:D2}
\Phi_{D_{n+2}}: \;\; X^2 + Y^2Z - h(Z,t) + 2\delta Y,
\end{equation}
where
\begin{eqnarray*}
f_{n+2}(U,t) &=& \prod_{i=1}^{n+2} (U+t_i) = f'_{n-m+1}(U) f''_{m+1}(U),\\
g_{n+2}(Z,t) &=& \prod_{i=1}^{n+2} (Z+t_i^2) = g'_{n-m+1}(Z) g''_{m+1}(Z),\\
h_{n+1}(Z,t) &=& \frac{1}{Z}(g - \delta^2).
\end{eqnarray*}
(Note that primes do not indicate derivatives!)  Recall that
$$\delta = f(0,t) = \prod_{i=1}^{n+2} t_i, \skop \delta' = \sqrt{g'(0)} = f'(0).$$
The allowed invariant functions after breaking to $D_{n-m+1} \oplus A_m$ will be:
$$g'_{n-m+1}(Z), \;\; f''_{m+1}(U), \;\; \and \;\; \delta'.$$
We can also replace $g'_{n-m+1}(Z)$ with
$$h'_{n-m}(Z) =\frac{1}{Z}(g'_{n-m+1}(Z) - \delta'^2).$$
Recalling the $P, Q$ decomposition for $f(U)$ into odd and even parts, we have
$f_{m+1}''(U) = U P''(-U^2) + Q''(-U^2).$  This yields
$$g''_{m+1} = Z P''(Z)^2 + Q''(Z)^2.$$
The importance of this is that we cannot directly use $g''_{m+1}$, but since
$f''_{m+1}$ is legal after the symmetry breaking, $P''$ and $Q''$ are.
The legal decomposition of $g_{n+2}$ after the blowup is then,
$$ g_{n+2}(Z) = g'_{n-m+1}(Z)(ZP''(Z)^2 + Q''(Z)^2).$$
Plugging this into the deformed $D_{n+2}$ equation yields:
$$\Phi_{D_{n+2}}: \;\; X^2 + Y^2Z - \frac{1}{Z}(g'(ZP''^2+Q''^2) - \delta'^2\delta''^2) + 2\delta'\delta''Y.$$
Note that without deformations ($t=0$),
\begin{eqnarray*}
m \;\;\even &\implies & g'(Z) = Z^{n-m+1},\;\;\; P''(Z) = (-1)^{m/2}Z^{m/2},\;\;\; Q''(Z)=0.\\
m \;\;\odd &\implies& g'(Z) = Z^{n-m+1}, \;\;\; P''(Z)=0,\;\;\;
Q''(Z)=(-1)^{(m+1)/2}Z^{(m+1)/2}.
\end{eqnarray*}

\noindent {\bf Quasi-homogeneity of the matrices}\\
\noindent Recalling the matrices of relations for these two cases:
in the case of $m$ odd, the relations matrix is
$$
\left( \begin{array}{cccc}
-iX & Z^{(m+1)/2} & -Y & 0\\
Z^{n+1-(m+1)/2} & iX & 0 & -Y\\
-YZ & 0 & -iX & Z^{(m+1)/2}\\
0 & -YZ & Z^{n+1-(m+1)/2} & iX
\end{array} \right),$$
and in the case of $m$ even, the relations matrix is
$$
\left(\begin{array}{cccc}
-iX & 0 & -Y & Z^{m/2}\\
0& iX &Z^{n-m/2}&-Y\\
-YZ&Z^{m/2+1}&-iX&0\\
Z^{n+1-m/2}&-YZ&0&iX\\
\end{array}\right),$$
we are tempted to make some natural guesses for the deformed
analogues. We need some clues from the quasi-homogeneity of the
matrix.  We assign weights in the equation
$$X^2 + Y^2Z - Z^{n+1}$$
as follows:
$$\deg Z = 2,\;\;\; \deg X = n+1, \;\;\; \deg Y = n, \;\;\
\deg(\mathrm{eqn}) = 2n+2.$$

\noindent A combined matrix which reduces to each of these in the $m$ odd and $m$ even cases
will have in each entry the degrees
$$
\left( \begin{array}{cccc}
n+1 & m+1 & n & m\\
2n-m+1 & n+1 & 2n-m & n\\
n+2 & m+2 & n+1 & m+1\\
2n-m+2 & n+2 & 2n-m+1 & n+1
\end{array} \right)$$

\noindent We know that $g'(Z) = Z^{n-m+1}+...$, and so $\deg g'(Z) =
2n-2m+2.$ Hence, from the deformed equation:
$$X^2 + Y^2Z - g'(Z)P''(Z)^2 - \frac{g'(Z)Q''(Z)^2}{Z},$$
we see that
$$\deg P''(Z) = m, \;\;\; \deg Q''(Z) = m+1.$$

\noindent We record for reference the degree of each element:
\begin{scriptsize}
$$\begin{array}{c|c|c|c|c|c|c|c|c|c|c|c}
\eqn & X   & Y & Z & g' & h'    & h'' & \delta' & \delta'' & P'' & S'' & Q''\\
2n+2 & n+1 & n & 2 & 2n-2m+2 & 2n-2m & 2m  & n-m+1   & m+1      & m   & m-1 & m+1\\
\end{array}$$
\end{scriptsize}
The fully deformed equation for the $D_{n+2}$ singularity is given
by
$$\eqn:\;\; X^2 + Y^2Z - h + 2\delta Y = 0,$$
where
\begin{eqnarray*}
Q &=& SZ+\delta,\\
h &=& (g-\delta^2)/Z = P^2 + (Q^2-\delta^2)/Z = P^2 + ZS^2 + 2\delta S\\
  &=& Zh'h''+h'\delta''^2+h''\delta'^2,\\
h'' &=& P''^2+ZS''^2+2\delta''S'',\\
Q'' &=& S''Z+\delta'',\\
\delta &=& \delta'\delta''.
\end{eqnarray*}
The degrees of the various elements that can go in the matrix are
$$
\begin{array}{c|c|c|c|c|c|c|c|c}
\mathrm{eqn} & X & Y & Z & h' & \delta' & P'' & S'' & Q''\\
2n+2 & n+1 & n & 2 & 2n-2m & n-m+1 & m & m-1 & m+1
\end{array}$$
Note that when $t=0$ we have $\delta'=0$ and
$$\begin{array}{c|ccc}
m \mathrm{\;even} & h'(Z)=Z^{n-m}, & P''(Z) = i^m Z^{m/2}, & Q''(Z)=0.\\
m \mathrm{\;odd} & h'(Z)=Z^{n-m}, & P''(Z) = 0, & Q''(Z)= i^{m+1}Z^{(m+1)/2}.
\end{array}$$

\noindent{\bf Finding fully deformed matrix of relations}\\
\noindent The fully deformed matrix of relations is found by requiring that
(1) it reduce properly to the non-deformed case, (2) the entries
respect quasi-homogeneity, (3) the relations in the ideal of $2
\times 2$ minors match those found in the non-deformed case. 

The most general matrix with these degrees and which matches at $t=0$ is thus
\begin{scriptsize}
$$
\left( \begin{array}{cccc}
-iX + a_1 \delta' P'' & i^{\pm(m+1)} Q'' & -Y + a_2\delta' S'' & i^{\pm m}P''\\
i^{\pm(m+1)} h' Q''+b_1 \delta' Y & iX + b_2 \delta' P'' & i^{\pm m} h'P'' & -Y + b_3\delta' S''\\
-YZ + c_1\delta' Q'' + c_2 \delta' S'' Z & i^{\pm m}ZP'' & -iX + c_3\delta' P'' & i^{\pm(m+1)}Q''\\
i^{\pm m}h'P''Z + d_1\delta' X + d_2 \delta'^2 P'' & -YZ + d_3\delta' Q'' + d_4 \delta' S'' Z & 
i^{\pm(m+1)} h' Q''+d_5 \delta' Y & iX + d_6 \delta' P''
\end{array} \right).$$
\end{scriptsize}

\noindent To find the right coefficients, we will want the minors of this matrix to satisfy the same relations that
the minors in the $t=0$ case satisfy.  For the top two rows and the bottom two rows, these are
\begin{enumerate}
\item $M_{12}^{14} = M_{12}^{23}$, or $[3]=[4]$.  This requires
      \begin{itemize}
      \item $b_3 = a_2$
      \item $b_2 = -a_1$
      \item $b_1 = 0$
      \end{itemize}
\item $M_{12}^{13} = (-1)^m h' M_{12}^{24}$, or $[2] = (-1)^m h'[5]$.
\item $M_{12}^{12} + Z M_{12}^{34} = \eqn$, or $[1] + Z [6] = \eqn$.
\item $M_{34}^{14} = M_{34}^{23}$, or $[3b]=[4b]$. 
\item $M_{34}^{13} = (-1)^m h' M_{34}^{24}$, or $[2b] = (-1)^m h'[5b]$.
\item $M_{34}^{12} + Z M_{34}^{34} = Z\eqn$, or $[1b] + Z [6b] = Z\eqn$.
\end{enumerate}
 
\noindent Imposing the first condition $[3]=[4]$, the matrix becomes
\begin{scriptsize}
$$
\left( \begin{array}{cccc}
-iX + a_1 \delta' P'' & i^{\pm(m+1)} Q'' & -Y + a_2\delta' S'' & i^{\pm m}P''\\
i^{\pm(m+1)} h' Q'' & iX - a_1 \delta' P'' & i^{\pm m} h'P'' & -Y + a_2 \delta' S''\\
-YZ + c_1\delta' Q'' + c_2 \delta' S'' Z & i^{\pm m}ZP'' & -iX + c_3\delta' P'' & i^{\pm(m+1)}Q''\\
i^{\pm m}h'P''Z + d_1\delta' X + d_2 \delta'^2 P'' & -YZ + d_3\delta' Q'' + d_4 \delta' S'' Z & 
i^{\pm(m+1)} h' Q''+d_5 \delta' Y & iX + d_6 \delta' P''
\end{array} \right).$$ 
\end{scriptsize}
The top two rows suggest that we can make the change of coordinates 
$$X \mapsto \til{X} = X+ia_1\delta'P'', \skop Y \mapsto \til{Y} = Y-a_2\delta'S'',$$
which would additionally require $$ c_3 = a_1, \skop d_6 = -a_1, \skop c_2 = a_2, \skop \and \skop d_4 = a_2.$$
Computing
\begin{eqnarray*}
\til{X}^2 &=& X^2 + 2ia_1\delta'XP'' - a_1^2\delta'^2P''^2 = X^2 + 2ia_1\delta'(\til{X} -ia_1\delta'P'')P'' - a_1^2\delta'^2P''^2\\
&=& X^2 + 2ia_1\delta'\til{X}P''+a_1^2\delta'^2P''^2,\\
\til{Y}^2Z &=& Y^2Z - 2a_2\delta'YS''Z + a_2^2\delta'^2S''^2Z,
\end{eqnarray*}
and expanding the original deformed equation
\begin{eqnarray*}
0 &=& X^2 + Y^2Z - h + 2\delta Y\\
&=& X^2 + Y^2Z - h'(ZP''^2+Q''^2)-\delta'(2YZS''-2YQ'')-\delta'^2(P''^2-ZS''^2+2S''Q'')
\end{eqnarray*}
we see that the coordinate change only makes sense if
$$a_1 = \pm i, \skop \and \skop a_2=1.$$
In this case the equation becomes
$$0 = \til{X}^2 + \til{Y}^2Z - h'(ZP''^2+Q''^2)+2\delta'(\til{Y}Q''\pm \til{X} P'')$$
and the matrix is
$$
\left( \begin{array}{cccc}
-i\til{X} & i^{\pm(m+1)} Q'' & -\til{Y} & i^{\pm m}P''\\
i^{\pm(m+1)} h' Q'' & i\til{X} & i^{\pm m} h'P'' & -\til{Y}\\
-\til{Y}Z + c_1\delta' Q'' & i^{\pm m}ZP'' & -i\til{X} & i^{\pm(m+1)}Q''\\
i^{\pm m}h'P''Z + d_1\delta' \til{X} & -\til{Y}Z + d_3\delta' Q'' & 
i^{\pm(m+1)} h' Q''+d_5 \delta' \til{Y} & i\til{X}
\end{array} \right).$$ 
Requiring conditions 2 and 3 forces us to choose $a_1 = (-1)^{m+1}i$, giving change of coords
$$X \mapsto \til{X} = X+(-1)^{m}\delta'P'', \skop Y \mapsto \til{Y} = Y-\delta'S'',$$
and deformed equation
\begin{equation}\label{eq:deformedD}
0 = \til{X}^2 + \til{Y}^2Z - h'(ZP''^2+Q''^2)+2\delta'(\til{Y}Q'' +(-1)^{m+1} \til{X} P'').
\end{equation}
Fixing the signs (and the remaining coefficients to zero) our matrix is now
\begin{equation}\label{eq:deformed-relations}
\mathcal{R} = \left( \begin{array}{cccc}
-i\til{X} & i^{-(m+1)} Q'' & -\til{Y} & i^{m}P''\\
i^{(m+1)} h' Q'' & i\til{X} & i^{-m} h'P'' & -\til{Y}\\
-\til{Y}Z & i^{m}P''Z & -i\til{X} & i^{-(m+1)}Q''\\
i^{-m}h'P''Z  & -\til{Y}Z  & i^{(m+1)} h' Q''& i\til{X}
\end{array} \right).
\end{equation}
Note that the odd and even cases are now treated simultaneously.  Together with the syzygies matrix
$\mathcal{S}$, this yields a matrix factorization for the deformed equation~\eqref{eq:D2}.

\subsection{The blowup}

We can find additional relations among the Pl\"ucker coordinates by restricting attention to 
the 2 $\times$ 2 minors in the top two rows of the relations matrix $\mathcal{R}$~\eqref{eq:deformed-relations}.
Together with the relation~\eqref{grass}
$$\gamma \varphi + \delta \alpha - \beta \varepsilon = 0,$$
this gives us the blowup.

The 2 $\times$ 2 minors for the top two rows of this matrix are
$$\begin{array}{cccccc}
\;[1] &=& \til{X}^2 - h'Q''^2, & [6] &=& \til{Y}^2 - h'P''^2 \\
 \;[2] &=& i^{-m+3}h'(\til{X}P''+i^{2m-2}\til{Y}Q''), & [5] &=& i^{m+3}(\til{X}P'' + i^{2m-2}\til{Y}Q'') \\
 \;[3] &=& i\til{X}\til{Y} - i^{2m+1}h'P''Q'', & [4] &=& i\til{X}\til{Y} - i^{2m+1}h'P''Q'' \\
\end{array}$$
and they satisfy the relations
\begin{itemize}
\item $[3]=[4]$
\item $[2]=(-1)^mh'[5]$, and
\item $[1]+Z[6]+2i^{(m-1)}\delta'[5] = [\eqn]$.
\end{itemize}
The 2 by 2 minors for the bottom two rows are
$$\begin{array}{cccccc}
\;[6b] &=& \til{X}^2 - h'Q''^2, & [1b] &=& Z^2(\til{Y}^2 - h'P''^2) \\
 \;[2b] &=& i^{m+3}h'Z(\til{Y}Q''+i^{2m-2}\til{X}P''), & [5b] &=& i^{-m+3}Z(\til{Y}Q''+i^{2m-2}\til{X}P'') \\
 \;[3b] &=& -i\til{X}\til{Y} + i^{2m+1}h'P''Q'', & [4b] &=& -i\til{X}\til{Y} + i^{2m+1}h'P''Q'' \\
\end{array}$$
and they satisfy almost the same relations
\begin{itemize}
\item $[3b]=[4b]$
\item $[2b]=(-1)^mh'[5b]$, and
\item $[1b]+Z[6b]+2i^{(m+1)}\delta'[5b] = Z[\eqn]$.
\end{itemize}

\noindent This motivates us to restrict attention to the Pl\"ucker coordinates
\begin{eqnarray}
\alpha &=& [4] = i\til{X}\til{Y} - i^{2m+1}h'P''Q''\\
\varepsilon &=& [5] = i^{m+3}\til{X}P''+i^{-m+1}\til{Y}Q''\\
\varphi &=& [6] = \til{Y}^2 - h'P''^2 \label{eq:plucker}
\end{eqnarray}
and the equation for the Grassmanian $G(2,4)\subseteq \PP^5$ ($[1]\varphi - [2]\varepsilon + [3]\alpha = 0$)
becomes 
$$(-Z\varphi-2i^{m-1}\delta'\varepsilon)\varphi + (-1)^{m+1}h'\varepsilon^2 + \alpha^2 = 0.$$
It will be useful to rewrite this as
\begin{equation}\label{eq:def-grass}
\alpha^2 - \varphi^2Z + (-1)^{m+1}h'\varepsilon^2 + 2i^{m+1}\delta'\varepsilon\varphi = 0.
\end{equation}

\noindent \textbf{Identification of residual singularities in the hyperplane}\\
\noindent The interesting charts for the blowup are $\varphi = 1$ and $\varepsilon = 1$.\skp

\noindent \framebox{$\varphi = 1$} 
In this case we have $\alpha = \varphi \alpha_1$, $\varepsilon = \varphi \varepsilon_1$, and
the equation~\eqref{eq:def-grass} becomes
\begin{equation}\label{eq:phi-grass}
\alpha_1^2 - Z + (-1)^{m+1}h'\varepsilon_1^2 + 2i^{m+1}\delta'\varepsilon_1 = 0.
\end{equation}
To see what kind of singularity this gives, define $G''(Z,i\alpha_1,t)$ as in~\eqref{eq:G-defn} by
$$i\alpha_1 P''(Z,t)+ Q''(Z,t) = (Z-\alpha_1^2)G''(Z,i\alpha_1,t)+ f''_{m+1}(i\alpha_1,t),$$
and notice that from~\eqref{eq:phi-grass} we have
$$Z-\alpha_1^2 = (-1)^{m+1}h'\varepsilon_1^2 + 2i^{m+1}\delta'\varepsilon_1.$$
Moreover, notice from~\eqref{eq:plucker} (dropping tildes) that
$$i\alpha_1 P'' + Q'' = \frac{i\alpha P'' + \varphi Q''}{\varphi} =
\frac{-XYP'' + i^{2m}h'P''^2Q'' + Y^2 Q'' - h'P''^2Q''}{Y^2-h'P''^2}.$$
When $m$ is even, this simplifies to
$$i\alpha_1 P'' + Q'' = i^{m+3}Y\varepsilon_1,$$
and we have
$$i^{m+3}Y\varepsilon_1 = \varepsilon_1((-1)^{m+1}h'\varepsilon_1 + 2i^{m+1}\delta')G'' + f''_{m+1}.$$
Defining
$$\Hat{Y} = i^{m+3}Y - ((-1)^{m+1}h'\varepsilon_1 + 2i^{m+1}\delta')G'',$$
we find that our equation in this chart is
$$\Hat{Y}\varepsilon_1 = f''_{m+1}(i\alpha_1,t),$$
which is precisely the equation for a deformed $A_m$ singularity!\skp

\noindent \framebox{$\varepsilon = 1$}
In this case we have $\alpha = \varepsilon \alpha_2$, $\varphi = \varepsilon \varphi_2$, 
and the equation~\eqref{eq:def-grass} becomes
$$\alpha_2^2 - \varphi_2^2Z + (-1)^{m+1}h' + 2i^{m+1}\delta'\varphi_2 = 0,$$
which we immediately recognize as the preferred versal form for a deformed $D_{n-m+1}$ singularity!

We thus see the expected residual singularities, since in the appropriate hyperplane
the blowup breaks $D_{n+2} \flecha D_{n-m+1} \oplus A_m$.  Notice, however, that the total space for
the resolved threefold
$$\alpha^2 - \varphi^2Z + (-1)^{m+1}h'\varepsilon^2 + 2i^{m+1}\delta'\varepsilon\varphi = 0,$$
is smooth.\skp

\setlength{\unitlength}{1 true in}

\begin{figure}[ht]

\begin{picture}(3,1)(1.65,.5)
\thicklines
\put(1.9,1){\circle*{.075}}
\put(1.9,1){\line(1,0){.5}}
\put(2.4,1){\circle*{.075}}
\put(2.4,1){\line(1,0){.4625}}
\put(2.9,1){\circle{.075}}
\put(2.9375,1){\line(1,0){.4625}}
\put(3.4,1){\circle*{.075}}
\put(3.4,1){\line(1,0){.5}}
\put(3.9,1){\circle*{.075}}
\put(3.9,1){\line(1,0){.5}}
\put(4.4,1){\circle*{.075}}
\put(4.4,1){\line(3,4){.3}}
\put(4.7,1.4){\circle*{.075}}
\put(4.4,1){\line(3,-4){.3}}
\put(4.7,.6){\circle*{.075}}
\put(1.775,1.05){\makebox(.25,.25){\footnotesize 1}}
\put(2.275,1.05){\makebox(.25,.25){\footnotesize 2}}
\put(2.775,1.05){\makebox(.25,.25){\footnotesize 2}}
\put(3.275,1.05){\makebox(.25,.25){\footnotesize 2}}
\put(3.775,1.05){\makebox(.25,.25){\footnotesize 2}}
\put(4.2,1.05){\makebox(.25,.25){\footnotesize 2}}
\put(4.575,1.45){\makebox(.25,.25){\footnotesize 1}}
\put(4.575,.3){\makebox(.25,.25){\footnotesize 1}}
\end{picture}

\bigskip
\bigskip

\begin{picture}(3,1)(1.65,.5)
\thicklines
\put(1.9,1){\circle*{.075}}
\put(1.9,1){\line(1,0){.5}}
\put(2.4,1){\circle*{.075}}
\put(3.4,1){\circle*{.075}}
\put(3.4,1){\line(1,0){.5}}
\put(3.9,1){\circle*{.075}}
\put(3.9,1){\line(1,0){.5}}
\put(4.4,1){\circle*{.075}}
\put(4.4,1){\line(3,4){.3}}
\put(4.7,1.4){\circle*{.075}}
\put(4.4,1){\line(3,-4){.3}}
\put(4.7,.6){\circle*{.075}}
\put(1.775,1.05){\makebox(.25,.25){\footnotesize 1}}
\put(2.275,1.05){\makebox(.25,.25){\footnotesize 2}}
\put(3.275,1.05){\makebox(.25,.25){\footnotesize 2}}
\put(3.775,1.05){\makebox(.25,.25){\footnotesize 2}}
\put(4.2,1.05){\makebox(.25,.25){\footnotesize 2}}
\put(4.575,1.45){\makebox(.25,.25){\footnotesize 1}}
\put(4.575,.3){\makebox(.25,.25){\footnotesize 1}}
\end{picture}

\caption[Small resolution via blowup of a length 2 node]{The blow up of the open length 2 node in the top Dynkin diagram provides a small
resolution for the length 2 threefold singularity. When restricted to the hyperplane section, the $D_8$ surface singularity splits into two lower-order surface singularities: $A_2$ and $D_5$.  Nevertheless, the total space
for the blow up (the threefold) is smooth.}
\end{figure}
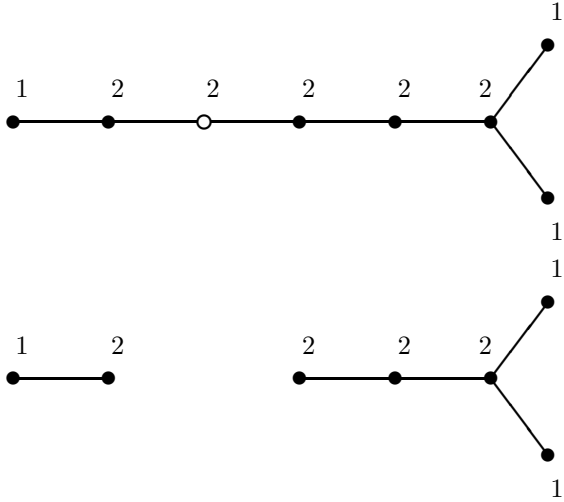

\setlength{\unitlength}{1 pt}

\noindent {\bf Transition functions}\\
\noindent The coordinates in our two charts are $(\alpha_1 = \alpha/\varphi, \varepsilon_1 = \varepsilon/\varphi,
Z, t)$ for the $\varphi=1$ chart, and $(\alpha_2 = \alpha/\varepsilon, \varphi_2 = \varphi/\varepsilon, Z, t)$ for 
the $\varepsilon = 1$ chart.  The transition functions between them can be written
\begin{eqnarray}\label{eq:length2transfns}
\varphi_2 &=& \varepsilon_1^{-1},\\
\alpha_2 &=& \varepsilon_1^{-1}\alpha_1,\\
Z &=& (-1)^{m+1} \varepsilon_1^2 h'(Z,t) + \alpha_1^2 + 2i^{m+1}\delta'\varepsilon_1,\\
t &=& t.
\end{eqnarray}

\section{Difficulties with the approach}
Although we were able to successfully perform explicit small resolutions for deformed $A_k$ and $D_{k+2}$ equations,
it is not at all clear how to use these blowups for finding the associated matrix model superpotentials.  The results
do not come to us packaged in such a way that comparison with Ferrari's framework can be made.  Only in the $A_1$ case were we able to make the connection, and here the answer was already known from Ferrari's example 1. 

The next section will develop an alternative ``top-down'' approach, beginning with the Intriligator--Wecht superpotentials, which we have conjectured provide the answer.  This will prove to be a more fruitful strategy, although we will still need to use our small resolution techniques in order to show that we get the right
singular space in the length 2 case.  The algorithm we devise for performing the blow-downs cannot guarantee that we have all global holomorphic functions, and hence each case must be checked by blowing back up and recovering the original space.
\end{onecolumn}

\begin{onecolumn}
\chapter{Small blow-downs from global holomorphic functions}\label{ch:algorithm}
\section{``Top-down'' approach}

We have seen in Chapter 5 the difficulties that arise in the ``bottom-up" approach.  It is hard to perform small resolutions explicitly (thus far we have only done so in the $A_{n-1}$ and $D_{n+2}$ cases), but even this is not enough.  Only in the $A_1$ case were we able to identify the superpotential from our small resolution, and in this case we already knew the answer!

While we know from the Katz--Morrison result \cite{morrison} what the singular geometries should be, it is hard to fit their small resolutions into Ferrari's framework, and hence we are unable to find a corresponding superpotential.
On the other hand, if we could guess the form of the superpotentials, then showing which singular geometries they
correspond to might not be so difficult.  As in Ferrari's original examples, we could start with the resolved geometry, find global holomorphic functions, and then use these to blow down the exceptional $\PP^1$'s.  We will 
call this the ``top-down'' approach:

\begin{itemize}
\item Guess $W$, then compute transition functions for resolved geometry $\Hat{\M}$.
\item Find global holomorphic functions (hard, need algorithm!).
\item Use these to construct the blow-down; relations among the global holomorphic functions describe
the singular space $\M_0$. 
\item Blow back up to check that we recover the original resolved geometry $\Hat{\M}$.
\end{itemize}

The last step is necessary because we do not yet know a way of determining whether or not we have found all
(independent) global holomorphic functions.  This is because the algorithm searches for global holomorphic functions which can be constructed from a given list of monomials.  While we can be certain we haven't ``missed'' any functions which can be made out of monomials in the list, it is, of course, impossible to check all possible monomials.  It is thus only in blowing back up that can we confirm that the singular space is correct.  This will require using the small resolution methods developed in Chapter~\ref{ch:smallres}.

Our guesses for $W(x,y)$ come from the Intriligator--Wecht classification of RG fixed points of $\mathcal{N}=1$ SQCD with adjoints.  We make the most naive conjecture possible:  that the ADE superpotentials in the classification correspond (under Ferrari's construction) to geometries with ADE singularities.  Moreover, in order to perform the blow-downs we devise an algorithm which searches for global holomorphic functions, thereby proving Proposition 2:

\begin{myproposition2}  There is an algorithm for blowing down the exceptional $\PP^1$'s when the resolved geometry is given by simple transition functions~\eqref{eq:transfns} as in Ferrari's construction.
\end{myproposition2}

\section{The Intriligator--Wecht classification}

Here we recall (see Introduction) the Intriligator--Wecht classification of RG fixed points of $\mathcal{N}=1$ SQCD with adjoints \cite{wecht}.  Using ``a-maximization'' and doing a purely field theoretic analysis, they classify all
relevant adjoint superpotential deformations for 4d $\N=1$ SQCD with
$\N_f$ fundamentals and $\N_a=2$ adjoint matter chiral superfields,
$X$ and $Y$.   The possible RG fixed points, together with the map of possible flows between fixed points, are summarized below:

\setlength{\unitlength}{1 true in}

\begin{figure}[h]
\begin{picture}(6,3)(0,0)
\put(.8,1.2){\makebox{
$\begin{array}{c|c}
\mathrm{type} & W(X,Y)\\
\hline
&\\
\Hat{O} & 0\\
\Hat{A} & \Tr Y^2\\
\Hat{D} & \Tr XY^2\\
\Hat{E} & \Tr Y^3\\
A_k & \Tr(X^{k+1}+Y^2)\\
D_{k+2} & \Tr(X^{k+1}+XY^2)\\
E_6 & \Tr(Y^3+X^4)\\
E_7 & \Tr(Y^3+YX^3)\\
E_8 & \Tr(Y^3+X^5)
\end{array}$}}
\put(3.7,0){\makebox{\includegraphics[scale=.7]{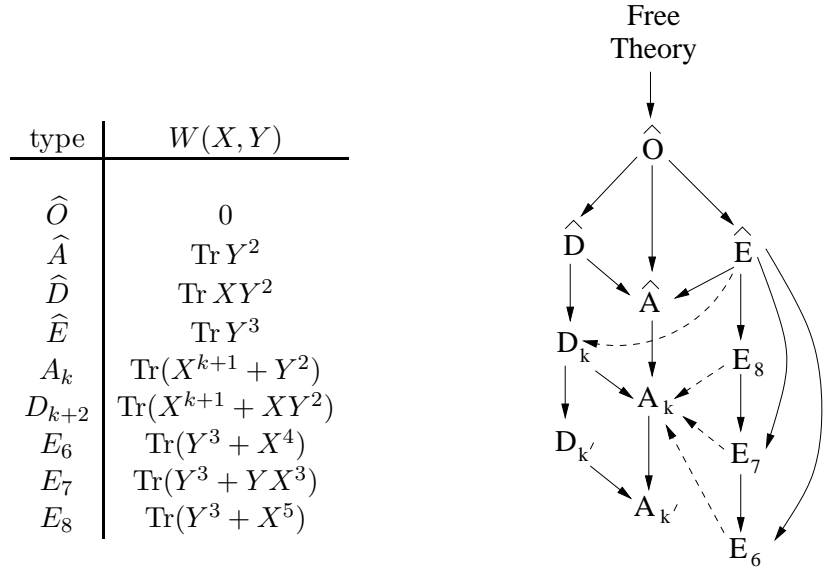}}}
\end{picture}
\caption[Intriligator--Wecht Classification of RG Fixed Points]{Intriligator--Wecht Classification 
of RG Fixed Points.  The diagram on the right shows the map of possible
flows between fixed points.  Dotted lines indicate flow to a particular value of $k$.  Note
that $k' < k$. \cite[pages 3-4]{wecht}}
\end{figure}

\setlength{\unitlength}{1 pt}

\begin{table}[h]
$$\begin{array}{c|c|c|c}
\mathrm{type} & W(x,y) & \partial_{w_1}E(\gamma,w_1) & \PP^1: (w_1,w_2), (v_1,v_2)\\
\hline
&&&\\
\Hat{O} & 0 & 0 & (x+\gamma y, 0), (\beta x + y, 0)\\
& & &\\
\Hat{A} & \12 y^2 & w_1 \leftrightarrow \gamma^2 w_1 & (x,0), (\beta x, x)\\
& & &\\
\Hat{D} & xy^2 & w_1^2 \leftrightarrow \gamma w_1^2 & (x,0), (\beta x, x^2)\\
& & &\\
\Hat{E} & \dfrac{1}{3} y^3 & \gamma^{-1}w_1^2 \leftrightarrow \gamma^2 w_1^2& (x,0), (\beta x, \beta x^2)\\
& & &\\
A_k & \dfrac{1}{k+1}x^{k+1}+\12 y^2 & \gamma^2 w_1^k + w_1 \leftrightarrow
\gamma^{1-k}w_1^k+\gamma^2w_1 & (0,0), (0,0)\\
& & &\\
D_{k+2} & \dfrac{1}{k+1}x^{k+1}+xy^2 & \gamma^2 w_1^k+w_1^2 \leftrightarrow
\gamma^{1-k} w_1^k + \gamma w_1^2 & (0,0), (0,0)\\
& & &\\
E_6 & \dfrac{1}{3}y^3+\dfrac{1}{4}x^4 & \gamma^{-1}w_1^2+\gamma^2 w_1^3 \leftrightarrow
\gamma^2 w_1^2 + \gamma^{-2} w_1^3 & (0,0), (0,0)\\
& & &\\
E_7 & \dfrac{1}{3}y^3+yx^3 & \gamma^{-1}w_1^2+\gamma w_1^3 \leftrightarrow
\gamma^2 w_1^2 + \gamma^{-1} w_1^3 & (0,0), (0,0)\\
& & &\\
E_8 & \dfrac{1}{3}y^3 + \dfrac{1}{5}x^5 & \gamma^{-1}w_1^2+\gamma^2 w_1^4
\leftrightarrow \gamma^2 w_1^2 + \gamma^{-3} w_1^4 & (0,0), (0,0)
\end{array}$$
\caption{Identification of corresponding resolved geometries}
\label{identify}
\end{table}

Using Table 3.2 from Chapter 3 we identify the corresponding resolved geometries, and these are shown in Table~\ref{identify}.  There are two possible perturbation terms $\partial_{w_1}E(\gamma,w_1)$ for each case due to the $x \leftrightarrow y$ symmetry, which exchanges
$$\gamma^n w_1^m \leftrightarrow \gamma^{3-n-m} w_1^m.$$
All of these examples are in the case $M=2$, with transition functions
$$\beta=\gamma^{-1}, \skop v_1=\gamma^{-1}w_1, \skop v_2=\gamma^3 w_2+\partial_{w_1}E(\gamma,w_1),$$
which gives a deformed $\OO(1)\oplus\OO(-3)$ bundle.
In the $\Hat{A}$ and $A_k$ cases, the linear perturbation term $w_1$ indicates that the normal bundle
to the $\PP^1$ is not really $\OO(1)\oplus\OO(-3)$, but rather a funny description for $\OO\oplus\OO(-2)$, or
$\OO(-1)\oplus\OO(-1)$ in the $A_1$ case.

\section{The problem}

We have seen from Ferrari's examples in Chapter~\ref{ch:ferrari} that we can blow down exceptional $\PP^1$'s
by finding global holomorphic functions (ghf's) on the resolved space $\Hat{\M}$.  Global holomorphic functions must
necessarily be constant on the $\PP^1$, so they provide natural coordinates for the blow-down.  Beginning
with transition functions
$$\beta = \gamma^{-1}, \skop v_1=\gamma^{-n}w_1, \skop
v_2=\gamma^{-m}w_2+\partial_{w_1}E(\gamma,w_1),$$
the idea is to find polynomials in $\CC[\beta,v_1,v_2]$ which are also polynomials in $\CC[\gamma,w_1,w_2]$. 

For example, in the case of $E_6$:
$$\beta = \gamma^{-1}, \skop v_1 = \gamma^{-1} w_1, \skop
v_2 = \gamma^3 w_2 + \gamma^2 w_1^3 + \gamma^{-1} w_1^2.$$
we immediately see from the transition functions that neither $\beta, v_1$ or $v_2$ are ghf's.  However, if 
we consider $\beta v_2$ we have
$$\beta v_2 = \gamma^2 w_2 + \gamma w_1^3 + \gamma^{-2} w_1^2.$$
The non-holomorphic term on the right can be moved to the left, using $v_1 = \gamma^{-1}w_1$, and so 
we construct the ghf
$$\beta v_2 - v_1^2 = \gamma^2 w_2 + \gamma w_1^3.$$
Multiplying by $\beta, v_1$ or $v_2$ yields new ghf's
\begin{eqnarray*}
\beta(\beta v_2 - v_1^2) &=& \gamma w_2 + w_1^3,\\
v_1( \beta v_2 - v_1^2) &=& \gamma w_1 w_2 + w_1^4,\\
v_2(\beta v_2 - v_1^2) &=& (\gamma^4 w_2 + \gamma^3 w_1^3 + w_1^2)(\gamma w_2 + w_1^3).
\end{eqnarray*}
Unfortunately, these functions all vanish on $\beta v_2 = v_1^2$, which is more than the $\PP^1$ we are 
trying to blow down (the defining equations for the $\PP^1$ are $v_1 = v_2 = 0$).  We must find a collection of ghf's which collapse the exceptional $\PP^1$ to a point, but provide
a birational isomorphism everywhere else.  

This rather ad-hoc method for finding global holomorphic functions turns out to be quite limited, as cleverness quickly runs dry with polynomials of higher degree.
We would like to make the search for ghf's more systematic.

\section{Reformulation of the problem}

\noindent We reformulate our problem of finding global holomorphic functions as an
ideal membership problem.  We begin by illustrating the reformulation in the
case of $E_6$:
$$\beta = \gamma^{-1}, \skop v_1 = \gamma^{-1} w_1, \skop
v_2 = \gamma^3 w_2 + \gamma^2 w_1^3 + \gamma^{-1} w_1^2.$$
We can think of this as describing a variety in $\CC^6$, defined
by the following ideal $I \subset \CC[\gamma,w_1,w_2,\beta,v_1,v_2]$:
$$I = \langle \beta\gamma-1, v_1-\beta w_1, v_2-\gamma^3 w_2 - \gamma^2 w_1^3 - \beta w_1^2 \rangle.$$
In order to blow down the exceptional $\PP^1$, we must find functions which are holomorphic
in each coordinate chart, and will therefore be constant on the $\PP^1$.  Such global holomorphic
functions correspond to elements of the ideal $I$ that can be written in the form\footnote{Note that in our $E_6$
example none of the defining generators of $I$ are of this form!}
$$f - g \in I, \:\;\; \where \;\;\; f \in \CC[\beta, v_1, v_2], \;\; g \in \CC[\gamma, w_1, w_2].$$
For each such element, the global holomorphic function is $f = g$.

In general, we begin with transition functions
$$\beta = \gamma^{-1}, \skop v_1=\gamma^{-n}w_1, \skop
v_2=\gamma^{-m}w_2+\partial_{w_1}E(\gamma,w_1),$$
and form the ideal
$$I = \langle \beta\gamma-1, v_1-\beta^n w_1, v_2-\beta^m w_2 - \partial_{w_1}\til{E}(\gamma,w_1)\rangle,$$
where $\til{E}(\gamma,w_1)$ is obtained from $E(\gamma,w_1)$ by replacing all instances of $\gamma^{-1}$ 
with $\beta$.

\section{The algorithm}

\def\pure{\mathrm{pure}}
\def\mixed{\mathrm{mixed}}

Consider a monomial $\beta^i v_1^j v_2^k \in \CC[\beta, v_1, v_2]$.  Using Groebner basis techniques,
we can easily reduce this modulo the ideal
$$I = \langle \beta\gamma-1, v_1-\beta^n w_1, v_2-\beta^m w_2 - \partial_{w_1}\til{E}(\gamma,w_1)\rangle,$$
which is determined by our particular geometry.  In general, we will find
$$\beta^i v_1^j v_2^k \stackrel{\mathrm{mod\;}I}{\equiv} \pure(\beta,v_1,v_2) + 
\pure (\gamma,w_1,w_2) + \mixed, $$
where ``pure'' and ``mixed'' stand for pure and mixed terms\footnote{We will refer to any monomial in $\CC[\gamma,w_1,w_2,\beta,v_1,v_2]$ which does not belong to either $\CC[\beta,v_1,v_2]$ or $\CC[\gamma,w_1,w_2]$ as a mixed term.} in the appropriate variables.  We can then
bring the $\pure(\beta,v_1,v_2)$ terms to the left hand side, ``updating'' our initial monomial to
the polynomial
$$\beta^i v_1^j v_2^k -\pure(\beta,v_1,v_2)\stackrel{\mathrm{mod\;}I}{\equiv} 
\pure (\gamma,w_1,w_2) + \mixed. $$
Now the challenge is to find a linear combination $f$ of such polynomials in $\CC[\beta,v_1,v_2]$ such that
the mixed terms cancel, and we are left with 
$$ f \stackrel{\mathrm{mod\;}I}{\equiv} g, \:\;\; \where \;\;\; f \in \CC[\beta, v_1, v_2], \;\; g \in \CC[\gamma, w_1, w_2].$$

The central idea (as in the Euclidean division algorithm) is to put a term order on the mixed terms we are
trying to cancel.  In this way, we can make sure we are cancelling mixed terms in an efficient manner, and the 
cancellation procedure terminates.  Because mixed terms (such as $\beta w_1$) correspond to ``poles'' in the $\gamma$ coordinate chart (such as $\gamma^{-1} w_1$), we use the weighted degree term order
\begin{mapleinput}
\mapleinline{active}{1d}{TP:=wdeg([1,1,1,-1,0,0],[b,v[1],v[2],g,w[1],w[2]]):}{}
\end{mapleinput}
\noindent which keeps track of the degree of the poles in $\gamma$.

Beginning with the superpotential, our algorithm thus consists of the following steps:
\begin{enumerate}
\item Compute transition functions following Ferrari's framework.  This gives an ideal 
$$I = \langle \beta\gamma-1, v_1-\beta^n w_1, v_2-\beta^m w_2 - \partial_{w_1}\til{E}(\gamma,w_1)\rangle \subset \CC[\gamma,w_1,w_2,\beta,v_1,v_2].$$
\item Find a Groebner basis $G$ for the ideal $I$, with respect to a term order $T$.
\item Generate a list $L$ of monomials in $\beta, v_1,$ and $v_2$ (up to some degree).
\item Reduce monomial $L[j]$ mod $I$, using $G$. What you have is
$$L[j] = \beta^i v_1^j v_2^k \stackrel{\mathrm{mod\;}I}{\equiv} \pure(\beta,v_1,v_2) + 
\pure (\gamma,w_1,w_2) + \mixed $$
where ``pure'' and ``mixed'' stand for pure and mixed terms in the appropriate variables.
Bring the $\pure(\beta,v_1,v_2)$ terms over to the LHS to make a polynomial
$$\beta^i v_1^j v_2^k - \pure(\beta,v_1,v_2) \in \CC[\beta,v_1,v_2].$$
      \begin{itemize}
      \item Record this polynomial in the array $F$ as $F[j,1]$.
      \item Record the leading term (with respect to the term order TP) of the ``mixed'' part as $F[j,2]$, and store the leading coefficient as $F[j,3]$.
      \end{itemize}
\item {\bf Reduction routine}
      \begin{itemize}
      \item Cycle through the list of previous polynomials $F[1..j-1,*]$ and cancel leading mixed terms as much as possible.
      \item The result is a new ``updated'' polynomial $F[j,1]$ which is {\em reduced} in the sense that its leading mixed term is as low as possible (with respect to the term order $TP$) due to cancellation with leading mixed terms from previous polynomials.
      \item Reduce the ``updated'' $F[j,1]$ modulo the ideal $I$ to update $F[j,2]$ and $F[j,3]$.  
      \item If the new leading mixed term $F[j,2]$ is 0, we have a {\bf global holomorphic function!}
      \end{itemize}
\item Determine which global holomorphic functions are ``new,'' so that the final list isn't redundant.
      \begin{itemize}
      \item Check that the new global holomorphic function $X_l$ is not in the ring $\CC[X_1,...,X_{l-1}]$
      generated by the previous functions!
      \item To do this we find a Groebner basis for the ideal $\langle X_1,...,X_l \rangle$ and compute partials
      to make sure we can't solve for the new function in terms of the previous ones.
      \end{itemize}
\item Find relations among the global holomorphic functions.  These will determine the (singular) 
geometry of the blow down.
\end{enumerate}

\noindent This completes the proof of Proposition 2.

\subsection{A shortcut}\label{shortcut}

In all of the Intriligator--Wecht cases, we can find weights for the variables $\beta,v_1,v_2$ and $\gamma,w_1,w_2$ such that the transition functions
$$\beta=\gamma^{-1}, \skop v_1=\gamma^{-1}w_1, \skop v_2=\gamma^3 w_2+\partial_{w_1}E(\gamma,w_1),$$
are quasi-homogeneous.  In the $E_6$ case, for example, the transition functions
$$\beta = \gamma^{-1}, \skop v_1 = \gamma^{-1} w_1, \skop v_2 = \gamma^3 w_2 + \gamma^2 w_1^3 +
\gamma^{-1} w_1^2,$$
are quasi-homogeneous if we assign the weights
$$\begin{array}{c|c|c|c|c|c}
\beta & v_1 & v_2 & \gamma & w_1 & w_2 \\
1 & 4 & 7 & -1 & 3 & 10
\end{array}.$$
In particular, this means all elements of the ideal 
$$I = \langle \beta\gamma-1, v_1-\beta w_1, v_2-\gamma^3 w_2 - \gamma^2 w_1^3 - \beta w_1^2 \rangle$$
are quasi-homogeneous in these weights, and all terms in the expression 
$$\beta^i v_1^j v_2^k \stackrel{\mathrm{mod\;}I}{\equiv} \pure(\beta,v_1,v_2) + 
\pure (\gamma,w_1,w_2) + \mixed, $$
will have the same weight.  

This immediately tells us that only combinations of monomials {\em of the same weight} can be used to cancel mixed
terms -- i.e. the global holomorphic functions we build will themselves be quasi-homogeneous.  This observation
cuts computational time immensely, since it means that in the {\bf reduction routine} we need only cycle through lists
of polynomials of the same weight in order to reduce the order of the mixed terms.  In particular, we can run the algorithm in parallel for different weights, restricting ourselves to lists of monomials in $\CC[\beta,v_1,v_2]$ which are all in the same weighted degree.

\subsection{Implementation}

At this point, we have only used the blow-down algorithm in order to find global holomorphic functions in
the Intriligator--Wecht cases.  Our Maple implementation, therefore, takes full advantage of the ``shortcut''~\ref{shortcut}, and has been coded specifically for the 2-matrix model case -- with deformed $\OO(1)\oplus\OO(-3)$ bundle.  The complete Maple code, along with a brief description of all functions and procedures, can be found in the appendix.

All of the global holomorphic functions listed in all of the cases we analyze  (Chapters~\ref{ch:length1},~\ref{ch:length2},~\ref{ch:length3},~\ref{ch:othercases}, and~\ref{ch:hatcases})
have been found using this implementation of the algorithm.

It is important to remember that the algorithm does not tell us when we are done looking for global holomorphic functions.  The only way to check that we have the right blow-down $\M_0$ is by inverting the procedure (blowing back up) and verifying that we recover the original transition functions for the resolved geometry $\Hat{\M}$.  It is because of this that we must still rely on the small resolution techniques developed in Chapter~\ref{ch:smallres}.
\end{onecolumn}

\begin{onecolumn}
\chapter{Length 1}\label{ch:length1}
\setlength{\unitlength}{1 true in}

\begin{figure}[h]

\begin{picture}(3,1)(1.65,.5)
\thicklines
\put(1.9,1){\circle*{.075}}
\put(1.9,1){\line(1,0){.5}}
\put(2.4,1){\circle*{.075}}
\put(2.4,1){\line(1,0){.4625}}
\put(2.9,1){\circle{.075}}
\put(2.9375,1){\line(1,0){.4625}}
\put(3.4,1){\circle*{.075}}
\put(3.4,1){\line(1,0){.5}}
\put(3.9,1){\circle*{.075}}
\put(3.9,1){\line(1,0){.5}}
\put(4.4,1){\circle*{.075}}
\put(1.775,1.05){\makebox(.25,.25){\footnotesize 1}}
\put(2.275,1.05){\makebox(.25,.25){\footnotesize 1}}
\put(2.775,1.05){\makebox(.25,.25){\footnotesize 1}}
\put(3.275,1.05){\makebox(.25,.25){\footnotesize 1}}
\put(3.775,1.05){\makebox(.25,.25){\footnotesize 1}}
\put(4.275,1.05){\makebox(.25,.25){\footnotesize 1}}
\end{picture}

\hspace*{\fill}

\begin{picture}(3,1)(1.65,.5)
\thicklines
\put(1.9,1){\circle*{.075}}
\put(1.9,1){\line(1,0){.5}}
\put(2.4,1){\circle*{.075}}
\put(3.4,1){\circle*{.075}}
\put(3.4,1){\line(1,0){.5}}
\put(3.9,1){\circle*{.075}}
\put(3.9,1){\line(1,0){.5}}
\put(4.4,1){\circle*{.075}}
\put(1.775,1.05){\makebox(.25,.25){\footnotesize 1}}
\put(2.275,1.05){\makebox(.25,.25){\footnotesize 1}}
\put(3.275,1.05){\makebox(.25,.25){\footnotesize 1}}
\put(3.775,1.05){\makebox(.25,.25){\footnotesize 1}}
\put(4.275,1.05){\makebox(.25,.25){\footnotesize 1}}
\end{picture}

\end{figure}

\vspace*{.5in}

\noindent In this chapter we prove the first part of {\bf Theorem 2}: \\
\noindent {\em The resolved geometry corresponding
to the Intriligator--Wecht superpotential 
$$W(x,y) = \dfrac{1}{k+1}x^{k+1} + \dfrac{1}{2}y^2,$$
corresponds to the singular geometry
$$ XY-T(Z^k-T) = 0.$$}
\noindent We will also discover that this potential is geometrically equivalent to 
$$W(x) = \dfrac{1}{k+1}x^{k+1}.$$


\section{The Case $A_k$}

\noindent {\bf The resolved geometry $\Hat{\MM}$}\\
\noindent From the Intriligator--Wecht superpotential
$$W(x,y) = \dfrac{1}{k+1}x^{k+1} + \dfrac{1}{2}y^2,$$
we compute the resolved geometry $\Hat{\MM}$ in terms of transition functions
$$\beta = \gamma^{-1}, \skop v_1 = \gamma^{-1} w_1, \skop v_2 = \gamma^3 w_2 + \gamma^2 w_1^k + w_1.$$
To find the $\PP^1$s, we substitute $w_1(\gamma) = x + \gamma y$ into the $v_2$ transition function
$$v_2(\beta) = \beta^{-3} w_2 + \beta^{-2}(x + \beta^{-1}y)^k + x + \beta^{-1}y.$$
If we choose
$$w_2(\gamma) = \dfrac{x^k - (x+\gamma y)^k}{\gamma},$$
in the $\beta$ chart the section is
$$v_1(\beta) = \beta x + y, \skop v_2(\beta) = \beta^{-2}x^k + \beta^{-1}y + x.$$
This is only holomorphic if
$$x^k = y = 0,$$
and so we have a single $\PP^1$ located at
$$w_1(\gamma) = w_2(\gamma) = 0, \skop v_1(\beta) = v_2(\beta) = 0.$$
This is exactly what we expect from computing critical points of the superpotential
$$\d W = x^k dx + y dy = 0.$$
\skp

\noindent{\bf Global holomorphic functions}\\
\noindent The transition functions
$$\beta = \gamma^{-1}, \skop v_1 = \gamma^{-1} w_1, \skop v_2 = \gamma^3 w_2 + \gamma^2 w_1^k + w_1,$$
are quasi-homogeneous if we assign the weights
$$\begin{array}{c|c|c|c|c|c}
\beta & v_1 & v_2 & \gamma & w_1 & w_2 \\
k-1 & k+1 & 2 & 1-k & 2 & 3k-1
\end{array}.$$
The global holomorphic functions will thus necessarily be quasi-homogeneous in these weights.
We find the following global holomorphic functions:
$$\begin{array}{c|ccccc}
2 & y_1 &=& v_2 &=& \gamma^3 w_2 + \gamma^2 w_1^k + w_1\\
k+1 & y_2 &=& \beta v_2 - v_1 &=& \gamma^2 w_2 + \gamma w_1^k\\
2k & y_3 &=& \beta^2 v_2-\beta v_1 &=& \gamma w_2 + w_1^k\\
3k-1 & y_4 &=& v_2^{k-1} v_1 - \beta^3 v_2 + \beta^2 v_1 &=&
\end{array}$$
These are the first 4 ``distinct'' functions produced by our algorithm, 
in the sense that none is contained in the ring generated by the other 3.
\skp

\noindent {\bf The singular geometry $\MM_0$}\\
\noindent We conjecture that the ring of 
global holomorphic functions is generated by $y_1,y_2,y_3$ and $y_4$, 
subject to the degree $4k$ relation
$$\MM_0: \;\;\; y_2 y_4 + y_3^2 + y_2^2 y_1^{k-1} - y_3 y_1^k = 0.$$
The functions $y_i$ give us a blow-down map whose image $\MM_0$
has an isolated $A_k$ singularity. To see this, consider the 
change of variables
$$\til{y}_4 = y_4 + y_2 y_1^{k-1}= \beta v_2^k - \beta^3 v_2 + \beta^2 v_1.$$
Note that like $y_4$, $\til{y}_4$ is also quasi-homogeneous of degree $3k-1$. 
The functions $y_1,y_2, y_3$ and $\til{y}_4$ now satisfy the simpler relation
$$\MM_0: \;\;\; y_2 \til{y}_4 + y_3(y_3 - y_1^k) = 0.$$
\skp

\noindent {\bf The blowup}\\
\noindent We now verify that we have identified the right singular space $\MM_0$ by inverting
the blow-down.  In the $\beta$ and $\gamma$ charts we find
$$\begin{array}{ccc|ccc}
\beta &=& y_3/y_2 & \gamma &=& y_2/y_3\\
v_1 &=& \beta y_1-y_2 = (y_3y_1 - y_2^2)/y_2 & w_1 &=& y_1 - \gamma y_2 = (y_1 y_3-y_2^2)/y_3\\
v_2 &=& y_1 & w_2 &=& \beta(y_3-w_1^k)\\
&&& &=& -\til{y}_4 + \dfrac{y_1^k - (y_1-\gamma y_2)^k}{\gamma}.
\end{array}$$
This suggests that we should blow up 
$$y_2 = y_3 = 0,$$
for the small resolution of $\MM_0$.  We introduce $\PP^1$ coordinates $[\beta:\gamma]$ such that
$\beta y_2 = \gamma y_3.$  The blow up in each chart is then
$$\begin{array}{c|c}
(\gamma = 1) & (\beta = 1)\\
&\\
y_3 = \beta y_2 & y_2 = \gamma y_3\\
\til{y_4} = \beta(y_1^k - \beta y_2) & y_3 = y_1^k - \gamma \til{y_4}\\
&\\
\mathrm{coords:\;\;} (y_1, y_2, \beta) & \mathrm{coords:\;\;}(y_1, \til{y_4}, \gamma)
\end{array}$$
\skp

\noindent{\bf Transition functions}\\ 
\noindent The transition functions between the $\beta$ and $\gamma$ charts are
$$\beta = \gamma^{-1},\skop y_1 = y_1, \skop y_2 = \gamma(y_1^k - \gamma \til{y}_4)
= -\gamma^2 \til{y}_4 + \gamma y_1^k.$$
Note that for $k > 1$, this is an $\OO\oplus\OO(-2)$ bundle over the exceptional $\PP^1$, and corresponds
to a superpotential with a single field ($M=1$):
$$W(x) = \dfrac{1}{k+1}x^{k+1}.$$
(For $k = 1$ the bundle is actually $\OO(-1)\oplus\OO(-1)$ and $W=0$.)

In terms of the original coordinates, the transition functions become
$$\beta = \gamma^{-1},\skop v_2 = \gamma^3 w_2 + \gamma^2 w_1^k + w_1, \skop
\beta v_2 - v_1 = \gamma^2 w_2 + \gamma w_1^k.$$
Substituting the second transition function into the third reveals $v_1 = \gamma^{-1} w_1,$ and so
we recover our original transition functions
$$\beta = \gamma^{-1}, \skop v_1 = \gamma^{-1} w_1, \skop v_2 = \gamma^3 w_2 + \gamma^2 w_1^k + w_1,$$
which define an $\OO(1)\oplus\OO(-3)$ bundle deformed by the two field ($M=2$) superpotential
$$W(x,y) = \dfrac{1}{k+1} x^{k+1} + \dfrac{1}{2}y^2.$$

\section{A puzzle}

\noindent \textbf{The problem}\\
\noindent We saw the $A_k$ case in Section~\ref{sec:A_k}, with geometry $\Hat{\M}$ given by the superpotential 
$$W(x) = \dfrac{1}{k+1}x^{k+1},$$
and hence corresponding to a deformed $\OO\oplus\OO(-2)$ bundle over the $\PP^1$, with one field. 
However, Intriligator and Wecht identify 
$$W(x,y) = \dfrac{1}{k+1}x^{k+1} + \dfrac{1}{2}y^2,$$
as corresponding to an $A_k$-type singularity, with an extra field $y$ which requires
that the transition functions look like
$$\beta = \gamma^{-1},\skop v_1 = \gamma^{-1} w_1, \skop v_2 = \gamma^3 w_2 + \gamma^2 w_1^k + w_1.$$
In particular, the geometry $\Hat{\M}$ looks like that of an $\OO(1)\oplus\OO(-3)$ bundle over the
exceptional $\PP^1$!  What's going on here?\skp

\noindent \textbf{Resolution of the problem}\\
\noindent For $n=1$ and $k=0$, Proposition 1 tells us that the superpotential
$$W(x,y) = \dfrac{1}{2} y^2$$
changes the bundle 
$$\OO(1)\oplus\OO(-3) \flecha \OO\oplus\OO(-2).$$
Hence the extra field $y$ from the Intriligator--Wecht potential (with purely quadratic contribution to the superpotential) can be ``integrated out.''  Its effect is to change the bundle for $A_k$ from $\OO(1)\oplus\OO(-3)$,
which is necessary for a two-field description, to reveal the true underlying $\OO\oplus\OO(-2)$ structure.
In other words, 
$$W(x,y) = \dfrac{1}{k+1}x^{k+1} + \dfrac{1}{2}y^2 \skop \and \skop W(x) = \dfrac{1}{k+1}x^{k+1}$$
are geometrically equivalent.\skp

\noindent {\bf Beyond Proposition 1}\\
\noindent Note that this is not just a straightforward application of Proposition 1, which implies that $W(x,y) = \dfrac{1}{2}y^2$ and $W(x) = 0$ are geometrically equivalent.  Beginning with transition functions for the 
Intriligator--Wecht $A_k$ superpotential $W(x,y) = \dfrac{1}{k+1}x^{k+1} + \dfrac{1}{2}y^2,$
$$\beta = \gamma^{-1}, \skop v_1 = \gamma^{-1} w_1, \skop v_2 = \gamma^3 w_2 + \gamma^2 w_1^k + w_1,$$
the change of coordinates suggested in the proof of Proposition 1 
$$\begin{array}{cccccc}
\til{v_1} &=& v_2, & \til{w_1} &=& w_1 + \gamma^3 w_2, \\
\til{v_2} &=& -v_1 + \beta v_2, & \til{w_2} &=& w_2,
\end{array}$$
does not yield the appropriate new transition functions.  Instead, the more complicated change of coordinates
$$\begin{array}{cccccc}
\til{v_1} &=& v_2, & \til{w_1} &=& w_1 + \gamma^3 w_2 + \gamma^2 w_1^k, \\
\til{v_2} &=& -v_1 + \beta v_2, & \til{w_2} &=& w_2 - \gamma^{-1}\left[(\gamma^3w_2 + \gamma^2w_1^k + w_1)^k-w_1^k\right],
\end{array}$$
is needed to give new transition functions
$$\beta = \gamma^{-1}, \skop \til{v_1} = \til{w_1}, \skop \til{v_2} = \gamma^2\til{w_2} + \gamma \til{w_1}^k,$$
corresponding to the superpotential $W(x) = \dfrac{1}{k+1}x^{k+1}.$

It would be interesting to try to generalize Proposition 1 to include examples such as this, where there
are additional terms in the superpotential besides the quadratic pieces which suggest a change in bundle structure.  Trying to understand what makes the change of coordinates possible in this case may give hints as to how the geometric picture for RG flow might be extended.  The ultimate goal would be to understand how ``integrating out'' the $y$ coordinate in a potential of the form $W(x,y) = f(x,y) + y^2$ affects other terms involving $y$.
\end{onecolumn}

\begin{onecolumn}
\chapter{Length 2}\label{ch:length2}
\setlength{\unitlength}{1 true in}

\begin{figure}[h]

\begin{picture}(3,1)(1.65,.5)
\thicklines
\put(1.9,1){\circle*{.075}}
\put(1.9,1){\line(1,0){.5}}
\put(2.4,1){\circle*{.075}}
\put(2.4,1){\line(1,0){.4625}}
\put(2.9,1){\circle{.075}}
\put(2.9375,1){\line(1,0){.4625}}
\put(3.4,1){\circle*{.075}}
\put(3.4,1){\line(1,0){.5}}
\put(3.9,1){\circle*{.075}}
\put(3.9,1){\line(1,0){.5}}
\put(4.4,1){\circle*{.075}}
\put(4.4,1){\line(3,4){.3}}
\put(4.7,1.4){\circle*{.075}}
\put(4.4,1){\line(3,-4){.3}}
\put(4.7,.6){\circle*{.075}}
\put(1.775,1.05){\makebox(.25,.25){\footnotesize 1}}
\put(2.275,1.05){\makebox(.25,.25){\footnotesize 2}}
\put(2.775,1.05){\makebox(.25,.25){\footnotesize 2}}
\put(3.275,1.05){\makebox(.25,.25){\footnotesize 2}}
\put(3.775,1.05){\makebox(.25,.25){\footnotesize 2}}
\put(4.2,1.05){\makebox(.25,.25){\footnotesize 2}}
\put(4.575,1.45){\makebox(.25,.25){\footnotesize 1}}
\put(4.575,.3){\makebox(.25,.25){\footnotesize 1}}
\end{picture}

\bigskip
\bigskip

\begin{picture}(3,1)(1.65,.5)
\thicklines
\put(1.9,1){\circle*{.075}}
\put(1.9,1){\line(1,0){.5}}
\put(2.4,1){\circle*{.075}}
\put(3.4,1){\circle*{.075}}
\put(3.4,1){\line(1,0){.5}}
\put(3.9,1){\circle*{.075}}
\put(3.9,1){\line(1,0){.5}}
\put(4.4,1){\circle*{.075}}
\put(4.4,1){\line(3,4){.3}}
\put(4.7,1.4){\circle*{.075}}
\put(4.4,1){\line(3,-4){.3}}
\put(4.7,.6){\circle*{.075}}
\put(1.775,1.05){\makebox(.25,.25){\footnotesize 1}}
\put(2.275,1.05){\makebox(.25,.25){\footnotesize 2}}
\put(3.275,1.05){\makebox(.25,.25){\footnotesize 2}}
\put(3.775,1.05){\makebox(.25,.25){\footnotesize 2}}
\put(4.2,1.05){\makebox(.25,.25){\footnotesize 2}}
\put(4.575,1.45){\makebox(.25,.25){\footnotesize 1}}
\put(4.575,.3){\makebox(.25,.25){\footnotesize 1}}
\end{picture}

\end{figure}

\vspace*{.5in}

\noindent In this chapter we prove the second part of {\bf Theorem 2}: \\
\noindent {\em The resolved geometry corresponding
to the Intriligator--Wecht superpotential 
$$W(x,y) = \dfrac{1}{k+1}x^{k+1}+xy^2$$
corresponds to the singular geometry
$$ X^2 - Y^2Z + T(Z^{k/2} - T)^2 = 0.$$}

\pagebreak 

\section{The Case $D_{k+2}$}

\noindent \textbf{The resolved geometry} $\Hat{\MM}$\\
\noindent From the Intriligator-Wecht superpotential
$$W(x,y) = \dfrac{1}{k+1}x^{k+1} + xy^2,$$
we compute the resolved geometry $\Hat{\MM}$ in terms of transition functions
$$\beta = \gamma^{-1}, \skop v_1 = \gamma^{-1} w_1, \skop v_2 = \gamma^3 w_2 + \gamma^2 w_1^k + w_1^2.$$
To find the $\PP^1$s, we substitute $w_1(\gamma) = x + \gamma y$ into the $v_2$ transition function
$$v_2(\beta) = \beta^{-3} w_2 + \beta^{-2}(x + \beta^{-1}y)^k + (x + \beta^{-1}y)^2.$$
If we choose
$$w_2(\gamma) = \dfrac{x^k - (x+\gamma y)^k}{\gamma},$$
in the $\beta$ chart the section is
$$v_1(\beta) = \beta x + y, \skop v_2(\beta) = \beta^{-2}(x^k+y^2) + \beta^{-1}(2xy) + x^2.$$
This is only holomorphic if
$$x^k+y^2 = 2xy = 0,$$
and so we have a single $\PP^1$ located at
$$w_1(\gamma) = w_2(\gamma) = 0, \skop v_1(\beta) = v_2(\beta) = 0.$$
This is exactly what we expect from computing critical points of the superpotential
$$\d W = (x^k+y^2) dx + 2xy \;dy = 0.$$
\skp

\noindent \textbf{Global holomorphic functions}\\
\noindent The transition functions
$$\beta = \gamma^{-1}, \skop v_1 = \gamma^{-1} w_1, \skop v_2 = \gamma^3 w_2 + \gamma^2 w_1^k + w_1^2,$$
are quasi-homogeneous if we assign the weights
$$\begin{array}{c|c|c|c|c|c}
\beta & v_1 & v_2 & \gamma & w_1 & w_2 \\
k-2 & k & 4 & 2-k & 2 & 3k-2
\end{array}.$$
The global holomorphic functions will thus necessarily be quasi-homogeneous in these weights.
We find the following global holomorphic functions for $k$ even:
$$\begin{array}{c|ccccc}
4 & X &=& v_2 &=& \gamma^3 w_2 + \gamma^2 w_1^k + w_1^2\\
2k & Y &=& \beta^2 v_2 - v_1^2 &=& \gamma w_2 + w_1^k\\
3k-2 & Z &=& \beta(X^{k/2}-Y) &=& \gamma^{-1}(X^{k/2}-Y)\\
3k & U &=& v_1(X^{k/2}-Y) &=& \gamma^{-1}w_1(X^{k/2}-Y),
\end{array}$$
and a similar set of global holomorphic functions for $k$ odd:
$$\begin{array}{c|ccccc}
4 & X &=& v_2 &=& \gamma^3 w_2 + \gamma^2 w_1^k + w_1^2\\
2k & Y &=& \beta^2 v_2 - v_1^2 &=& \gamma w_2 + w_1^k\\
3k-2 & Z &=& v_1 X^{\lfloor k/2 \rfloor}-\beta Y &=& 
\gamma^{-1}(w_1 X^{\lfloor k/2 \rfloor}- Y)\\
3k & U &=& \beta X^{\lceil k/2 \rceil}-v_1 Y &=&
\gamma^{-1}( X^{\lceil k/2 \rceil}-w_1 Y).
\end{array}$$
\skp

\noindent \textbf{The singular geometry} $\MM_0$\\
\noindent We conjecture that the ring of 
global holomorphic functions is generated by $X,Y,Z$ and $U$, 
subject to the degree $6k$ relation
$$\begin{array}{ccccccc}
\MM_0:& \;\;\; & U^2 - XZ^2 +Y(X^{k/2}-Y)^2 &=& 0,& \;\;\; & k \; \mathrm{even},\\
\MM_0:& \;\;\; & U^2 - XZ^2 -Y(X^k - Y^2) &=& 0,& \;\;\; & k \; \mathrm{odd}.
\end{array}$$
The functions $X,Y,Z$ and $U$ give us a blow-down map whose image $\MM_0$
has an isolated $D_{k+2}$ singularity.\skp

\noindent {\bf Review of length 2 blowup}\\
\noindent Before doing the blowup to see that we have the right blow down, we review some results from Chapter~\ref{ch:smallres}.  There we found
small resolutions of length 2 singularities by using deformations of matrix factorizations for  $D_{n+2}$
surface singularities.  

Recall (dropping tildes!) the deformed $D_{n+2}$ equation~\eqref{eq:deformedD}
$$0 = X^2 + Y^2Z - h'(ZP''^2+Q''^2)+2\delta'(YQ'' +(-1)^{m+1} XP''),$$
where for $t=0$ we have $\delta'=0$ and
$$\begin{array}{c|ccc}
m \mathrm{\;even} & h'(Z)=Z^{n-m}, & P''(Z) = i^m Z^{m/2}, & Q''(Z)=0.\\
m \mathrm{\;odd} & h'(Z)=Z^{n-m}, & P''(Z) = 0, & Q''(Z)= i^{m+1}Z^{(m+1)/2}.
\end{array}$$
The blowup was given by the equation for the Grassmanian $G(2,4)\subseteq \PP^5$
$$\alpha^2 - \varphi^2Z + (-1)^{m+1}h'\varepsilon^2 + 2i^{m+1}\delta'\varepsilon\varphi = 0,$$
in terms of Pl\"ucker coordinates
\begin{eqnarray*}
\alpha &=& iXY - i^{2m+1}h'P''Q''\\
\varepsilon &=&  i^{m+3}XP''+i^{-m+1}YQ''\\
\varphi &=& Y^2 - h'P''^2. 
\end{eqnarray*}
The interesting charts were $\varphi = 1$ and $\varepsilon = 1$, with transition functions
\begin{eqnarray}\label{eq:tf}
\varphi_2 &=& \varepsilon_1^{-1},\\
\alpha_2 &=& \varepsilon_1^{-1}\alpha_1,\\
Z &=& (-1)^{m+1} \varepsilon_1^2 h'(Z,t) + \alpha_1^2 + 2i^{m+1}\delta'\varepsilon_1,\\
t &=& t.
\end{eqnarray}
\skp

\noindent \textbf{The blowup}\\
\noindent As usual, we verify that we have the right singular space by inverting the blow down. 
A simple change of coordinates shows our $\M_0$ equations to be
\begin{eqnarray*}
 X^2 + Y^2Z - T(Z^{k/2}-T)^2 &=& 0 \skop (k \mathrm{\;\;even} )\\
 X^2 + Y^2Z + T(Z^k - T^2) &=& 0 \skop (k \mathrm{\;\;odd}). 
\end{eqnarray*}

\noindent {\bf $\mathbf{k}$ even}\\
\noindent For $k$ even, the equation for $\M_0$ corresponds to
$$\delta'=0, \skop h'=T, \skop P''=0, \skop Q''=i^k(Z^{k/2}-T), \skop n=m=k-1,$$
and Pl\"ucker coordinates
$$\alpha = iXY, \skop \varepsilon = -Y(Z^{k/2}-T), \skop \varphi = Y^2.$$
The connection with the original transition function coordinates is
$$\begin{array}{c|ccccc}
4 & Z &=& v_2 &=& \gamma^3 w_2 + \gamma^2 w_1^k + w_1^2\\
2k & T &=& \beta^2 v_2 - v_1^2 &=& \gamma w_2 + w_1^k\\
3k-2 & Y &=& \beta(Z^{k/2}-T) &=& \gamma^{-1}(Z^{k/2}-T)\\
3k & iX &=& v_1(Z^{k/2}-T) &=& \gamma^{-1}w_1(Z^{k/2}-T).
\end{array}$$
Note that
$$\begin{array}{ccc}
\beta &=& \dfrac{Y}{Z^{k/2}-T} = -\dfrac{\varphi}{\varepsilon} = -\varphi_2 \\
v_1 &=& \dfrac{iX}{Z^{k/2}-T} = -\dfrac{\alpha}{\varepsilon} = -\alpha_2\\
v_2 &=& Z \\
\gamma &=& \dfrac{Z^{k/2} - T}{Y} = -\dfrac{\varepsilon}{\varphi} = -\varepsilon_1\\
w_1 &=& \dfrac{iX}{Y} = \dfrac{\alpha}{\varphi} = \alpha_1 \\
w_2&=& \beta(T-w_1^k) 
= -\varepsilon_2^{-1}(T-\alpha_2^k)
\end{array}$$
In particular, the transition functions~\eqref{eq:tf} become
\begin{eqnarray*}
\beta &=& \gamma^{-1},\\
v_1 &=& \gamma^{-1}w_1,\\
v_2 &=& \gamma^2 T + w_1^2 = \gamma^3 w_2 + \gamma^2 w_1^k + w_1^2.
\end{eqnarray*}

\noindent {\bf $\mathbf{k}$ odd}\\
Here the equation for $\M_0$ corresponds to
$$\delta' = 0, \skop h' = -T, \skop P''=i^{k-1}Z^{(k-1)/2}, \skop Q''=iT, \skop n=m=k-1,$$
and Pl\"ucker coordinates
$$\alpha = iXY + i^{3k+3}T^2Z^{(k-1)/2}, \skop \varepsilon = i^{2k+1}XZ^{(k-1)/2} + i^{3k+3}TY, \skop
\varphi = Y^2 + i^{2k+2}TZ^{k-1}.$$
$$\begin{array}{c|ccccc}
4 & Z &=& v_2 &=& \gamma^3 w_2 + \gamma^2 w_1^k + w_1^2\\
2k & T &=& \beta^2 v_2 - v_1^2 &=& \gamma w_2 + w_1^k\\
3k-2 & Y &=& v_1 Z^{(k-1)/2}-\beta T &=& 
\gamma^{-1}(w_1 Z^{(k-1)/2}- T)\\
3k & iX &=& \beta Z^{(k-1)/2}-v_1 T &=&
\gamma^{-1}( Z^{(k-1)/2}-w_1 T).
\end{array}$$
Note that
\begin{eqnarray*}
iX + Y &=& (Z^{(k-1)/2}-T)(\beta + v_1)\\
\alpha + \phi &=& (iX + Y)Y + (Z^{(k-1)/2}+i^{k+1}T)Z^{(k-1)/2}T,
\end{eqnarray*}
and recall that for $m$ even
$$i\alpha_1 P'' + Q'' = i^{m+3} Y \varepsilon_1.$$

\section{Transition functions for more general length 2 blowup}

\noindent \textbf{Generalization to} $\mathbf{P''=0.}$\\
\noindent Note that $P''=0$ implies $m$ odd.  The deformed equation simplifies to
$$0 = X^2 + Y^2Z - h'Q''^2+2\delta'YQ'',$$
with Pl\"ucker coordinates
$$\alpha = iXY, \skop \varepsilon = i^{-m+1}YQ'', \skop \varphi = Y^2.$$
Following the example of the Intriligator--Wecht case for $k$ even (and $m$ odd), we conjecture that the connection with the transition function coordinates 
in the blowup is
$$\begin{array}{c|ccccc}
4 & Z &=& v_2 &=& \gamma^3 w_2 + \gamma^2 w_1^k + w_1^2 \\
2k & T &=& \beta^2 v_2 - v_1^2 &=& \gamma w_2 + w_1^k \\
3k-2 & Y &=& i^{-m-1} \beta Q'' &=& i^{-m-1}\gamma^{-1}Q''\\
3k & iX &=& i^{-m-1} v_1 Q'' &=& i^{-m-1}\gamma^{-1}w_1 Q'',
\end{array}$$
and we can express the upstairs coordinates as
$$\begin{array}{ccc}
\beta &=& i^{m+1}Y/Q'' = -\varphi/\varepsilon = -\varphi_2 \\
v_1 &=& i^{m+2}X/Q'' = -\alpha/\varepsilon = -\alpha_2\\
v_2 &=& Z \\
\gamma &=& i^{-m-1}Q''/Y = -\varepsilon/\varphi = -\varepsilon_1\\
w_1 &=& iX/Y = \alpha/\varphi = \alpha_1 \\
w_2&=& \beta((-1)^{m+1}h'-w_1^k-2i^{m+1}\beta\delta') \\
&=& -\varepsilon_2^{-1}((-1)^{m+1}h'-\alpha_2^k+2i^{m+1}\varepsilon_2^{-1}\delta'). 
\end{array}$$
\noindent From the original deformed $D_{n+2}$ equation 
$$0 = X^2 + Y^2Z - h'(ZP''^2+Q''^2)+2\delta'(YQ'' +(-1)^{m+1} X P''),$$
we know that ($P''=0$)
\begin{eqnarray*}
v_2 = Z &=& \dfrac{h'Q''^2 - X^2 - 2\delta'YQ''}{Y^2}\\
&=& \gamma^2 h' + w_1^2 -2\delta' i^{m+1}\gamma\\
&=& \gamma^2 (\gamma w_2 + w_1^k+2i^{m+1}\gamma^{-1}\delta') + w_1^2 -2\delta' i^{m+1}\gamma\\
&=& \gamma^3 w_2 + \gamma^2 w_1^k + w_1^2.
\end{eqnarray*}
Equivalently, the equation for the Grassmannian  $G(2,4)\subseteq \PP^5$  
$$(-Z\varphi-2i^{m-1}\delta'\varepsilon)\phi + (-1)^{m+1}h'\varepsilon^2 + \alpha^2 = 0,$$
tells us that
\begin{eqnarray*}
v_2 = Z &=& \dfrac{\alpha^2+(-1)^{m+1}h'\varepsilon^2+2i^{m+1}\delta'\varepsilon\varphi}{\varphi^2}\\
&=& \left(\dfrac{\alpha}{\varphi}\right)^2+(-1)^{m+1}h'\left(\dfrac{\varepsilon}{\varphi}\right)^2+
2i^{m+1}\delta'\left(\dfrac{\varepsilon}{\varphi}\right)\\
&=& w_1^2+(-1)^{m+1}h'\gamma^2 -2 i^{m+1}\gamma\delta'\\
&=& \gamma^2 (\gamma w_2 + w_1^k+2i^{m+1}\gamma^{-1}\delta') + w_1^2 -2\delta' i^{m+1}\gamma\\
&=& \gamma^3 w_2 + \gamma^2 w_1^k + w_1^2.
\end{eqnarray*}
Notice that this upstairs transition functions does not depend on $h'$ or $Q''$!  However, since 
$$\delta' = \delta'(T),\skop \and\skop h'=h'(Z,T),$$ it is not clear that we can in general write 
this as transition functions between two $\CC^3$ charts.

Assuming that these coordinates suffice to describe the blowup, we get exactly the 
Intriligator--Wecht transition functions:
$$\beta = \gamma^{-1}, \skop v_1 = \gamma^{-1}w_1, \skop v_2 = \gamma^3 w_2 + \gamma^2 w_1^k + w_1^2.$$
In other words, the actual values for $m$, $Q''$, $\delta'$ and 
$h'$ don't matter, as long as the corresponding upstairs coordinates are legitimate (and this is
where they may matter).

\noindent {\bf Full generality}\\
\noindent We have Pl\"ucker coordinates
\begin{eqnarray*}
\alpha &=& iXY - i^{2m+1}h'P''Q''\\
\varepsilon &=& i^{m+3}XP''+i^{-m+1}YQ''\\
\varphi &=& Y^2 - h'P''^2,
\end{eqnarray*}
and equation for the Grassmanian $G(2,4)\subseteq \PP^5$ 
$$(-Z\varphi-2i^{m-1}\delta'\varepsilon)\varphi + (-1)^{m+1}h'\varepsilon^2 + \alpha^2 = 0.$$
It is easy to check that this is exactly equivalent to the original deformed equation
$$0 = X^2 + Y^2Z - h'(ZP''^2+Q''^2)+2\delta'(YQ'' +(-1)^{m+1} X P'').$$
Definining upstairs coordinates
$$\begin{array}{ccc|ccc}
\beta &=& -\varphi/\varepsilon & \gamma &=& -\varepsilon/\varphi \\
v_1 &=& -\alpha/\varepsilon & w_1 &=& \alpha/\varphi \\ 
v_2 &=& Z & w_2 &=& \beta[(-1)^{m+1}h'-w_1^k-2i^{m+1}\beta\delta']
\end{array}$$
and using the equation for $G(2,4)$ yields transition function
\begin{eqnarray*}
v_2 = Z &=& \dfrac{\alpha^2+(-1)^{m+1}h'\varepsilon^2+2i^{m+1}\delta'\varepsilon\varphi}{\varphi^2}\\
&=& \left(\dfrac{\alpha}{\varphi}\right)^2+(-1)^{m+1}h'\left(\dfrac{\varepsilon}{\varphi}\right)^2+
2i^{m+1}\delta'\left(\dfrac{\varepsilon}{\varphi}\right)\\
&=& w_1^2+(-1)^{m+1}h'\gamma^2 -2 i^{m+1}\gamma\delta'\\
&=& \gamma^2 (\gamma w_2 + w_1^k+2i^{m+1}\gamma^{-1}\delta') + w_1^2 -2\delta' i^{m+1}\gamma\\
&=& \gamma^3 w_2 + \gamma^2 w_1^k + w_1^2.
\end{eqnarray*}

Note that the choice of $k$ was arbitrary in the definition of the upstairs coordinates.
Moreover, we don't know if these are enough coordinates to fully describe the resolved geometry.
However, in cases where it is enough (and this may depend on $k$!), 
we see that the upstairs transition functions are all
the same as in the Intriligator--Wecht case:
$$\beta = \gamma^{-1}, \skop v_1 = \gamma^{-1}w_1, \skop v_2 = \gamma^3 w_2 + \gamma^2 w_1^k + w_1^2.$$
The question is whether we can solve for $X, Y$ and $T$ in the $\varepsilon=1$ chart ($\beta,v_1,v_2$) and similarly
can we restrict to $\gamma,w_1$ and $w_2$ in the $\varphi = 1$ chart.  In the $k$ even Intriligator--Wecht case, with
$\delta'=P''=0$, and $h'=T$, this certainly seems to be the case.

\end{onecolumn}

\begin{onecolumn}
\chapter{Length 3}\label{ch:length3}
\setlength{\unitlength}{1 true in}

\begin{figure}[h]

\begin{picture}(3,1)(1.9,.5)
\thicklines
\put(1.9,1){\circle*{.075}}
\put(1.9,1){\line(1,0){.5}}
\put(2.4,1){\circle*{.075}}
\put(2.4,1){\line(1,0){.4625}}
\put(2.9,1){\circle*{.075}}
\put(2.9,1){\line(0,-1){.4625}}
\put(2.9,.5){\circle{.075}}
\put(2.9375,1){\line(1,0){.4625}}
\put(3.4,1){\circle*{.075}}
\put(3.4,1){\line(1,0){.5}}
\put(3.9,1){\circle*{.075}}
\put(3.9,1){\line(1,0){.5}}
\put(4.4,1){\circle*{.075}}
\put(4.4,1){\line(1,0){.5}}
\put(4.9,1){\circle*{.075}}
\put(1.775,1.05){\makebox(.25,.25){\footnotesize 2}}
\put(2.275,1.05){\makebox(.25,.25){\footnotesize 4}}
\put(2.775,1.05){\makebox(.25,.25){\footnotesize 6}}
\put(2.925,.375){\makebox(.25,.25){\footnotesize 3}}
\put(3.275,1.05){\makebox(.25,.25){\footnotesize 5}}
\put(3.775,1.05){\makebox(.25,.25){\footnotesize 4}}
\put(4.275,1.05){\makebox(.25,.25){\footnotesize 3}}
\put(4.775,1.05){\makebox(.25,.25){\footnotesize 2}}
\end{picture}

\bigskip
\bigskip

\begin{picture}(3,1)(1.9,.5)
\thicklines
\put(1.9,1){\circle*{.075}}
\put(1.9,1){\line(1,0){.5}}
\put(2.4,1){\circle*{.075}}
\put(2.4,1){\line(1,0){.4625}}
\put(2.9,1){\circle*{.075}}
\put(2.9375,1){\line(1,0){.4625}}
\put(3.4,1){\circle*{.075}}
\put(3.4,1){\line(1,0){.5}}
\put(3.9,1){\circle*{.075}}
\put(3.9,1){\line(1,0){.5}}
\put(4.4,1){\circle*{.075}}
\put(4.4,1){\line(1,0){.5}}
\put(4.9,1){\circle*{.075}}
\put(1.775,1.05){\makebox(.25,.25){\footnotesize 2}}
\put(2.275,1.05){\makebox(.25,.25){\footnotesize 4}}
\put(2.775,1.05){\makebox(.25,.25){\footnotesize 6}}
\put(3.275,1.05){\makebox(.25,.25){\footnotesize 5}}
\put(3.775,1.05){\makebox(.25,.25){\footnotesize 4}}
\put(4.275,1.05){\makebox(.25,.25){\footnotesize 3}}
\put(4.775,1.05){\makebox(.25,.25){\footnotesize 2}}
\end{picture}

\end{figure}

\setlength{\unitlength}{1pt}

\vspace*{.5in}

\noindent In this chapter we prove the third part of {\bf Theorem 2}: \\
\noindent {\em The resolved geometry corresponding
to the Intriligator--Wecht superpotential 
$$W(x,y) = \dfrac{1}{3}y^3+yx^3$$
corresponds to the singular geometry
$$ X^2 - Y^3 + Z^5 + 3TYZ^2 + T^3Z = 0. $$}

\pagebreak

\section{The Case $E_7$}

\noindent \textbf{The resolved geometry} $\Hat{\MM}$\\
From the Intriligator-Wecht superpotential
$$W(x,y) = \dfrac{1}{3}y^3 + yx^3,$$
we compute the resolved geometry $\Hat{\MM}$ in terms of transition functions
$$\beta = \gamma^{-1}, \skop v_1 = \gamma^{-1} w_1, \skop v_2 = \gamma^3 w_2 + \gamma w_1^3 + 
\gamma^{-1} w_1^2.$$
To find the $\PP^1$s, we substitute $w_1(\gamma) = x + \gamma y$ into the $v_2$ transition function
$$v_2(\beta) = \beta^{-3} w_2 + \beta^{-1}(x + \beta^{-1}y)^3 + \beta(x + \beta^{-1}y)^2.$$
If we choose
$$w_2(\gamma) = \dfrac{x^3 + 3\gamma x^2y - (x+\gamma y)^3}{\gamma^2},$$
in the $\beta$ chart the section is
$$v_1(\beta) = \beta x + y, \skop v_2(\beta) = \beta^{-2}(3x^2 y) + \beta^{-1}(x^3 + y^2)+2xy+\beta x^2.$$
This is only holomorphic if
$$3x^2y = x^3+y^2 = 0,$$
and so we have a single $\PP^1$ located at
$$w_1(\gamma) = w_2(\gamma) = 0, \skop v_1(\beta) = v_2(\beta) = 0.$$
This is exactly what we expect from computing critical points of the superpotential
$$\d W = 3x^2y\; dx + (y^2+x^3)dy = 0.$$
\skp

\noindent \textbf{Global holomorphic functions}\\
\noindent The transition functions
$$\beta = \gamma^{-1}, \skop v_1 = \gamma^{-1} w_1, \skop v_2 = \gamma^3 w_2 + \gamma w_1^3 +
\gamma^{-1} w_1^2,$$
are quasi-homogeneous if we assign the weights
$$\begin{array}{c|c|c|c|c|c}
\beta & v_1 & v_2 & \gamma & w_1 & w_2 \\
1 & 3 & 5 & -1 & 2 & 8
\end{array}.$$
The global holomorphic functions will thus necessarily be quasi-homogeneous in these weights.
We find the following global holomorphic functions:
$$\begin{array}{c|ccccc}
6 & X &=& \beta v_2-v_1^2 && \\
8 & Y &=& v_1v_2 - \beta^2X && \\
10 & Z &=& v_2^2-\beta v_1 X && \\
15 & F &=& v_2^3 - 2v_1^3X + (\beta^3-3v_1)X^2. && 
\end{array}$$
\skp

\noindent \textbf{The singular geometry} $\MM_0$\\
We conjecture that the ring of 
global holomorphic functions is generated by $X,Y,Z$ and $F$, 
subject to the degree $30$ relation
$$\MM_0: \;\;\; F^2 -Z^3 + X^5 + 3X^2YZ + XY^3 = 0.$$
Do the functions $X,Y,Z$ and $F$ give us a blow-down map whose image $\MM_0$
has an isolated $E_7$ singularity?\skp

\noindent \textbf{The blowup}\\
\noindent We now verify that we have identified the right singular space $\MM_0$ by inverting
the blow-down.  In the $\beta$ and $\gamma$ charts we find
$$\begin{array}{ccc|cccc}
\beta &=& \dfrac{Y^2+XZ}{F} & \gamma &=& \dfrac{F}{Y^2+XZ}&\\
&&&&&&\\
v_1 &=& \dfrac{X^3 + YZ}{F} & w_1 &=& \dfrac{X^3 + YZ}{Y^2 + XZ}&\\
&&&&&&\\
v_2 &=& \dfrac{Z^2 - X^2Y}{F} & w_2 &=& \dfrac{X^4 - Y^3}{Y^2 + XZ}& = w_1X-Y.
\end{array}$$
This suggests that we should rewrite the equation for $\MM_0$ as
$$\MM_0: \;\;\; F^2  + XY(Y^2 + XZ) + X^2(X^3+YZ) - Z(Z^2 - X^2Y) = 0,$$ 
and that we can obtain $\Hat{\MM}$ by blowing up
$$F = Y^2+XZ = X^3+YZ = Z^2-X^2Y = 0.$$

\noindent \textbf{The locus} $\C$\\
\noindent Let $\S$ denote the surface
$$\S:\;\; F=Y^2+XZ=0.$$
Our $\Hat{M}$ coordinate patches $(\beta,v_1,v_2)$ and $(\gamma,w_1,w_2)$ cover
everything except the locus 
$$\C = \S \cap \MM_0.$$
The intersection of $\S$ with the 3-fold $\MM_0$ yields the new equation
$$X^5 + 2X^2YZ - Z^3 = 0.$$
(This was obtained by finding the Groebner basis for the ideal generated by $F$, $Y^2+XZ$, and
the equation for $\MM_0$.)  
\skp

\noindent For $X \neq 0$ we can write
$$Z = -\dfrac{Y^2}{X},$$
and so the equations for $\C$ become
$$\C: \;\;F = Y^2 + XZ = (X^4 - Y^3)^2 = 0, \skop (X \neq 0).$$
We can parametrize this curve by
\begin{eqnarray*}
X &=& t^3,\\
Y &=& t^4,\\
Z &=& -t^5,\\
F &=& 0.
\end{eqnarray*}
From this we see that blowing up $\C \subset \MM_0$ is equivalent to blowing up
$$F = Y^2+XZ = X^3+YZ = Z^2-X^2Y = X^4 - Y^3 = 0.$$
In this case we would have additional coordinates $v_3, w_3$ for the blowup:
$$\begin{array}{ccc|cccc}
\beta &=& \dfrac{Y^2+XZ}{F} & \gamma &=& \dfrac{F}{Y^2+XZ}&\\
&&&&&&\\
v_1 &=& \dfrac{X^3 + YZ}{F} & w_1 &=& \dfrac{X^3 + YZ}{Y^2 + XZ}&\\
&&&&&&\\
v_2 &=& \dfrac{Z^2 - X^2Y}{F} & w_2 &=& \dfrac{X^4 - Y^3}{Y^2 + XZ}&\\
&&&&&&\\
v_3 &=& \dfrac{X^4-Y^3}{F} & w_3 &=& \dfrac{Z^2 - X^2Y}{Y^2+XZ}.&
\end{array}$$
Note that both $v_3$ and $w_3$ add nothing new, as we can solve for them
in terms of the other coordinates:
$$\begin{array}{ccccc}
v_3 &=& \gamma^{-1}w_2 &=& \beta^4 v_2 - \beta^3 v_1^2 - v_1^3,\\
w_3 &=& \beta^{-1} v_2 &=& \gamma^4 w_2 + \gamma^2 w_1^3 + w_1^2.
\end{array}$$
These are precisely the additional coordinates we introduced in our resolution of the ideal sheaf (see Section 3.5).
Finally, with these identifications we find transition functions
$$\beta = \gamma^{-1}, \skop v_1 = \gamma^{-1} w_1, \skop v_2 = \gamma^3 w_2 + \gamma w_1^3 + 
\gamma^{-1} w_1^2,$$
which are exactly the ones we started with.

\section{Matrix factorization for a length 3 singularity}

Our singular geometry is in preferred versal form of type $E_8$, and the presence of the $T^3Z$ term in the polynomial\footnote{Note the change of coordinates!}
$$ X^2 - Y^3 + Z^5 + 3TYZ^2 + T^3Z = 0, $$
tells us that our geometry has a length 3 singularity which can be resolved by blowing up the open node in
the following $E_8$ diagram:\footnote{This uses results in the proof of the Katz--Morrison classification \cite{morrison}, as discussed in the Introduction.}

\setlength{\unitlength}{1 true in}
\begin{picture}(3,1)(1.9,.5)
\thicklines
\put(1.9,1){\circle*{.075}}
\put(1.9,1){\line(1,0){.5}}
\put(2.4,1){\circle*{.075}}
\put(2.4,1){\line(1,0){.4625}}
\put(2.9,1){\circle*{.075}}
\put(2.9,1){\line(0,-1){.4625}}
\put(2.9,.5){\circle{.075}}
\put(2.9375,1){\line(1,0){.4625}}
\put(3.4,1){\circle*{.075}}
\put(3.4,1){\line(1,0){.5}}
\put(3.9,1){\circle*{.075}}
\put(3.9,1){\line(1,0){.5}}
\put(4.4,1){\circle*{.075}}
\put(4.4,1){\line(1,0){.5}}
\put(4.9,1){\circle*{.075}}
\put(1.775,1.05){\makebox(.25,.25){\footnotesize 2}}
\put(2.275,1.05){\makebox(.25,.25){\footnotesize 4}}
\put(2.775,1.05){\makebox(.25,.25){\footnotesize 6}}
\put(2.925,.375){\makebox(.25,.25){\footnotesize 3}}
\put(3.275,1.05){\makebox(.25,.25){\footnotesize 5}}
\put(3.775,1.05){\makebox(.25,.25){\footnotesize 4}}
\put(4.275,1.05){\makebox(.25,.25){\footnotesize 3}}
\put(4.775,1.05){\makebox(.25,.25){\footnotesize 2}}
\end{picture}
\setlength{\unitlength}{1pt}

An alternative method for finding the small resolution for this length 3 singularity should thus be possible by deforming the $E_8$ matrix factorization corresponding to this node.  After some trial and error (and help from Maple), we find the following deformed matrix factorization for our singular equation:
$$\left(\begin{array}{cc} \psi & -xI \\ xI & \varphi \end{array}\right), \skop 
\left(\begin{array}{cc} \varphi & xI \\ -xI & \psi \end{array}\right),$$
where we have used Kn\"orrer's periodicity \cite[Chapter 12]{yoshino} and the fact that
$$\varphi = \left(\begin{array}{ccc} -y & z^2 & tz \\ t & -y & z^2 \\ z & t & -y \end{array}\right),
\skop \psi = \left(\begin{array}{ccc} y^2-tz^2 & yz^2+t^2z & z^4+tyz \\
 z^3+ty & y^2-tz^2 & yz^2+t^2z \\ yz+t^2 & z^3+ty & y^2-tz^2 \end{array}\right),$$
gives a matrix factorization for
$$-y^3 + z^5 + 3tyz^2 + t^3z = 0.$$

Using the approach in Chapter 5, this should allow us to perform the small resolution, and hence to extend
Theorem 1 to include a length 3 singularity.  We leave this for future work.
\end{onecolumn}

\begin{onecolumn}
\chapter{Other ADE cases}\label{ch:othercases}
There are two ADE cases left in the Intriligator--Wecht classification: $E_6$ and $E_8$.
Using the algorithm from Chapter~\ref{ch:algorithm}, we were able to find several global holomorphic
functions in each case, but not enough to perform the blow-down.  It is not clear whether or not
these cases will yield singular spaces of lengths 4, 5 or 6, and for the moment we are skeptical
that these small resolutions admit such simple descriptions.  Nevertheless, we record here our progress,
in the hope that it might be useful for future work.

\section{The Case $E_6$}

\noindent \textbf{The resolved geometry} $\Hat{\MM}$\\
\noindent From the Intriligator-Wecht superpotential
$$W(x,y) = \dfrac{1}{3}y^3 + \dfrac{1}{4}x^4,$$
we compute the resolved geometry $\Hat{\MM}$ in terms of transition functions
$$\beta = \gamma^{-1}, \skop v_1 = \gamma^{-1} w_1, \skop v_2 = \gamma^3 w_2 + \gamma^2 w_1^3 + 
\gamma^{-1} w_1^2.$$
To find the $\PP^1$s, we substitute $w_1(\gamma) = x + \gamma y$ into the $v_2$ transition function
$$v_2(\beta) = \beta^{-3} w_2 + \beta^{-2}(x + \beta^{-1}y)^3 + \beta(x + \beta^{-1}y)^2.$$
If we choose
$$w_2(\gamma) = \dfrac{x^3 - (x+\gamma y)^3}{\gamma},$$
in the $\beta$ chart the section is
$$v_1(\beta) = \beta x + y, \skop v_2(\beta) = \beta^{-2}x^3 + \beta^{-1}y^2+2xy+\beta x^2.$$
This is only holomorphic if
$$x^3 = y^2 = 0,$$
and so we have a single $\PP^1$ located at
$$w_1(\gamma) = w_2(\gamma) = 0, \skop v_1(\beta) = v_2(\beta) = 0.$$
This is exactly what we expect from computing critical points of the superpotential
$$\d W = x^3 dx + y^2 dy = 0.$$
\skp

\noindent \textbf{Global holomorphic functions}\\
\noindent The transition functions
$$\beta = \gamma^{-1}, \skop v_1 = \gamma^{-1} w_1, \skop v_2 = \gamma^3 w_2 + \gamma^2 w_1^3 +
\gamma^{-1} w_1^2,$$
are quasi-homogeneous if we assign the weights
$$\begin{array}{c|c|c|c|c|c}
\beta & v_1 & v_2 & \gamma & w_1 & w_2 \\
1 & 4 & 7 & -1 & 3 & 10
\end{array}.$$
The global holomorphic functions will thus necessarily be quasi-homogeneous in these weights.
We find (so far!) the following global holomorphic functions:
$$\begin{array}{c|ccc}
8 & y_1 &=& \beta v_2-v_1^2 \\
9 & y_2 &=& \beta y_1 \\
12 & y_3 &=& v_1 y_1 \\
15 & y_4 &=& v_2 y_1 \\
19 & y_5 &=& (v_1 v_2-\beta^3 y_1)y_1\\
22 & y_6 &=& (v_2^2 - \beta^2 v_1 y_1)y_1\\
29 & y_7 &=& (v_2^3 - 2\beta v_1^3 y_1 + \beta^5 y_1^2)y_1 
\end{array}$$
But there may be many more...
\skp

\noindent \textbf{The singular geometry} $\MM_0$\\
\noindent The functions $y_1,...,y_7$ satisfy 21 independent relations.  The first few are
$$\begin{array}{c|ccc}
24 & r_1 &=& y_1^3 - y_2y_4 + y_3^2 = 0 \\
27 & r_2 &=& y_2^3 - y_3y_4 + y_1y_5 = 0 \\
30 & r_3 &=& -y_4^2 + y_1 y_6 + y_3 y_2^2 = 0 \\
31 & r_4 &=& -y_6y_2+y_1^2y_4 + y_3y_5 = 0 \\
34 & r_5 &=& -y_3y_6+y_4y_5+y_1^2y_2^2=0\\
37 & r_6 &=& y_6y_4-y_1y_7-y_5y_2^2+2y_2y_3y_1^2=0\\
38 & r_7 &=& -y_2y_7+y_5^2-2y_1^2y_6+3y_1y_4^2=0\\
39 & r_8 &=& -y_4^2y_2 + y_4y_3^2+y_1y_6y_2-y_1y_3y_5 =0
\end{array}$$
and the remaining occur in degrees
$$41, 44, 45, 46, 52, 58, 62, 64, 66, 70, 94, 100, 138.$$ 
From a simple generating function we can generate all the monomials in each degree.
Here we list only up to degree 40, and only in cases where there is more than
one monomial in that degree.
$$\begin{array}{|c|c|c|c}
24 & y_1^3,\;y_3^2,\;y_2y_4 & 34 & y_1^2 y_2^2, \; y_4y_5, \; y_3y_6\\
27 & y_2^3,\; y_3y_4,\; y_1y_5 & 35 & y_2^3,\; y_1 y_3 y_4,\; y_1^2 y_5\\
28 & y_1^2 y_3, \; y_2 y_5 & 36 & y_2^4,\; y_3 y_1^3,\; y_3^3,\; y_2 y_3 y_4,\;y_1y_2y_5\\
29 & y_1y_2y_3,\; y_7 & 37 & y_2y_3y_1^2,\; y_5 y_2^2,\; y_4y_6,\; y_1 y_7\\
30 & y_3y_2^2,\; y_4^2,\; y_1y_6 & 38 & y_2^2 y_1 y_3, \; y_1 y_4^2,\; y_5^2,\;y_1^2 y_6,\; y_2y_7\\
31 & y_1^2y_4,\; y_3y_5,\; y_2y_6 & 39 & y_2^3y_3,\; y_4y_1^3,\; y_4y_3^2,\; y_2y_4^2,\;y_1y_3y_5,\;y_1y_2y_6\\
32 & y_1^4,\;y_1y_3^2,\;y_1y_2y_4 & 40 & y_1^5,\;y_1^2y_3^2,\;y_2y_1^2y_4,\;y_2y_3y_5,\;y_2^2y_6\\
33 & y_1^3y_2,\; y_3^2y_2,\; y_2^2y_4 &
\end{array}$$
Note that in almost every case, the relations $r_i$ involve every possible monomial in their respective
degrees.

\section{The Case $E_8$}

\noindent \textbf{The resolved geometry} $\Hat{\MM}$\\
\noindent From the Intriligator-Wecht superpotential
$$W(x,y) = \dfrac{1}{3}y^3 + \dfrac{1}{5}x^5,$$
we compute the resolved geometry $\Hat{\MM}$ in terms of transition functions
$$\beta = \gamma^{-1}, \skop v_1 = \gamma^{-1} w_1, \skop v_2 = \gamma^3 w_2 + \gamma^2 w_1^4 + 
\gamma^{-1} w_1^2.$$
To find the $\PP^1$s, we substitute $w_1(\gamma) = x + \gamma y$ into the $v_2$ transition function
$$v_2(\beta) = \beta^{-3} w_2 + \beta^{-2}(x + \beta^{-1}y)^4 + \beta(x + \beta^{-1}y)^2.$$
If we choose
$$w_2(\gamma) = \dfrac{x^4 - (x+\gamma y)^4}{\gamma},$$
in the $\beta$ chart the section is
$$v_1(\beta) = \beta x + y, \skop v_2(\beta) = \beta^{-2}x^4 + \beta^{-1}y^2+2xy+\beta x^2.$$
This is only holomorphic if
$$x^4 = y^2 = 0,$$
and so we have a single $\PP^1$ located at
$$w_1(\gamma) = w_2(\gamma) = 0, \skop v_1(\beta) = v_2(\beta) = 0.$$
This is exactly what we expect from computing critical points of the superpotential
$$\d W = x^4 dx + y^2 dy = 0.$$
\skp

\noindent \textbf{Global holomorphic functions}\\
\noindent The transition functions
$$\beta = \gamma^{-1}, \skop v_1 = \gamma^{-1} w_1, \skop v_2 = \gamma^3 w_2 + \gamma^2 w_1^4 +
\gamma^{-1} w_1^2,$$
are quasi-homogeneous if we assign the weights
$$\begin{array}{c|c|c|c|c|c}
\beta & v_1 & v_2 & \gamma & w_1 & w_2 \\
2 & 5 & 8 & -2 & 3 & 14
\end{array}.$$
The global holomorphic functions will thus necessarily be quasi-homogeneous in these weights.
We find (so far!) the following global holomorphic functions:
$$\begin{array}{c|ccc}
10 & y_1 &=& \beta v_2-v_1^2 \\
12 & y_2 &=& \beta y_1 \\
15 & y_3 &=& v_1 y_1 \\
18 & y_4 &=& v_2 y_1 \\
26 & y_5 &=& (v_2^2-\beta^3 y_1)y_1 
\end{array}$$
But there may be many more...
\skp

\noindent \textbf{The singular geometry} $\MM_0$\\
\noindent The functions $y_1,...,y_5$ satisfy 2 independent relations
$$\begin{array}{c|ccc}
30 & r_1 &=& y_1^3 - y_2y_4 + y_3^2 = 0 \\
36 & r_2 &=& y_2^3 - y_4^2 + y_1y_5 = 0 
\end{array}$$

\noindent \textbf{Remarks}
\begin{itemize}
\item Each relation uses every possible monomial in its degree.  
\item $r_1$ and $r_2$ are almost identical to the first two relations in the $E_6$ case.  Moreover,
the first two relations $r_1,r_2$ in the $E_6$ case are the only relations we would have if we only
considered the first five functions $y_1,..,y_5$.  
\item The functions $y_1,...,y_5$ are almost identical to the first five functions in the $E_6$ case!\skp
\end{itemize}

\end{onecolumn}

\begin{onecolumn}
\chapter{The `Hat' cases}\label{ch:hatcases}
In this section we analyze the extra `hat' cases in 
the Intriligator--Wecht classification of superpotentials.  In particular, we prove
{\bf Theorem 3}:\skp

\noindent {\em The singular geometries corresponding to the $\Hat{O},\Hat{A},\Hat{D}$ and $\Hat{E}$ cases in the Intriligator--Wecht classification of superpotentials are:
\begin{table}[h]
$$\begin{array}{c|c|c|c}
\mathrm{type} & W(x,y) & \partial_{w_1}E(\gamma,w_1) & \mathrm{singular\;\;geometry\;\;} \M_0\\
\hline
&&&\\
\Hat{O} & 0 & 0 & \CC^3/\ZZ_3\\
& & &\\
\Hat{A} & \12 y^2 & \gamma^2 w_1 & \CC \times \CC^2/\ZZ_2\\
& & &\\
\Hat{D} & xy^2 & \gamma w_1^2 & X^2 + Y^2Z - T^3 = 0\\
& & &\\
\Hat{E} & \dfrac{1}{3} y^3 & \gamma^2 w_1^2 & \mathrm{Spec} (\CC[a,b,u,v]/\ZZ_2)/(b^4 - u^2 - av).\\
&&&\\
\hline
\end{array}$$
\end{table}
}
\pagebreak 

\section{The Case $\Hat{O}$}

\noindent \textbf{The resolved geometry} $\Hat{\MM}$\\
\noindent From the Intriligator-Wecht superpotential
$$W(x,y) = 0,$$
we compute the resolved geometry $\Hat{\MM}$ in terms of transition functions
$$\beta = \gamma^{-1}, \skop v_1 = \gamma^{-1} w_1, \skop v_2 = \gamma^3 w_2.$$
To find the $\PP^1$s, we substitute $w_1(\gamma) = x + \gamma y$ into the $v_2$ transition function
$$v_2(\beta) = \beta^{-3} w_2.$$
If we choose
$$w_2(\gamma) = 0,$$
in the $\beta$ chart the section is
$$v_1(\beta) = \beta x + y, \skop v_2(\beta) = 0.$$
This is holomorphic for all $x$ and $y$,
and so we have a 2-parameter family of $\PP^1$s located at
$$w_1(\gamma) = x+\gamma y, w_2(\gamma) = 0, \skop v_1(\beta) = \beta x + y, v_2(\beta) = 0.$$
This is exactly what we expect from computing critical points of the superpotential
$$\d W = 0.$$
\skp

\noindent \textbf{Global holomorphic functions}\\
\noindent The transition functions
$$\beta = \gamma^{-1}, \skop v_1 = \gamma^{-1} w_1, \skop v_2 = \gamma^3 w_2$$
are quasi-homogeneous if we assign the weights
$$\begin{array}{c|c|c|c|c|c}
\beta & v_1 & v_2 & \gamma & w_1 & w_2 \\
1 & d+1 & e-3 & -1 & d & e
\end{array}.$$
Notice the freedom in choosing $d$ and $e$: there is a two-dimensional lattice of possible weight 
assignments. The global holomorphic functions will necessarily be quasi-homogeneous in these weights.
We find global holomorphic functions:
$$X_{ij} = \beta^{i} v_1^j v_2 = \gamma^{3-i-j} w_1^j w_2, \skop i,j \geq 0 ,\;\;\; i+j \leq 3.$$
\skp

\noindent \textbf{The singular geometry} $\MM_0$\\
\noindent If we rewrite our functions in a homogeneous manner as
$$\til{X}_{ij} = a^{3-i-j} b^i c^j, \skop i,j \geq 0 ,\;\;\; i+j \leq 3,$$
we can now identify the ring of global holomorphic functions as homogeneous
polynomials of degree 3 in 3 variables.  In other words, the ring is isomorphic to
$$\CC[a,b,c]^{\ZZ_3},$$
and our singular variety is simply
$$\MM_0: \;\;\; \CC^3/\ZZ_3.$$
\skp

\noindent \textbf{The blowup}\\
\noindent We now verify that we have identified the right singular space $\MM_0$ by inverting
the blow-down.  In the $\beta$ and $\gamma$ charts we find
$$\begin{array}{ccc|cccc}
\beta &=& X_{10}/X_{00} = \til{X}_{10}/\til{X}_{00} = b/a & \gamma &=& a/b &\\
&&&&&&\\
v_1 &=& X_{01}/X_{00} = \til{X}_{01}/\til{X}_{00} = c/a & w_1 &=& c/b&\\
&&&&&&\\
v_2 &=& X_{00} = \til{X}_{00} = a^3 & w_2 &=& b^3 &
\end{array}$$
which gives transition functions
$$\beta = \gamma^{-1}, \skop v_1 = \gamma^{-1} w_1, \skop v_2 = \gamma^3 w_2.$$
These are precisely what we started with!\skp

\noindent \textbf{Remark}\\
\noindent In the resolved $\Hat{\M}$ geometry, what we have is a $\PP^1$ inside a $\PP^2$ (or any other del Pezzo surface).  If you have a $\PP^2$ inside a Calabi-Yau and blow it down, you get $\CC^3/\ZZ^3$ as the singular point.

\section{The Case $\Hat{A}$}

\noindent \textbf{The resolved geometry} $\Hat{\MM}$\\
\noindent From the Intriligator-Wecht superpotential
$$W(x,y) = \dfrac{1}{2}x^2,$$
we compute the resolved geometry $\Hat{\MM}$ in terms of transition functions
$$\beta = \gamma^{-1}, \skop v_1 = \gamma^{-1} w_1, \skop v_2 = \gamma^3 w_2 + \gamma^2 w_1.$$
To find the $\PP^1$s, we substitute $w_1(\gamma) = x + \gamma y$ into the $v_2$ transition function
$$v_2(\beta) = \beta^{-3} (w_2 + y) + \beta^{-2}x.$$
If we choose
$$w_2(\gamma) = -y,$$
in the $\beta$ chart the section is
$$v_1(\beta) = \beta x + y, \skop v_2(\beta) = \beta^{-2}x.$$
This is only holomorphic $x=0$. Since $y$ is free,
we have a 1-parameter family of $\PP^1$s located at
$$w_1(\gamma) = \gamma y, w_2(\gamma) = -y, \skop v_1(\beta) = y, v_2(\beta) = 0.$$
This is exactly what we expect from computing critical points of the superpotential
$$\d W = x\d x = 0.$$
\skp

\noindent \textbf{Global holomorphic functions}\\
\noindent The transition functions
$$\beta = \gamma^{-1}, \skop v_1 = \gamma^{-1} w_1, \skop v_2 = \gamma^3 w_2 + \gamma^2 w_1$$
are quasi-homogeneous if we assign the weights
$$\begin{array}{c|c|c|c|c|c}
\beta & v_1 & v_2 & \gamma & w_1 & w_2 \\
1 & d+1 & d-2 & -1 & d & d+1
\end{array}.$$
Notice the freedom in choosing $d$: there is a one-dimensional lattice of possible weight 
assignments. The global holomorphic functions will necessarily be quasi-homogeneous in these weights.
We find global holomorphic functions:
$$\begin{array}{c|ccc}
d-2 & y_1 &=& v_2 = \gamma^3 w_2 + \gamma^2 w_1\\
d-1 & y_2 &=& \beta v_2 = \gamma^2 w_2 + \gamma w_1\\
d   & y_3 &=& \beta^2 v_2 = \gamma w_2 + w_1\\
d+1 & y_4 &=& \beta^3 v_2 - v_1 = w_2
\end{array}$$
\skp

\noindent \textbf{The singular geometry} $\MM_0$\\
\noindent The functions $y_i$ satisfy the single degree $2d-2$ relation
$$y_2^2 - y_1 y_3 = 0,$$
with $y_4$ free.  In other words, our singular geometry $\MM_0$ is a curve
of $A_1$ singularities, parametrized by $y_4$.
\skp

\noindent\textbf{The blowup}\\
\noindent We now verify that we have identified the right singular space $\MM_0$ by inverting
the blow-down.  In the $\beta$ and $\gamma$ charts we find
$$\begin{array}{ccc|cccc}
\beta &=& y_2/y_1 & \gamma &=& y_1/y_2 &\\
&&&&&&\\
v_1 &=& (y_2y_3-y_1y_4)/y_1 & w_1 &=& (y_2y_3-y_1y_4)/y_2&\\
&&&&&&\\
v_2 &=& y_1 & w_2 &=& y_4. &
\end{array}$$
This gives transition functions
$$\beta = \gamma^{-1}, \skop v_1 = \gamma^{-1} w_1, \skop v_2 = \gamma^3 w_2 + \gamma^2 w_1,$$
as expected.


\section{The Case $\Hat{D}$}

\noindent \textbf{The resolved geometry} $\Hat{\MM}$\\
\noindent From the Intriligator-Wecht superpotential
$$W(x,y) = x^2 y,$$
we compute the resolved geometry $\Hat{\MM}$ in terms of transition functions
$$\beta = \gamma^{-1}, \skop v_1 = \gamma^{-1} w_1, \skop v_2 = \gamma^3 w_2 + \gamma w_1^2.$$
To find the $\PP^1$s, we substitute $w_1(\gamma) = x + \gamma y$ into the $v_2$ transition function
$$v_2(\beta) = \beta^{-3}(w_2+y^2) + \beta^{-2}(2xy) + \beta^{-1}x^2.$$
If we choose
$$w_2(\gamma) = -y^2,$$
in the $\beta$ chart the section is
$$v_1(\beta) = \beta x + y, \skop v_2(\beta) = \beta^{-2}(2xy) + \beta^{-1}x^2.$$
This is only holomorphic if $x=0$. Since $y$ is free,
we have a 1-parameter family of $\PP^1$s located at
$$w_1(\gamma) = \gamma y, w_2(\gamma) = -y^2, \skop v_1(\beta) = y, v_2(\beta) = 0.$$
This is exactly what we expect from computing critical points of the superpotential
$$\d W = 2xy \;\d x + x^2\d y = 0.$$
\skp

\noindent \textbf{Global holomorphic functions}\\
\noindent The transition functions
$$\beta = \gamma^{-1}, \skop v_1 = \gamma^{-1} w_1, \skop v_2 = \gamma^3 w_2 + \gamma w_1^2$$
are quasi-homogeneous if we assign the weights
$$\begin{array}{c|c|c|c|c|c}
\beta & v_1 & v_2 & \gamma & w_1 & w_2 \\
1 & d+1 & 2d-1 & -1 & d & 2d+2
\end{array}.$$
Notice the freedom in choosing $d$: there is a one-dimensional lattice of possible weight 
assignments. The global holomorphic functions will necessarily be quasi-homogeneous in these weights.
We find global holomorphic functions:
$$\begin{array}{c|ccc}
2d-1 & y_1 &=& v_2 = \gamma^3 w_2 + \gamma w_1^2\\
2d & y_2 &=& \beta v_2 = \gamma^2 w_2 + w_1^2\\
3d   & y_3 &=& v_1 v_2 = \gamma^2 w_1 w_2 + w_1^3\\
2d+2 & y_4 &=& \beta^3 v_2 - v_1^2 = w_2
\end{array}$$
\skp

\noindent \textbf{The singular geometry} $\MM_0$\\
\noindent The functions $y_i$ satisfy the single degree $6d$ relation
$$y_3^2 - y_2^3 + y_1^2 y_4 =0.$$
\skp

\noindent \textbf{The blowup}\\
We now verify that we have identified the right singular space $\MM_0$ by inverting
the blow-down.  In the $\beta$ and $\gamma$ charts we find
$$\begin{array}{ccc|cccc}
\beta &=& y_2/y_1 & \gamma &=& y_1/y_2 &\\
&&&&&&\\
v_1 &=& y_3/y_1 & w_1 &=& y_3/y_2&\\
&&&&&&\\
v_2 &=& y_1 & w_2 &=& y_4, &
\end{array}$$
with transition functions
$$\beta = \gamma^{-1}, \skop v_1 = \gamma^{-1} w_1, \skop v_2 = \gamma^3 w_2 + \gamma w_1^2.$$

\section{The Case $\Hat{E}$}

\noindent \textbf{The resolved geometry} $\Hat{\MM}$\\
\noindent From the Intriligator-Wecht superpotential
$$W(x,y) = \dfrac{1}{3}x^3,$$
we compute the resolved geometry $\Hat{\MM}$ in terms of transition functions
$$\beta = \gamma^{-1}, \skop v_1 = \gamma^{-1} w_1, \skop v_2 = \gamma^3 w_2 + \gamma^2 w_1^2.$$
To find the $\PP^1$s, we substitute $w_1(\gamma) = x + \gamma y$ into the $v_2$ transition function
$$v_2(\beta) = \beta^{-3}(w_2+2xy+\gamma y^2) + \beta^{-2}x^2.$$
If we choose
$$w_2(\gamma) = -2xy-\gamma y^2,$$
in the $\beta$ chart the section is
$$v_1(\beta) = \beta x + y, \skop v_2(\beta) = \beta^{-2}x^2.$$
This is only holomorphic if $x=0$. Since $y$ is free,
we have a 1-parameter family of $\PP^1$s located at
$$w_1(\gamma) = \gamma y, w_2(\gamma) = -\gamma y^2, \skop v_1(\beta) = y, v_2(\beta) = 0.$$
This is exactly what we expect from computing critical points of the superpotential
$$\d W = x^2 \;\d x = 0.$$
\skp

\noindent \textbf{Global holomorphic functions}\\
\noindent The transition functions
$$\beta = \gamma^{-1}, \skop v_1 = \gamma^{-1} w_1, \skop v_2 = \gamma^3 w_2 + \gamma^2 w_1^2$$
are quasi-homogeneous if we assign the weights
$$\begin{array}{c|c|c|c|c|c}
\beta & v_1 & v_2 & \gamma & w_1 & w_2 \\
1 & d+1 & 2d-2 & -1 & d & 2d+1
\end{array}.$$
Notice the freedom in choosing $d$: there is a one-dimensional lattice of possible weight 
assignments. The global holomorphic functions will necessarily be quasi-homogeneous in these weights.
We find global holomorphic functions:
$$\begin{array}{c|ccc}
2d-2 & y_1 &=& v_2 \\
2d-1 & y_2 &=& \beta v_2 \\
2d   & y_3 &=& \beta^2 v_2 \\
3d-1 & y_4 &=& v_1 v_2\\
3d   & y_5 &=& \beta v_1 v_2\\
4d   & y_6 &=& v_1^2 v_2\\
4d+1 & y_7 &=& \beta(\beta^4 v_2-v_1^2)v_2 = \beta v_3 v_2\\
5d+1 & y_8 &=& v_1(\beta^4 v_2-v_1^2)v_2 = v_1 v_3 v_2\\
6d+2 & y_9 &=& (\beta^4 v_2-v_1^2)^2 v_2 = v_3^2 v_2
\end{array}$$
where we have defined
$$v_3 = \beta^4 v_2-v_1^2 = \gamma^{-1}w_2.$$
Note that this is the same definition made in our resolution of the ideal sheaf (Section 3.5). 
\skp

\noindent \textbf{The singular geometry} $\MM_0$\\
\noindent The functions $y_i$ satisfy a total of 20 distinct relations, most of which are obvious.
To simplify things, consider the monomial mapping
$$\beta^i v_1^j v_3^k v_2 \longmapsto a^{2-i-j-k} b^i c^j f^k.$$
Our functions now become
$$\begin{array}{c|ccc}
2d-2 & y_1 &=& a^2 \\
2d-1 & y_2 &=& ab \\
2d   & y_3 &=& b^2 \\
3d-1 & y_4 &=& ac\\
3d   & y_5 &=& bc\\
4d   & y_6 &=& c^2\\
4d+1 & y_7 &=& bf\\
5d+1 & y_8 &=& cf\\
6d+2 & y_9 &=& f^2
\end{array}$$
Note that
$$\begin{array}{ccccc}
\beta &=& y_2/y_1 &=& b/a\\
v_1 &=& y_4/y_1 &=& c/a\\
v_2 &=& y_1 &=& a^2\\
v_3 &=& y_7/y_2 &=& f/a,
\end{array}$$
so the relation defining $v_3$ becomes
$$af = b^4 - c^2.$$
This means we can add the function $af$ to our list, together with the relation:
$$y_{10} = af = b^4 - c^2.$$
\skp

\noindent Now the functions $y_1,...,y_{10}$ are exactly the 10 monomials of degree 2 in 4 variables, together
with the above relation.  The ring of global holomorphic functions is thus
$$(\CC[a,b,c,f]/\ZZ_2)/(af-b^4+c^2),$$
where the $\ZZ_2$ acts diagonally as -1.  In other words, we have a hypersurface in a $\ZZ_2$
quotient space:
$$(b^2+c)(b^2-c)=af \skop \mathrm{in} \;\;\;\CC^4/\ZZ_2.$$
We can immediately see from this equation that a small resolution, where we blow up an ideal of
the form
$$b^2+c=a=0,$$
won't work, since the $\ZZ_2$ action interchanges $b^2+c$ and $b^2-c$.  We will need to do
a big blow up of the origin instead.
\skp

\noindent \textbf{The blowup}\\
\noindent We now verify that we have identified the right singular space $\MM_0$ by inverting
the blow-down.  In the $\beta$ and $\gamma$ charts we find
$$\begin{array}{ccc|cccc}
\beta &=& b/a & \gamma &=& a/b &\\
&&&&&&\\
v_1 &=& c/a & w_1 &=& c/b &\\
&&&&&&\\
v_2 &=& a^2 & w_2 &=& f/b &\\
&&&&&&\\
v_3 &=& f/a & w_3 &=& b^2 &
\end{array}$$
We will perform the big blowup of the origin, with corresponding $\PP^3$ coordinates:
$$\begin{array}{ccccccccc}
a &=& b &=& c &=& f &=& 0.\\
\alpha && \delta && \rho && \nu &&
\end{array}$$
Note that all eight coordinates switch sign under the $\ZZ_2$ action.\skp
The blowup has four coordinate charts
$$\begin{array}{c|c|c|c}
\alpha = 1 & \delta = 1 & \rho = 1 & \nu = 1\\
&&&\\
a=a & a = \alpha_2 b & a = \alpha_3 c & a = \alpha_4 f\\
b=\delta_1 a & b=b & b=\delta_3 c & b = \delta_4 f\\
c = \rho_1 a & c=\rho_2 b & c=c & c = \rho_4 f\\
f = \nu_1 a & f = \nu_2 b & f = \nu_3 c & f=f\\
&&&\\
\nu_1 = \delta_1^4 a^2-\rho_1^2 & b^2 = \alpha_2\nu_2+\rho_2^2 & \alpha_3\nu_3 = \delta_3^4c^2-1 &
\alpha_4 = \delta_4^4 f^2 - \rho_4^2\\
&&&\\
(a^2,\delta_1,\rho_1) & (\alpha_2,\rho_2,\nu_2) & (\alpha_3,\delta_3,c^2,\nu_3) & (\delta_4,\rho_4,f^2)  
\end{array}$$
\skp

\noindent \textbf{Remarks}
\begin{itemize}
\item The functions $\alpha_i,\delta_i,\rho_i,$ and $\nu_i$ are all invariant under the $\ZZ_2$ action,
since they are all ratios of functions which change sign:
$$\delta_1 = \delta/\alpha, \;\;\rho_2 = \rho/\delta, \;\;...\; \mathrm{etc.}$$
\item Because $a,b,c,$ and $f$ all change sign under the $\ZZ_2$ action, we must take their invariant
counterparts $a^2, b^2, c^2,$ and $f^2$ when we list the final coordinates for each chart.
\item In the $\delta=1$ chart, we solve for $b^2$ instead of $b$, because $b$ is not an invariant function.
\item The blow up is nonsingular.  In the $\alpha = 1, \delta = 1,$ and $\nu=1$ charts we see this
because we are left with three coordinates and no relations, so these charts are all isomorphic to $\CC^3$.
In the $\rho = 1$ chart, we have a hypersurface in $\CC^4$ defined by the non-singular equation
$$\alpha_3 \nu_3 = \delta_3^4 c^2 -1.$$
\end{itemize}

\noindent \textbf{Transition functions}\\
\noindent Between the first two charts $\alpha = 1$ and $\delta = 1$, we have transition functions
\begin{eqnarray*}
\delta_1 &=& \delta/\alpha = \alpha_2^{-1}\\
\rho_1 &=& \rho/\alpha = (\delta/\alpha)(\rho/\delta) = \alpha_2^{-1} \rho_2\\
a^2 &=& \alpha_2^2 b^2 = \alpha_2^2 (\alpha_2 \nu_2 + \rho_2^2)\\
&=& \alpha_2^3\nu_2 + \alpha_2^2\rho_2^2
\end{eqnarray*}
Notice that
$$\begin{array}{ccccc|ccccc}
\delta_1 &=& b/a &=& \beta & \alpha_2 &=& a/b &=& \gamma\\
\rho_1 &=& c/a &=& v_1 & \rho_2 &=& c/b &=& w_1\\
a^2 &=& a^2 &=& v_2 & \nu_2 &=& f/b &=& w_2,
\end{array}$$
and so our transition functions are really
$$\beta = \gamma^{-1}, \skop v_1 = \gamma^{-1} w_1, \skop v_2 = \gamma^3 w_2 + \gamma^2 w_1^2.$$
These are exactly the $\Hat{E}$ transition functions we started with!

\section{Comparison with ADE cases}

We have seen that the singular geometries corresponding to Intriligator and Wecht's `hat' cases are given
by

$$\begin{array}{c|c|c}
\Hat{O} & W(x,y) = 0 & \CC^3/\ZZ_3\\
&&\\
\Hat{A} & W(x,y) = \dfrac{1}{2}y^2 & \CC[X,Y,Z,T]/(XY-Z^2) \cong \CC \times \CC^2/\ZZ_2\\
& & \mathrm{geometry\;has\;curve\;of\;} A_1 \mathrm{singularities}\\
&&\\
\Hat{D} & W(x,y) = xy^2 & y_1^2 - y_2^3 + y_3^2 y_4 = 0\\
&&\\
& & \mathrm{recall\;geometry\;for\;} D_{k+2}:\\
&& y_1^2 + y_2^3 + y_3^2y_4 + y_4^k y_2 = 0\\
&&\\
\Hat{E} & W(x,y) = \dfrac{1}{3}y^3 & (\CC[a,b,u,v]/\ZZ_2)/(b^4 - u^2 - av)\\
& & \mathrm{This\;is\;a\;hypersurface\;in\;} \CC^4/\ZZ_2.
\end{array}$$

Note that in both the $\Hat{A}$ and $\Hat{D}$ cases, the resulting equations can be obtained from the $A_k$ and
$D_{k+2}$ equations by dropping the $k$-dependent terms.  In other words, we are tempted to think of $\Hat{A}$ and $\Hat{D}$ as the $k \rightarrow \infty$ limit.  Perhaps in trying to come up with an analogous statement for $\Hat{E}$ we can learn something about the ``missing'' $E_6$ and $E_8$ cases.  In particular, it will be interesting to understand the role of these spaces in a geometric model for RG flow.
\end{onecolumn}

\begin{onecolumn}
\chapter{Future directions}\label{ch:future}
We have learned many things about the interplay between string dualities and Gorenstein threefold singularities.
The full story, however, is far from complete.  We end by posing a series of questions for the future which are beyond the scope of this work.

From the Intriligator--Wecht classification, we still need to understand the $E_6$ and $E_8$ cases.  
Using our algorithm we have found many global holomorphic functions, but not enough to 
give us the blow-down.  There are also questions which arise from the extra `hat' cases. 
What is the interpretation of the $\Hat{O},\Hat{A},\Hat{D},$ and $\Hat{E}$ geometries from the string
theory perspective?  The $\PP^1$s are no longer isolated; do they correspond to D-branes wrapping families of $\PP^1$s?  Moreover, what is the role of higher order terms in the superpotential?  

Do we have a geometric model for RG flow?  Proposition 1 suggests that the geometry might encode something about the RG fixed points of the corresponding matrix models or gauge theories.  Can Proposition 1 be extended?  The discussion
in Section 7.2 suggests that this should be possible, since we can find an appropriate change of coordinates to make
it work in the $A_k$ cases.  Finding more general coordinate changes which can show how the rest of the terms
in the superpotential are affected when a bundle-changing coordinate is ``integrated out'' is a necessary step in
developing this kind of geometric picture.

Furthermore, Intriligator and Wecht have a chart of all possible flows between the RG fixed points.  We can make a similar chart based on our geometric framework.  Do they match?
Finally, what is the role of fundamentals?   Our entire analysis involves only adjoint fields, which correspond geometrically to parameters of the $\PP^1$ deformation space.  Intriligator and Wecht only find the ADE classification for superpotentials involving 2 adjoint fields, but their paper also analyzes many cases with fundamentals.  Is it possible to have a geometric interpretation for these fields?

As far as Ferrari's construction is concerned, there are many open ends to be explored.  Can we generalize Ferrari's framework to include perturbation terms for both $v_1$ and $v_2$ transition functions?  Can we generalize for cases where the geometry is specified by more than two charts?  This would enable more flexibility in 
identifying superpotentials in a ``bottom-up'' approach.  On the other hand, the techniques developed in \cite{aspinwall} in principle allow computation of the superpotential in general.  In cases where the superpotential cannot be easily identified in the transition functions, perhaps this approach should be used instead.

Moreover, in all of our new cases there is still
work to be done to complete the remaining steps in Ferrari's program.  For example, what is the solution to the matrix model corresponding to the length 3 singularity?  And what can we learn about the matrix models corresponding to the `hat' cases?  Although the singularities are no longer isolated, is it still possible to compute resolvents from the geometry?  If Ferrari's conjecture about the Calabi-Yau geometry encoding the solution to the matrix model is correct, we should now be able to solve the matrix models corresponding to the length 3 and `hat' cases.  If solutions are already known (or can be computed using traditional matrix model techniques), these examples will provide new tests to the conjecture.

Finally, there is the idea of Conjecture 1 -- that small resolutions for Gorenstein threefold singularities
might be obtainable by deforming matrix factorizations in the surface case.  In Theorem 1 we were able to prove this in the length 1 and length 2 cases, and with the help of the Intriligator--Wecht $E_7$ case we are very close to showing this in length 3 as well.  An ambitious, but tangible goal is to extend this result to lengths 4, 5 and 6.  In particular, this may lead us to an alternative (and more concrete) proof of the original Katz--Morrison classification \cite{morrison}.
\end{onecolumn}


\appendix


\chapter{Maple code for blow-down algorithm}\label{ch:maple}
\section{Description of the code}

\noindent The actual maple routine has the following parts:

\begin{figure}[t]
\begin{picture}(4000,2500)(0,0)
\putsquare<1`1`1`1;2000`1100>(1000,1200)[\framebox{Input}`\framebox{\shortstack{Ideal \& \\ Groebner basis}}`\framebox{{\tt Polysearch}}`\framebox{{\tt xreduce}};\mbox{perturbation \& weights}`\mbox{L,X}`\mbox{\shortstack{reduction\\ functions}}`\mbox{L[j], array F}]
\putsquare<0`1`0`-1;2000`1000>(1000,200)[\phantom{\framebox{{\tt Polysearch}}}`\phantom{\framebox{{\tt xreduce}}}`\framebox{{\tt newfun}}`\framebox{\shortstack{Procedures \\ on Arrays}}; `\mbox{\shortstack{f,X\\ if mix(f)=0}}``\mbox{array routines}]
\putmorphism(1000,1240)(1,0)[\phantom{\framebox{{\tt Polysearch}}}`\phantom{\framebox{{\tt xreduce}}}`\mbox{f = L[j] + new terms}]{2000}{-1}a
\putmorphism(1040,1200)(0,1)[\phantom{\framebox{{\tt newfun}}} `\phantom{\framebox{\shortstack{Procedures \\ on Arrays}}}`\mbox{updated X}]
{1000}{-1}r
\end{picture}
\caption{An algorithm for constructing the blow down}
\label{fig:algorithm}
\end{figure}
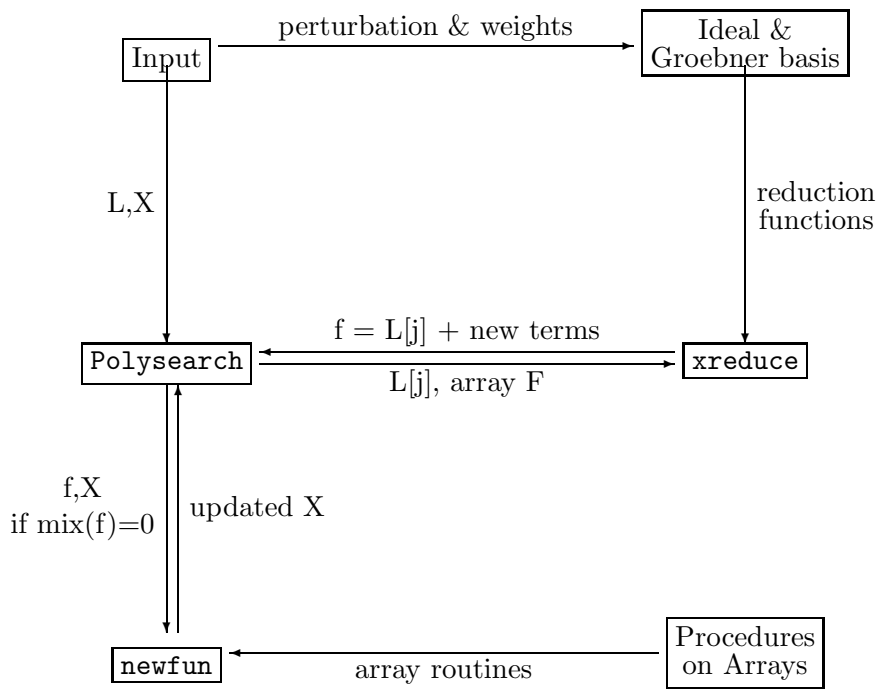

\begin{itemize}
\item Input information
      \begin{itemize}
      \item perturbation terms for ADE cases
      \item weights for $\beta,v_1,v_2$ in each case
      \item initial (empty) seed array $X$.
      \end{itemize}
\item Ideal-related definitions
      \begin{itemize}
      \item ideal in $\CC[\beta,v_1,v_2,\gamma,w_1,w_2]$
      \item term order T (lexdeg, eliminate $v_2$)
      \item term order TP (wdeg, degrees reflect poles)
      \item groebner basis G for ideal, wrt T
      \item {\tt modG, arraymodG, shift} and {\tt mix}, all compute reductions modulo the ideal G
      \end{itemize} 
\item Lists of monomials
      \begin{itemize}
      \item {\tt `genfun'} creates generating function of q-homog monomials
      \item {\tt `List'} extracts all monomials in weighted degree n (uses {\tt `genfun'})
      \end{itemize}
\item Procedures on arrays (used exclusively for {\tt `newfun'} and {\tt `invert'})
      \begin{itemize}
      \item {\tt `concat'} takes an array and adds an element to the end
      \item {\tt `yarray'} creates an array of indeterminates $[y_1,...,y_n]$
      \item {\tt `equal0'} converts an array A into a list $\{A_i=0\}$
      \item {\tt `subtract'} takes the difference between two arrays
      \item {\tt `parcials'} takes the derivative of all in functions in an array A wrt a given variable,
             then evaluates at a point `P' (a list of equations) and returns a list
      \item {\tt `zero'} creates a list of m elements which are all 0
      \item {\tt `extract'} picks out the elements which are pure in $\gamma,w_1,w_2$ from an array G.
      \end{itemize}
\item Main routines
      \begin{itemize}
      \item {\tt `xreduce'} is *the* main routine: it inputs a function and a list {\tt F} of other functions
      {\tt F[i,1]} together with their leading mixed terms { \tt F[i,2]} and leading mixed coefficients {\tt F[i,3]}.
      The function is reduced as much as possible by {\tt F}, and what is returned is a polynomial in which
      the original function is modified by terms {\tt F[i,1]} which are subtracted in order to reduce the leading
      mixed term.  (Uses {\tt shift, mix,} and {\tt TP}.)
      \item {\tt `polysearch'} receives an array `L' of polynomials and a list $X$ of global holomorphic functions
      as input.  It uses {\tt `xreduce'} to cycle through each polynomial and maximally reduce its leading mixed term.  
      It then uses the results to build the
      array {\tt F[i,j]} of maximally reduced polynomials and their mixed terms, all in the same weighted degree.     When the mixed term is 0, {\tt `newfun'} is used to check if the global homogeneous function is in the ring 
      generated by the list $Y$ of prior ghf's.  The final output is the updated list $Y$.
      (Uses {\tt shift, mix, TP, `xreduce'}, and {\tt `newfun'}.)
      \item {\tt `newfun'} checks if the global holomorphic function `f' is new, in the sense that it is not in the 
      ring generated by previous ghf's.  It returns the (potentially updated) list of truly independent ghf's.
      (Uses {\tt `yarray', `concat', `subtract', `extract', `parcials', `equal0',`zero'}.)
      \end{itemize}
\item Execution commands
      \begin{itemize}
      \item {\tt `findpoly'} searches for global holomorphic functions within certain bounds on the weighted degree.
      (Uses {\tt `List'} and {\tt `polysearch'})
      \item {\tt `invert'} solves for original variables ($\beta,v_1,v_2$ or $\gamma,w_1,w_2$) in terms of the new
      ones: global holomorphic functions fed in through the array A. (Uses {\tt `yarray'} and {\tt `subtract'}.)
      \end{itemize}
\end{itemize}

\section{Maple code}

\bigskip 

\begin{flushleft}
\begin{scriptsize}
\def\emptyline{\vspace{12pt}}
\DefineParaStyle{maplegroup}
\DefineParaStyle{Author}
\DefineParaStyle{subtitle}
\DefineParaStyle{Bullet Item}
\DefineParaStyle{Dash Item}
\DefineParaStyle{Diagnostic}
\DefineParaStyle{Error}
\DefineParaStyle{Heading 1}
\DefineParaStyle{Heading 2}
\DefineParaStyle{Heading 3}
\DefineParaStyle{Heading 4}
\DefineParaStyle{Normal}
\DefineParaStyle{Text Output}
\DefineParaStyle{Title}
\DefineParaStyle{Warning}
\DefineCharStyle{2D Math}
\DefineCharStyle{Help Heading}
\DefineCharStyle{Maple Input}
\DefineCharStyle{FixedWidth}
\DefineCharStyle{endFixedWidth}
\begin{mapleinput}
\mapleinline{active}{1d}{with(Groebner):}{}
\end{mapleinput}\begin{mapleinput}
\mapleinline{active}{1d}{extraterm:=0:}{}
\end{mapleinput}

\medskip 

THE FOLLOWING MUST BE CHANGED ACCORDING TO ADE  CASE\\

\medskip 

These are the perturbation terms for each of the ADE cases.
\begin{mapleinput}
\mapleinline{active}{1d}{Ohat:= 0: Ahat:= g^2*w[1]: Dhat:= g*w[1]^2: Ehat:= g^2*w[1]^2:}{}
\end{mapleinput}\begin{mapleinput}
\mapleinline{active}{1d}{Ak:= k -> g^2*w[1]^k + w[1]: Dk:= k -> g^2*w[1]^k + w[1]^2:}{}
\end{mapleinput}\begin{mapleinput}
\mapleinline{active}{1d}{E6:= b*w[1]^2 + g^2*w[1]^3: E7:= b*w[1]^2 + g*w[1]^3: E8:= b*w[1]^2 + g^2*w[1]^4:}{}
\end{mapleinput}These are the weights [deg(b),deg(v[1]), deg(v[2])] in each of the cases.\begin{mapleinput}
\mapleinline{active}{1d}{deg_Ohat:= [1,1,1]: deg_Ahat:= [1,4,1]: deg_Dhat:=[1,2,1]: deg_Ehat:=[1,3,2]:}{}
\end{mapleinput}\begin{mapleinput}
\mapleinline{active}{1d}{deg_Ak:= k -> [k-1,k+1,2]: deg_Dk:= k -> [k-2,k,4]:}{}
\end{mapleinput}\begin{mapleinput}
\mapleinline{active}{1d}{deg_E6:= [1,4,7]: deg_E7:= [1,3,5]: deg_E8:= [2,5,8]:}{}
\end{mapleinput}This is the perturbation term in the transition functions, and what distinguishes our different cases.\begin{mapleinput}
\mapleinline{active}{1d}{pterm:=Ak(4);}{}
\end{mapleinput}"d" captures the weights for b, v[1], v[2] in order to make the transition functions quasi-homogeneous. \begin{mapleinput}
\mapleinline{active}{1d}{d:=deg_Ak(4);}{}
\end{mapleinput}The first holomorphic function goes into the "seed" array X.\begin{mapleinput}
\mapleinline{active}{1d}{X:=Array([]):}{}
\end{mapleinput}

\medskip 

IDEAL AND REDUCTION MODULO IDEAL\\

\medskip 

We define the ideal coming from the transition functions, in which we will look for global holomorphic functions.\begin{mapleinput}
\mapleinline{active}{1d}{ideal:= [b*g-1, v[1]-b*w[1],v[2]-g^3*w[2]-pterm];}{}
\end{mapleinput}Our chosen term order for computing groebner basis and reduction mod the ideal.\begin{mapleinput}
\mapleinline{active}{1d}{T:=lexdeg([v[2]],[w[2],v[1],w[1],b,g]):}{}
\end{mapleinput}The term order for keeping track of the polar degree, which we use to order the cancellations. \begin{mapleinput}
\mapleinline{active}{1d}{TP:=wdeg([1,1,1,-1,0,0],[b,v[1],v[2],g,w[1],w[2]]):}{}
\end{mapleinput}We compute the groebner basis of `ideal' with respect to the term order T.\begin{mapleinput}
\mapleinline{active}{1d}{G:=gbasis(ideal,T);}{}
\end{mapleinput}The function "modG" computes the reduction of a function f modulo the ideal G.\begin{mapleinput}
\mapleinline{active}{1d}{modG:= f -> normalf(f,G,T):}{}
\end{mapleinput}The procedure "arraymodG" computes the reduction of every function in the array "A" modulo the ideal G.\begin{mapleinput}
\mapleinline{active}{1d}{arraymodG:=proc(A::Array)}{}
\end{mapleinput}\begin{mapleinput}
\mapleinline{active}{1d}{   local n,i,B;}{}
\end{mapleinput}\begin{mapleinput}
\mapleinline{active}{1d}{   n:=ArrayNumElems(A);}{}
\end{mapleinput}\begin{mapleinput}
\mapleinline{active}{1d}{   B:=Array(1..n);}{}
\end{mapleinput}\begin{mapleinput}
\mapleinline{active}{1d}{   for i to n do}{}
\end{mapleinput}\begin{mapleinput}
\mapleinline{active}{1d}{\hspace{.25in}B[i]:=modG(A[i]);}{}
\end{mapleinput}\begin{mapleinput}
\mapleinline{active}{1d}{   end do;}{}
\end{mapleinput}\begin{mapleinput}
\mapleinline{active}{1d}{   return B;}{}
\end{mapleinput}\begin{mapleinput}
\mapleinline{active}{1d}{end:}{}
\end{mapleinput}The function "mix" extracts the mixed terms from a polynomial, by subtracting off the pure terms from each chart.\begin{mapleinput}
\mapleinline{active}{1d}{mix:= f -> modG(f) - subs(g=0,w[1]=0,w[2]=0,modG(f)) - subs(b=0,v[1]=0,v[2]=0,modG(f)):}{}
\end{mapleinput}The function "shift" just brings terms in b,v[1],v[2] over to the left-hand-side.\begin{mapleinput}
\mapleinline{active}{1d}{shift:= f -> f - subs(g=0,w[1]=0,w[2]=0,modG(f)):}{}
\end{mapleinput}

\medskip 

GENERATING MONOMIALS OF SAME WEIGHTED HOMOGENEOUS DEGREE\\

\medskip 

The function "genfun" creates a generating function for the quasi-homogeneous monomials in each degree up to n.\begin{mapleinput}
\mapleinline{active}{1d}{genfun:= n -> convert(map(expand,series(1/((1-b*u^d[1])*(1-v[1]*u^d[2])*(1-v[2]*u^d[3])),u,n+1)),polynom):}{}
\end{mapleinput}The function "List" uses the above generating function to extract the monomials in weighted degree n.\begin{mapleinput}
\mapleinline{active}{1d}{List:= n -> Array([op(1+coeff(genfun(n),u,n))]):}{}
\end{mapleinput}

\medskip 

PROCEDURES ON ARRAYS

\medskip 

The "concat" procedure takes an array "A" and tacks on an extra element "a" onto the end.
\begin{mapleinput}
\mapleinline{active}{1d}{concat:= proc(A::Array, a)}{}
\end{mapleinput}\begin{mapleinput}
\mapleinline{active}{1d}{   local B,n;}{}
\end{mapleinput}\begin{mapleinput}
\mapleinline{active}{1d}{   n:= ArrayNumElems(A);}{}
\end{mapleinput}\begin{mapleinput}
\mapleinline{active}{1d}{   B:= Array(1...n+1);}{}
\end{mapleinput}\begin{mapleinput}
\mapleinline{active}{1d}{   B[1..n]:= A[1..n];}{}
\end{mapleinput}\begin{mapleinput}
\mapleinline{active}{1d}{   B[n+1]:= a;}{}
\end{mapleinput}\begin{mapleinput}
\mapleinline{active}{1d}{   return B;}{}
\end{mapleinput}\begin{mapleinput}
\mapleinline{active}{1d}{end:}{}
\end{mapleinput}

The "yarray" procedure creates an array with indeterminates y[1],...,y[n].
\begin{mapleinput}
\mapleinline{active}{1d}{yarray:= proc(n::integer)}{}
\end{mapleinput}\begin{mapleinput}
\mapleinline{active}{1d}{   local Y,i;}{}
\end{mapleinput}\begin{mapleinput}
\mapleinline{active}{1d}{   Y:= Array(1...n);}{}
\end{mapleinput}\begin{mapleinput}
\mapleinline{active}{1d}{   for i from 1 to n do}{}
\end{mapleinput}\begin{mapleinput}
\mapleinline{active}{1d}{\hspace{.25in}Y[i]:= y[i];}{}
\end{mapleinput}\begin{mapleinput}
\mapleinline{active}{1d}{   end do;}{}
\end{mapleinput}\begin{mapleinput}
\mapleinline{active}{1d}{   return Y;}{}
\end{mapleinput}\begin{mapleinput}
\mapleinline{active}{1d}{end:}{}
\end{mapleinput}

The procedure "equal0" converts an array into a list  {A[i]=0}.
\begin{mapleinput}
\mapleinline{active}{1d}{equal0:= proc(A::Array)}{}
\end{mapleinput}\begin{mapleinput}
\mapleinline{active}{1d}{   local B,i,n;}{}
\end{mapleinput}\begin{mapleinput}
\mapleinline{active}{1d}{   n:=ArrayNumElems(A);}{}
\end{mapleinput}\begin{mapleinput}
\mapleinline{active}{1d}{   B:=Array(1...n);}{}
\end{mapleinput}\begin{mapleinput}
\mapleinline{active}{1d}{   for i from 1 to n do}{}
\end{mapleinput}\begin{mapleinput}
\mapleinline{active}{1d}{\hspace{.25in}B[i]:= A[i]=0;}{}
\end{mapleinput}\begin{mapleinput}
\mapleinline{active}{1d}{   end do;}{}
\end{mapleinput}\begin{mapleinput}
\mapleinline{active}{1d}{   return convert(B,'list');}{}
\end{mapleinput}\begin{mapleinput}
\mapleinline{active}{1d}{end:}{}
\end{mapleinput}

The "subtract" procedure takes the difference between two arrays.
\begin{mapleinput}
\mapleinline{active}{1d}{subtract:= proc(A::Array, B::Array)}{}
\end{mapleinput}\begin{mapleinput}
\mapleinline{active}{1d}{   local C,i,n;}{}
\end{mapleinput}\begin{mapleinput}
\mapleinline{active}{1d}{   n:= ArrayNumElems(A);}{}
\end{mapleinput}\begin{mapleinput}
\mapleinline{active}{1d}{   C:= Array(1...n);}{}
\end{mapleinput}\begin{mapleinput}
\mapleinline{active}{1d}{   for i from 1 to n do}{}
\end{mapleinput}\begin{mapleinput}
\mapleinline{active}{1d}{\hspace{.25in}C[i]:= A[i]-B[i]}{}
\end{mapleinput}\begin{mapleinput}
\mapleinline{active}{1d}{   end do;}{}
\end{mapleinput}\begin{mapleinput}
\mapleinline{active}{1d}{   return C;}{}
\end{mapleinput}\begin{mapleinput}
\mapleinline{active}{1d}{end:}{}
\end{mapleinput}

The procedure "parcials" takes the derivative of all of the functions in A w.r.t. the variable "var", then evaluates at "P".
\begin{mapleinput}
\mapleinline{active}{1d}{parcials:= proc(A::Array,P::list,var)}{}
\end{mapleinput}\begin{mapleinput}
\mapleinline{active}{1d}{   local i,n,t;}{}
\end{mapleinput}\begin{mapleinput}
\mapleinline{active}{1d}{   n:=ArrayNumElems(A); t:=Array(1..n);}{}
\end{mapleinput}\begin{mapleinput}
\mapleinline{active}{1d}{   for i from 1 to n do}{}
\end{mapleinput}\begin{mapleinput}
\mapleinline{active}{1d}{\hspace{.25in}t[i]:=eval(subs(P,diff(A[i],var)));}{}
\end{mapleinput}\begin{mapleinput}
\mapleinline{active}{1d}{   end do;}{}
\end{mapleinput}\begin{mapleinput}
\mapleinline{active}{1d}{   return convert(t,'list'); }{}
\end{mapleinput}\begin{mapleinput}
\mapleinline{active}{1d}{end:}{}
\end{mapleinput}

The procedure "zero" creates an array of 0's.
\begin{mapleinput}
\mapleinline{active}{1d}{zero:=proc(m::integer)}{}
\end{mapleinput}\begin{mapleinput}
\mapleinline{active}{1d}{   local A,i; A:=Array(1..m);}{}
\end{mapleinput}\begin{mapleinput}
\mapleinline{active}{1d}{   for i to m do}{}
\end{mapleinput}\begin{mapleinput}
\mapleinline{active}{1d}{\hspace{.25in}A[i]:=0;}{}
\end{mapleinput}\begin{mapleinput}
\mapleinline{active}{1d}{   end do;}{}
\end{mapleinput}\begin{mapleinput}
\mapleinline{active}{1d}{   return convert(A,'list');}{}
\end{mapleinput}\begin{mapleinput}
\mapleinline{active}{1d}{end:}{}
\end{mapleinput}

The procedure "extract" picks out elements which are pure in (g,w[1],w[2]) from the array G.
\begin{mapleinput}
\mapleinline{active}{1d}{extract:=proc(G::Array)}{}
\end{mapleinput}\begin{mapleinput}
\mapleinline{active}{1d}{   local n,i,j,h,H;}{}
\end{mapleinput}\begin{mapleinput}
\mapleinline{active}{1d}{   n:=ArrayNumElems(G); H:=Array(1..n); j:=1;}{}
\end{mapleinput}\begin{mapleinput}
\mapleinline{active}{1d}{   for i to n do}{}
\end{mapleinput}\begin{mapleinput}
\mapleinline{active}{1d}{\hspace{.25in}h:= subs(v[2]=0,v[1]=0,b=0,G[i]);    }{}
\end{mapleinput}\begin{mapleinput}
\mapleinline{active}{1d}{\hspace{.25in}if G[i]=h then}{}
\end{mapleinput}\begin{mapleinput}
\mapleinline{active}{1d}{\hspace{.5in}H[j]:=G[i];}{}
\end{mapleinput}\begin{mapleinput}
\mapleinline{active}{1d}{\hspace{.5in}j:=j+1;}{}
\end{mapleinput}\begin{mapleinput}
\mapleinline{active}{1d}{\hspace{.25in}end if;}{}
\end{mapleinput}\begin{mapleinput}
\mapleinline{active}{1d}{   end do;}{}
\end{mapleinput}\begin{mapleinput}
\mapleinline{active}{1d}{   return H[1...j-1];}{}
\end{mapleinput}\begin{mapleinput}
\mapleinline{active}{1d}{end:}{}
\end{mapleinput}

\medskip 

THE MAIN REDUCTION ROUTINES

\medskip 

The procedure "xreduce" is the main routine.  It reduces the function mod the array "F", returning a modified function "f" whose mixed terms have lower degree.
\begin{mapleinput}
\mapleinline{active}{1d}{xreduce:= proc (function,F::Array)}{}
\end{mapleinput}\begin{mapleinput}
\mapleinline{active}{1d}{   local i,n,f,h,lt,lc,subtraction;}{}
\end{mapleinput}\begin{mapleinput}
\mapleinline{active}{1d}{   n:= ArrayNumElems(F)/3;}{}
\end{mapleinput}\begin{mapleinput}
\mapleinline{active}{1d}{   f:= shift(function); }{}
\end{mapleinput}\begin{mapleinput}
\mapleinline{active}{1d}{   h:= mix(f);}{}
\end{mapleinput}\begin{mapleinput}
\mapleinline{active}{1d}{   while h <> 0 do}{}
\end{mapleinput}\begin{mapleinput}
\mapleinline{active}{1d}{\hspace{.25in}lt:= leadterm(h,TP); }{}
\end{mapleinput}\begin{mapleinput}
\mapleinline{active}{1d}{\hspace{.25in}lc:= leadcoeff(h,TP);}{}
\end{mapleinput}\begin{mapleinput}
\mapleinline{active}{1d}{\hspace{.25in}subtraction:=0;}{}
\end{mapleinput}\begin{mapleinput}
\mapleinline{active}{1d}{\hspace{.25in}for i from n to 1 by -1 do}{}
\end{mapleinput}\begin{mapleinput}
\mapleinline{active}{1d}{\hspace{.5in}if F[i,2] = lt then}{}
\end{mapleinput}\begin{mapleinput}
\mapleinline{active}{1d}{\hspace{.75in}f:= f-lc/F[i,3]*F[i,1];}{}
\end{mapleinput}\begin{mapleinput}
\mapleinline{active}{1d}{\hspace{.75in}h:= mix(f);}{}
\end{mapleinput}\begin{mapleinput}
\mapleinline{active}{1d}{\hspace{.75in}subtraction:= 1;}{}
\end{mapleinput}\begin{mapleinput}
\mapleinline{active}{1d}{\hspace{.5in}end if;}{}
\end{mapleinput}\begin{mapleinput}
\mapleinline{active}{1d}{\hspace{.5in}if subtraction = 1 then }{}
\end{mapleinput}\begin{mapleinput}
\mapleinline{active}{1d}{\hspace{.75in}break; }{}
\end{mapleinput}\begin{mapleinput}
\mapleinline{active}{1d}{\hspace{.5in}end if;}{}
\end{mapleinput}\begin{mapleinput}
\mapleinline{active}{1d}{\hspace{.25in}end do;}{}
\end{mapleinput}\begin{mapleinput}
\mapleinline{active}{1d}{\hspace{.25in}if subtraction = 0 then}{}
\end{mapleinput}\begin{mapleinput}
\mapleinline{active}{1d}{\hspace{.5in}h:= h - lc*lt;}{}
\end{mapleinput}\begin{mapleinput}
\mapleinline{active}{1d}{\hspace{.25in}end if;}{}
\end{mapleinput}\begin{mapleinput}
\mapleinline{active}{1d}{   end do;}{}
\end{mapleinput}\begin{mapleinput}
\mapleinline{active}{1d}{   return f;}{}
\end{mapleinput}\begin{mapleinput}
\mapleinline{active}{1d}{end:}{}
\end{mapleinput}

The procedure "polysearch" receives an array "L" of polynomials as input, uses "xreduce", and finds the global holomorphic functions.
\begin{mapleinput}
\mapleinline{active}{1d}{polysearch:= proc(L::Array, X::Array, dw::integer)}{}
\end{mapleinput}\begin{mapleinput}
\mapleinline{active}{1d}{   local i,j,n,f,F,sl,nf,Y,t;}{}
\end{mapleinput}\begin{mapleinput}
\mapleinline{active}{1d}{   n:= ArrayNumElems(L); }{}
\end{mapleinput}\begin{mapleinput}
\mapleinline{active}{1d}{   F:= Array(1..n,1..3); }{}
\end{mapleinput}\begin{mapleinput}
\mapleinline{active}{1d}{   F[1,1]:= r(L[1]); }{}
\end{mapleinput}\begin{mapleinput}
\mapleinline{active}{1d}{   sl:= s(F[1,1]);}{}
\end{mapleinput}\begin{mapleinput}
\mapleinline{active}{1d}{   F[1,2]:= leadterm(sl,TP); }{}
\end{mapleinput}\begin{mapleinput}
\mapleinline{active}{1d}{   F[1,3]:= leadcoeff(sl,TP);}{}
\end{mapleinput}\begin{mapleinput}
\mapleinline{active}{1d}{   Y:= X;}{}
\end{mapleinput}\begin{mapleinput}
\mapleinline{active}{1d}{   print('n'=n); }{}
\end{mapleinput}\begin{mapleinput}
\mapleinline{active}{1d}{   j:= 2;}{}
\end{mapleinput}\begin{mapleinput}
\mapleinline{active}{1d}{   for i from 2 to n do}{}
\end{mapleinput}\begin{mapleinput}
\mapleinline{active}{1d}{\hspace{.25in}print(i,j); print('wdeg'=dw);}{}
\end{mapleinput}\begin{mapleinput}
\mapleinline{active}{1d}{\hspace{.25in}f:= xreduce(L[i],F[1..j-1,1..3]);}{}
\end{mapleinput}\begin{mapleinput}
\mapleinline{active}{1d}{\hspace{.25in}if f <> 0 then}{}
\end{mapleinput}\begin{mapleinput}
\mapleinline{active}{1d}{\hspace{.5in}sl:= mix(f);}{}
\end{mapleinput}\begin{mapleinput}
\mapleinline{active}{1d}{\hspace{.5in}F[j,1]:= f;}{}
\end{mapleinput}\begin{mapleinput}
\mapleinline{active}{1d}{\hspace{.5in}if sl = 0 then}{}
\end{mapleinput}\begin{mapleinput}
\mapleinline{active}{1d}{\hspace{.75in}print('f'=F[j,1]);}{}
\end{mapleinput}\begin{mapleinput}
\mapleinline{active}{1d}{\hspace{.75in}print('holomorphic');}{}
\end{mapleinput}\begin{mapleinput}
\mapleinline{active}{1d}{\hspace{.75in}nf:= newfun(F[j,1],Y);}{}
\end{mapleinput}\begin{mapleinput}
\mapleinline{active}{1d}{\hspace{.75in}Y:= nf[1]; t:= nf[2]; }{}
\end{mapleinput}\begin{mapleinput}
\mapleinline{active}{1d}{\hspace{.75in}print(t);}{}
\end{mapleinput}\begin{mapleinput}
\mapleinline{active}{1d}{\hspace{.5in}else}{}
\end{mapleinput}\begin{mapleinput}
\mapleinline{active}{1d}{\hspace{.75in}F[j,2]:= leadterm(sl,TP);}{}
\end{mapleinput}\begin{mapleinput}
\mapleinline{active}{1d}{\hspace{.75in}F[j,3]:= leadcoeff(sl,TP);}{}
\end{mapleinput}\begin{mapleinput}
\mapleinline{active}{1d}{\hspace{.75in}print('h'=F[j,3]*F[j,2]);}{}
\end{mapleinput}\begin{mapleinput}
\mapleinline{active}{1d}{\hspace{.5in}end if;}{}
\end{mapleinput}\begin{mapleinput}
\mapleinline{active}{1d}{\hspace{.5in}j:= j+1;}{}
\end{mapleinput}\begin{mapleinput}
\mapleinline{active}{1d}{\hspace{.25in}end if;}{}
\end{mapleinput}\begin{mapleinput}
\mapleinline{active}{1d}{   end do;}{}
\end{mapleinput}\begin{mapleinput}
\mapleinline{active}{1d}{   return Y;}{}
\end{mapleinput}\begin{mapleinput}
\mapleinline{active}{1d}{end:}{}
\end{mapleinput}

The procedure "newfun" checks if the global holomorphic function "f" is new, in the sense that it is not a multiple of previous functions.  It returns the (potentially updated) list of truly distinct global holomorphic functions Y, which was fed in through X.
\begin{mapleinput}
\mapleinline{active}{1d}{newfun:= proc (f,X::Array)}{}
\end{mapleinput}\begin{mapleinput}
\mapleinline{active}{1d}{   local Y,yideal,TY,G,GY,A,n,m,t;}{}
\end{mapleinput}\begin{mapleinput}
\mapleinline{active}{1d}{   n:=ArrayNumElems(X)+1;}{}
\end{mapleinput}\begin{mapleinput}
\mapleinline{active}{1d}{   Y:=yarray(n); A:=concat(X,f);}{}
\end{mapleinput}\begin{mapleinput}
\mapleinline{active}{1d}{   TY:=lexdeg([v[2],v[1],b],convert(Y,'list'));}{}
\end{mapleinput}\begin{mapleinput}
\mapleinline{active}{1d}{   yideal:= convert(subtract(Y,A),'list'); }{}
\end{mapleinput}\begin{mapleinput}
\mapleinline{active}{1d}{   G:= gbasis(yideal,TY); GY:=extract(Array(G)); m:=ArrayNumElems(GY); t:=parcials(GY,equal0(Y),y[n]);}{}
\end{mapleinput}\begin{mapleinput}
\mapleinline{active}{1d}{   if t = zero(m) then return [A,'New!',GY,t]; else return [X,'not_new',GY,t]; end if;}{}
\end{mapleinput}\begin{mapleinput}
\mapleinline{active}{1d}{end:}{}
\end{mapleinput}\begin{mapleinput}
\mapleinline{active}{1d}{}{}
\end{mapleinput}

\medskip 

EXECUTION COMMANDS

\medskip 

The procedure "findpoly" searches for global holomorphic polynomials between certain bounds ("m" and "n") on the weighted degree.  It takes as inputthe array "X", which contains previously found functions.
\begin{mapleinput}
\mapleinline{active}{1d}{findpoly:= proc(start,step,n,X)}{}
\end{mapleinput}\begin{mapleinput}
\mapleinline{active}{1d}{   local i, j, Y, L;}{}
\end{mapleinput}\begin{mapleinput}
\mapleinline{active}{1d}{   Y:= X;}{}
\end{mapleinput}\begin{mapleinput}
\mapleinline{active}{1d}{   for i from 0 to n do}{}
\end{mapleinput}\begin{mapleinput}
\mapleinline{active}{1d}{\hspace{.25in}j:= start+step*i;}{}
\end{mapleinput}\begin{mapleinput}
\mapleinline{active}{1d}{\hspace{.25in}print('o');}{}
\end{mapleinput}\begin{mapleinput}
\mapleinline{active}{1d}{\hspace{.25in}print('o');}{}
\end{mapleinput}\begin{mapleinput}
\mapleinline{active}{1d}{\hspace{.25in}print('wdeg' = j);}{}
\end{mapleinput}\begin{mapleinput}
\mapleinline{active}{1d}{\hspace{.25in}print('o');}{}
\end{mapleinput}\begin{mapleinput}
\mapleinline{active}{1d}{\hspace{.25in}L:= List(j);}{}
\end{mapleinput}\begin{mapleinput}
\mapleinline{active}{1d}{\hspace{.25in}Y:= polysearch(L,Y,j);}{}
\end{mapleinput}\begin{mapleinput}
\mapleinline{active}{1d}{\hspace{.25in}print('Y'=Y);}{}
\end{mapleinput}\begin{mapleinput}
\mapleinline{active}{1d}{   end do;}{}
\end{mapleinput}\begin{mapleinput}
\mapleinline{active}{1d}{return Y;}{}
\end{mapleinput}\begin{mapleinput}
\mapleinline{active}{1d}{end:}{}
\end{mapleinput}

\medskip 

\begin{mapleinput}
\mapleinline{active}{1d}{invert:=proc(A::Array, type)}{}
\end{mapleinput}\begin{mapleinput}
\mapleinline{active}{1d}{   local n,TA,yideal,G;}{}
\end{mapleinput}\begin{mapleinput}
\mapleinline{active}{1d}{   n:=ArrayNumElems(A);}{}
\end{mapleinput}\begin{mapleinput}
\mapleinline{active}{1d}{   if type = beta then TA:=lexdeg([v[2],v[1],b],convert(yarray(n),'list')); else      TA:=lexdeg([w[2],w[1],g],convert(yarray(n),'list')); end if;}{}
\end{mapleinput}\begin{mapleinput}
\mapleinline{active}{1d}{   yideal:= convert(subtract(yarray(n),A),'list');}{}
\end{mapleinput}\begin{mapleinput}
\mapleinline{active}{1d}{   G:= gbasis(yideal,TA);}{}
\end{mapleinput}\begin{mapleinput}
\mapleinline{active}{1d}{   return G;}{}
\end{mapleinput}\begin{mapleinput}
\mapleinline{active}{1d}{end:}{}
\end{mapleinput}

\medskip 

SAMPLE COMPUTATION

\medskip 

\begin{mapleinput}
\mapleinline{active}{1d}{A:= findpoly(1,1,11,X);}{}
\end{mapleinput}\begin{mapleinput}
\mapleinline{active}{1d}{nf:=newfun(0,A);}{}
\end{mapleinput}\begin{mapleinput}
\mapleinline{active}{1d}{invert(A,beta);}{}
\end{mapleinput}\begin{mapleinput}
\mapleinline{active}{1d}{B:=invert(arraymodG(A),gamma);}{}
\end{mapleinput}\begin{mapleinput}
\mapleinline{active}{1d}{ArrayNumElems(Array(B));}{}
\end{mapleinput}\begin{mapleinput}
\mapleinline{active}{1d}{solve(B[15]=0,w[2]);}{}
\end{mapleinput}\begin{mapleinput}
\mapleinline{active}{1d}{}{}
\end{mapleinput}
\end{scriptsize}
\end{flushleft}


\bibliography{thesis}

\biography
Carina Curto was born on April 15, 1978, of Argentine parents.  
She was raised in Iowa City, IA.  In 2000 she graduated from Harvard, 
and received an NSF graduate fellowship to support her PhD studies at Duke.  

\end{document}